\documentclass[11pt,a4paper]{article}
\usepackage{times}
\usepackage{dsfont}
\usepackage{fullpage}
\usepackage{titlesec}
\usepackage[normalem]{ulem}
\usepackage{graphicx,cite}
\usepackage{amssymb,amsmath,amsthm}
\usepackage{amssymb,bm}
\usepackage{array}
\usepackage{amsfonts}
\usepackage{enumitem}
\usepackage{color}

\usepackage{ifthen}

\usepackage{setspace}
\usepackage{amsmath,amsfonts}
\usepackage{amssymb}
\usepackage{graphicx}

\usepackage{amsmath,amscd,amssymb,amsgen,amsfonts,amsbsy,amsthm}
\usepackage{mathrsfs}
\usepackage{bm}
\usepackage{bbm}
\usepackage{color}
\usepackage{textcomp}
\usepackage{float}
\usepackage{latexsym,graphicx}
\usepackage{url}
\usepackage{todonotes}
\usepackage{authblk}

\theoremstyle{plain}
\newtheorem{definition}{Definition}

\newtheorem{theorem}{Theorem}
\newtheorem{lemma}[theorem]{Lemma}
\newtheorem{corollary}[theorem]{Corollary}
\newtheorem{proposition}[theorem]{Proposition}

\theoremstyle{remark}

\newtheorem{remark}{Remark}

\newtheorem{example}{Example}

\begin{document}
	\title{On the stability of redundancy models}
	\author[1,3]{E. Anton}
	\author[1,2,3,4]{U. Ayesta}
	\author[5]{M. Jonckheere}
	\author[1,3]{I.M. Verloop}
	\affil[1]{CNRS, IRIT, 2 rue Charles Camichel, 31071 Toulouse, France}
	\affil[2]{IKERBASQUE - Basque Foundation for Science, 48011 Bilbao, Spain}
	\affil[3]{Universit\'e de Toulouse, INP, 31071 Toulouse, France}
	\affil[4]{UPV/EHU, University of the Basque Country, 20018 Donostia, Spain}
	\affil[5]{Instituto de C\'alculo - Conicet, Facultad de Ciencias Exactas y Naturales, Universidad de Buenos Aires (1428) Pabell\'on II, Ciudad Universitaria Buenos Aires, Argentina. }
	\date{}
	\maketitle
	
\begin{abstract}
We investigate the stability condition of redundancy-$d$ multi-server systems. Each server has its own queue and implements popular scheduling disciplines such as First-Come-First-Serve (FCFS), Processor Sharing (PS), and Random Order of Service (ROS).
New jobs arrive according to a Poisson process and copies of each job are sent to $d$ servers chosen uniformly at random. The service times of jobs are assumed to be exponentially distributed. A job departs as soon as one of its copies finishes service. 
Under the assumption that all $d$ copies are i.i.d., we show that for PS and ROS (for FCFS it is already known) sending redundant copies does not reduce the stability region.
Under the assumption that the $d$ copies are identical, we show that (i) ROS does not reduce the stability region, (ii) FCFS reduces the stability region, which can be characterized through an associated saturated system, and (iii) PS severely reduces the stability region, which coincides with the system where all copies have to be \emph{fully} served.
The proofs are based on careful characterizations of scaling limits of the underlying stochastic process. Through simulations we obtain interesting insights on the system's performance for non-exponential service time distributions and heterogeneous server speeds. 
\end{abstract}


%


\section{Introduction}
\label{sec:intro}
The main motivation to investigate redundancy models comes from empirical evidence suggesting that redundancy can help improve the performance of real-world
applications. 
For example Vulimiri et al. \cite{Vulimiri13} illustrate
the advantages of redundancy in a DNS query network where a host computer can query multiple DNS
servers simultaneously to resolve a name. Dean and Barroso \cite{Dean13} note that Google's big
table services use redundancy in order to improve latency. While there are several variants  of a redundancy-based system, the  general notion of
redundancy is to  create multiple copies of the same job that will be sent to a subset 
of servers. By allowing for redundant copies, the aim is to minimize the system latency
by exploiting the variability in  the queue lengths and the capacity  of the different servers. Several
recent works, both empirically (\cite{Ananthanarayanan10,Ananthanarayanan13,Dean13,Vulimiri13})
and theoretically (\cite{Joshi15,Shah16,Gardner16,Lee17a,Lee17b,Gardner17b}), have
provided indications that redundancy can help in reducing the response time of a system.

Most of the literature on performance evaluation of redundancy systems has been carried out under the assumption of i.i.d. copies. Only very recently, a few works that relax this assumption have appeared, see Section~\ref{sec:related} for more details. In particular, Gardner et al. \cite{Gardner17b} recently showed that the i.i.d.\ assumption can be unrealistic, and that it might lead to theoretical results that do not reflect the results of replication schemes in real-life computer systems.

Motivated by the above, in this paper we aim to study the impact that the modeling assumptions have on the performance of the redundancy-$d$ model. 
In particular, we study the dependence of the stability condition on e.g.\ the number of redundant copies, the type of copies (i.i.d.\ copies or identical copies) and the service policy implemented in the servers. 
To some extent, stability is a theoretical notion, since in reality a system will induce stability, for example by limiting the number of accepted jobs. However, stability, or the lack thereof, gives an indication of the quality of the performance that can be expected in practice.

Before detailing our main contributions, we describe a known result that provides a starting point for our work.  
Gardner et al.
\cite{Gardner16,Gardner17} and Bonald and Comte \cite{Bonald17a} show that under the assumption of  exponential service times and  i.i.d. copies, and when the First-Come-First-Serve (FCFS) discipline is implemented in all servers,  the stability region  is not reduced due to adding redundant copies. 
This might seem counter-intuitive at first, as redundancy induces a waste on resources on the $d-1$ servers that work on copies that do not end up being completely finished.  
The reason why the stability is not reduced is due to the assumption of exponential service times and independent copies. Hence, as soon as all servers are busy,
the instantaneous copy departure
rate (and hence job departure rate) is the maximum possible.

\subsection{Main contributions}

We briefly describe the redundancy-$d$ model we consider. There are $K$ servers each with their own queue. The scheduling discipline implemented in all servers is either First-Come-First-Serve (FCFS), Processor Sharing (PS), or Random Order of Service (ROS). New jobs arrive according to a Poisson process at rate~$\lambda$ and $d\leq K$ copies are sent to $d$ servers chosen uniformly at random. The service times of jobs are assumed to be exponentially distributed with parameter $\mu$. A job's service is completed as soon as one of its copies finishes its service. In the absence of redundancy ($d=1$) and for any work-conserving scheduling policy implemented in the servers, the sufficient and necessary condition for stability is $\rho := \frac{\lambda}{\mu K} <1$.
In this paper, we show that adding redundancy can impact the stability condition. An overview of our main results can be found in Table~\ref{tab:freq}. 

\begin{table}[h]\caption{Summary of stability conditions}
	\label{tab:freq}
	\begin{center}
		\begin{tabular}{|c|c|c|c|}
			\hline
			&  \textbf{FCFS} & \textbf{PS}    & \textbf{ROS} \\
			\hline
			\textbf{i.i.d.\ copies} & $\rho<1$   & $\rho<1$ (Prop~\ref{theorem1IID})    &  $\rho<1$ (Prop~\ref{theorem1IID}) \\
			\hline
			\textbf{identical copies} &  $\rho<\bar\ell/K$  (Prop~\ref{prop:FCFS_stab})&$\rho<1/d$ (Prop~\ref{prop:unstab}) & $\rho<1$ (Prop~\ref{prop:ROSID})  \\
			\hline
		\end{tabular}
	\end{center}
\end{table}

In the case of \emph{i.i.d.\ copies}, we prove that with  both PS and ROS, the stability region is not reduced. Hence, even though copies are unnecessarily served, the system remains stable. This statement might lead the reader to think that for any work-conserving policy with i.i.d copies the system is stable under the condition $\rho<1$. This is however not the case, and we present a counterexample based on a priority policy. 

Surprisingly at first sight, in the case of \emph{identical copies}, we prove that the stability condition heavily depends on the scheduling discipline employed by the servers. 

When implementing the Random Order of Service (ROS) discipline, redundancy does not impact the stability condition, that is, it is stable whenever $\rho<1$.

The  stability condition for FCFS with identical copies is given by $\rho<\bar\ell/K$, where $\bar\ell$ denotes the mean number of jobs in service  in an associated saturated system that we characterize. It holds that $\bar\ell <K - (d-1)$, which follows easily by noting that among the jobs being served, the oldest job in the system is  served simultaneously in $d$ servers.
In the particular case of $d=K-1$, the stability region becomes $\rho<2/K$, i.e., it reduces by a factor $2/K$. 
Although in general we cannot obtain closed-form expressions for $\bar\ell$, we prove that $\bar\ell/K$, and hence the stability region, increases as the number of servers, $K$, increases. In the limit, as $K\to\infty$, we numerically observed that $\bar\ell/K$ converges to some $c<1$.
Furthermore, we numerically observe  that  $\bar\ell/K$, and hence the stability region,  decreases in the number of redundant copies, $d$. 

Under PS with identical copies, the system is stable if and only if $\rho <1/d$. In particular, the stability region reduces as the number of redundant copies increases. In fact, under PS the stability condition is the same as in a system in which all $d$ copies have to be fully served, and hence represents the worst possible reduction in the system's stability region. 

Through a light-traffic analysis, we obtain an approximation for the mean number of jobs  under either FCFS, PS or ROS  with identical copies, and find that $d^*=\lfloor K/2\rfloor$ is  the value of $d$ that minimizes the mean number of jobs. This shows that, although the stability region is reduced,  redundancy does improve the performance for sufficiently low loads.

Through simulations, we explore the stability region for general service requirement distributions. Our numerical results indicate the following. 
First,  for i.i.d.\ copies and either FCFS or PS,  the stability region increases in the number of copies $d$, and in the variability of the service requirements. However, with ROS the stability condition seems to be invariant to the service requirement distribution.
Second, if one considers instead  identical copies and either FCFS or ROS,   the performance  deteriorates as the service time variability and/or $d$ increases. Third, for identical copies and PS, the performance  deteriorates as $d$ increases but seems to be nearly insensitive to the service time distribution beyond its mean service time.
Finally, we consider heterogeneous server speeds and present a preliminary analysis and numerics, and observe that for sufficiently heterogeneous servers, the stability region under both FCFS and PS increases in $d$.

In summary, the main takeaway message from our work is that the stability region strongly depends on the modeling assumptions. As shown in the theoretical results, the stability condition depends on the scheduling discipline deployed in servers and  on the correlation structure between copies. Our simulation results illustrate that both the service requirement distribution as well as the service speeds have an important impact on the performance of the system. 
In particular, we believe that our analysis serves as a warning that redundancy needs to be implemented with care in order to prevent an unnecessary degradation of the performance. 

The techniques to prove the results are largely based on sufficiently precise characterizations of scaling limits of the Markov processes describing the number of jobs present in the system, combined with stochastic comparison results.

The rest of the paper is organized as follows. 
First, in Section~\ref{sec:related} we discuss   related work.  
Section~\ref{sec:description} describes the model. Section~\ref{sec:IID} presents the stability results for i.i.d.\ copies.  Sections~\ref{Sec:FCFSdes},~\ref{Sec:Des}~and~\ref{sec:ROS} focus on identical copies  and present the stability results for  FCFS, PS and ROS, respectively. Section~\ref{sec:num} is devoted to insights obtained through simulations and includes the light-traffic analysis. Section~\ref{sec:conclusion}  concludes the paper with some final remarks.
For the sake of readability, the proofs are deferred to the Appendix.

\section{Related work} \label{sec:related}
In redundancy systems with cancel-on-complete  ($c.o.c.$, as considered in this paper), once one of the copies has completed service, the other copies are deleted and the job is said to have received service. 
Most of the recent literature on redundancy has focused on $c.o.c.$ and i.i.d. copies with FCFS as service policy implemented in the servers.  For example, under these assumption, a thorough performance analysis  has been carried out by Gardner et al. 
\cite{Gardner16,Gardner17}, and as mentioned in the introduction,  the stability condition has been fully characterized in \cite{Gardner16,Bonald17a}. 
In Gardner et al. \cite{Gardner16}, the authors consider a class-based
model where redundant copies of an arriving job type are dispatched to a type-specific subset of servers,
and show that the steady-state distribution has a product form. In Gardner et al. \cite{Gardner17}, the
previous result is applied to analyze a multi-server model with homogeneous servers where incoming jobs are
dispatched to randomly selected $d$ servers. 
An important insight obtained there is that stability is not affected by $d$ and that the mean job delay in the system reduces as the redundancy degree $d$ increases.

In a recent study, Gardner et al. ~\cite{GHR19},  the impact of the scheduling policy employed in the server is investigated for i.i.d.\ copies  and exponential service. The authors show that for FCFS the performance might not improve as the number of redundant copies increases, while for other policies  as proposed in that paper,  redundancy does improve the performance.

In Koole and Righter \cite{KR08} the authors show that with FCFS and certain service time distributions (including exponential), the best policy is to replicate as much as possible. Raaijmakers et al.\cite{Raaijmakers2018}, consider FCFS and i.i.d.\ copies, and consider  non-exponential distributed service requirements. As opposed to exponential service requirements, they show that the stability region increases (without bound) in both the number of copies, $d$, and in the parameter that describes the variability of the service requirements. 

Very recently, preliminary results on redundancy without the i.i.d.\ assumption have been published. Gardner et al. \cite{Gardner17b} propose a model in which the service time of a redundant copy is decoupled into two components, one related to the inherent job size of the task, and the other related to the server's slowdown. The paper also proposes a load balancing scheme that in case all servers are busy, it would only dispatch one copy per job.   Such a dispatching policy, under the assumption that the dispatcher has the information regarding the status of servers, would be stable under the condition $\rho <1$. 
Hellemans and van Houdt \cite{HvH18b} consider   identical copies and FCFS, and develop a numerical method  to compute the workload and response time distribution when the number of servers tend to infinity. In order for this method to work, the authors assume the parameters to be such that the system is stable. However, no characterization of the stability region is given in \cite{HvH18b}.

As opposed to $c.o.c.$,  in redundancy systems with cancel-on-start ($c.o.s.$), once one of the copies starts being served, the other copies are deleted.
Up till now, $c.o.s.$\  
has received far less attention than $c.o.c.$. 
The main reason for this comes from the fact that in practice, redundancy aims at exploiting server's speed variability, which is a task that $c.o.c.$ achieves better.  
From the stability point of view, $c.o.s.$ does not bring any extra work to the system, and thus, its stability region is the same as in the non-redundant system. The steady-state distribution of $c.o.s$ has been recently analyzed in Ayesta et al.\cite{ABV18}, and the equivalence of the $c.o.s.$ redundancy model with two other parallel-service models has been shown in Adan et al. ~\cite{IvoAdan2017}. A thorough analysis of $c.o.s.$ in the mean-field regime has been derived in Hellemans and van Houdt \cite{HH18}.


\section{Description of the model}
\label{sec:description}
\begin{figure}[h]
	\begin{center}
		\includegraphics[scale=.5]{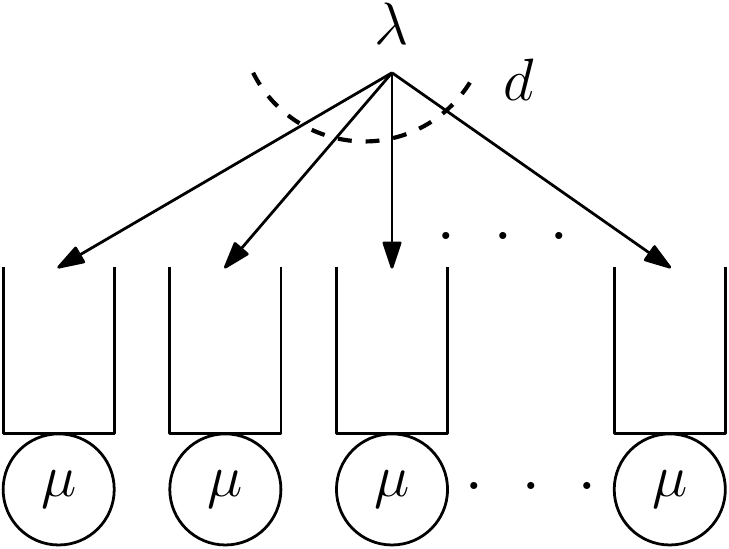}
	\end{center}
	\caption{ Redundancy-$d$ model}
	\label{fig:model}
\end{figure}

As briefly introduced in Section~\ref{sec:intro}, the redundancy-$d$ model consists of $K$ homogeneous servers each with  capacity 1, see Figure~\ref{fig:model}. Jobs arrive according to a Poisson process with rate $\lambda$. An arriving job chooses $d$ servers out of $K$ uniformly at random and sends $d$ copies to these servers. 

We consider two possible correlation structures between copies of the same job: 
\begin{itemize}
	\item The $d$ copies of one job are all i.i.d.\  and have exponentially distributed service requirements with mean $1/\mu$. We refer to this as the redundancy model with i.i.d.\ copies.
	\item The $d$  copies are exact replicas and hence have all the same service requirement, which is exponentially distributed  with mean $1/\mu$. We refer to this as the redundancy model with  \emph{identical copies}.
\end{itemize}
When one of the $d$ copies of a certain job completes its service, the rest of the copies are immediately removed.

We denote by $S$ the set of all servers, $S=\{1,\ldots,K\}$. 
Each job will be assigned a type label, $c=\{s_1,\ldots,s_d\}$, with $s_1,\ldots, s_d\in S, s_i\neq s_j, i\neq j$, to indicate the $d$ servers to which a copy is sent.
We denote by $\mathcal{C}$ the set of all types, that is,
$\mathcal{C}:=\{ \{s_1,\ldots,s_d\} \subset S  \ : \ s_i\neq s_j, \ \quad \forall i\neq j\},$
and $\vert \mathcal{C} \vert = \binom{K}{d}$. 
We denote by $\mathcal{C}(s)$ the subset of types that 
are served at server~$s$, that is,
$\mathcal{C}(s)= \{ c\in \mathcal{C} \ : \ s\in c\}.$
The number of types 
served at server~$s$ equals the number of possible ways to choose  $d-1$ servers out of the remaining $K-1$ servers, that is $\vert \mathcal{C}(s) \vert = \binom{K-1}{d-1}$. 

We denote by $N_c(t)$  the number of type-$c$ jobs at time $t$ and $\vec N(t)=(N_c(t),c\in \mathcal{C})$. 
Furthermore, we denote by $M_s(t) := \sum_{c\in \mathcal{C}(s)} N_c(t)$, $s=1,\ldots,K$, the number of copies per server, and $\vec M(t)=(M_1(t),\ldots,M_K(t))$.
For the $i$-th type-$c$ job, let $b_{cis}$ denote the realization of the service requirement of its copy in server~$s$, $i=1,\ldots, N_c(t)$, $s\in c$. Note that in case the copies are identical, then $b_{cis}=b_{ci}$ for all $s\in c$. We let $a_{cis}(t)$ denote the attained service in server~$s$ of the $i$-th type-$c$ job at time $t$. 
We denote by $ A_c(t)=(a_{cis}(t))_{is}$ a matrix on $\mathbb{R}_+$ of dimension $N_c(t)\times d$. Note that the number of type-$c$ jobs increases by one at rate $\frac{\lambda}{\binom{K}{d}}$, which implies that a row composed of zeros is added to $A_c(t)$. When one element $a_{cis}(t)$ in matrix $A_c(t)$ reaches the required service $b_{cis}$, the corresponding job departs and all of its copies are removed from the system. Hence, row~$i$ in matrix $A_c(t)$ is removed. The rate at which the attained service $a_{cis}(t)$  increases is determined by the employed scheduling policy in that server.

Within a server, a service discipline determines how the capacity of the server is shared among the copies. In this paper, we mostly focus on three service disciplines: \emph{(i)} First-Come-First-Serve (FCFS), where copies within a server are served in order of arrival, \emph{(ii)}  Processor Sharing (PS), where each copy in server~$s$ receives capacity $1/M_s(t)$, and \emph{(iii)} Random Order of Service (ROS), where an idle server chooses uniformly at random a new copy from its queue. All these three policies have in common that they schedule only based on $\{N_c(t), A_c(t),c\in\mathcal{C}\}_{t\geq 0}.$ Hence, the latter is a Markovian descriptor of the system.
As to distinguish between the different policies, we will add a superscript $\{FCFS, PS, ROS\}$ to the process $\vec N(t)$. 

We call the system stable when the process $\vec N(t)$ is positive recurrent, and unstable when the process  $\vec N(t)$ is transient. 
We define the total traffic load  by $\rho:=\frac{\lambda}{\mu K}$. Note that without redundancy, i.e., $d=1$, the system is stable if and only if $\rho<1$ for any work-conserving policy employed in the servers. We further note that for both i.i.d.\ copies and identical copies, the stability region might reduce, but cannot increase. This follows since under exponentially distributed services and homogeneous servers, the total departure rate is at most $K\mu$, while the total arrival rate is $\lambda$. Hence, $\rho<1$ is a necessary stability condition for any value of $d$.  In the remainder of the paper, we will determine the exact stability conditions under the various redundancy models considered. 

\section{Independent identically distributed copies}
\label{sec:IID}
In this section we analyze the stability of the redundancy-$d$ model  when copies of a job are i.i.d..
For FCFS, it was recently proved that $\rho<1$ is the stability condition  with i.i.d.\ copies (\cite{Bonald17a, Gardner16}), that is, the stability condition is not impacted by the redundancy parameter~$d$. 
In Section~\ref{sec:iid_fair}, we prove that the same result holds for PS and ROS.
This result does however not extend to any arbitrary work-conserving policy, as we will show through a counterexample in Section~\ref{sec:priority}. Appendix A contains the proofs of all results obtained in this section. 

\subsection{PS and ROS }
\label{sec:iid_fair}

In this section, we study the policies PS and ROS and prove that their stability condition is $\rho<1$ under the i.i.d.\ copies assumption. 
An intuitive explanation for this result is the following. 
Under both PS and ROS,   on average a fraction $N_c(t)/M_s(t)$ of server~$s$ is dedicated to type-$c$ jobs at time~$t$. Since copies are i.i.d.\, the departure rate of  type-$c$ jobs is given by the sum of the departure rates in 
the $d$ servers (in the set $c$) the job is sent to, that is, $\mu
\left(  \sum_{\tilde s\in c} \frac{N_c(t)}{M_{\tilde s}(t)} \right)$. Now, summing over all jobs types that have a copy  
in server~$s$, we obtain as total departure rate  from server~$s$, 
\begin{equation}
\label{eq:tdr}
\mu
\left( \sum_{c\in \mathcal C(s)}\sum_{\tilde s\in c} \frac{N_c(t)}{M_{\tilde s}(t)} \right).
\end{equation}
For a given time~$t$, let $s_{max}$ be a server containing the largest number of copies, i.e., $M_{s_{max}}(t) \geq M_s(t)$, for all~$s$. It then follows that the departure rate from a server with the largest number of copies equals
$$
\mu
\left( \sum_{c\in \mathcal C(s_{max})}\sum_{\tilde s\in c} \frac{N_c(t)}{M_{\tilde s}(t)} \right)\geq 
\mu
\frac{1}{M_{s_{max}}(t)}\left( \sum_{c\in \mathcal C(s_{max})}\sum_{\tilde s\in c}  N_c(t) \right) 
=
\mu \frac{d}{M_{s_{max}}(t)} \sum_{c\in \mathcal C(s_{max})} N_c(t) =\mu d.
$$ 
The arrival rate  of copies to a server equals $\frac{d}{K}\lambda$.  If $\rho<1$, then $\lambda \frac{d}{K}<\mu d$, hence the total arrival rate to a server with the largest number of copies is smaller than its departure rate, which allows us to prove stability.

In order to make the above exact, we investigate the fluid-scaled system.  The fluid-scaling consists in studying the rescaled sequence of systems indexed by parameter $r$. For $r>0$, denote by $N_{c}^{IID,r}(t)$ the system where the initial state satisfies $N_{c}^{IID}(0)=rn_c(0)$, for all $c\in \mathcal C$. The superscript $IID$ refers to the system under either PS or ROS in the system with i.i.d.\ copies. 
Using standard arguments, see \cite{Bramson06}, we can write for the fluid-scaled number of jobs per type 
\begin{equation}
\label{eq:frelationiid}
\frac{N_c^{IID,r}(rt)}{r}= n_c(0) +   \frac{1}{r}\tilde A_c(rt)  - \frac{1}{r}\tilde S_{c}(T_c^{IID,r}(rt)),
\end{equation}
where $\tilde A_c(t)$ and $\tilde S_{c}(t)$ are independent Poisson processes having rates $ \frac{\lambda}{\binom{K}{d}}$ and $\mu$, respectively,  and $T^{IID,r}_{c}(t)=\sum_{s\in c} T^{IID,r}_{s,c}(t)$, where $T^{IID,r}_{s,c}(t)$ is the cumulative amount of capacity spend on serving type-$c$ jobs in server~$s\in c$ during the time interval $(0,t]$. 

In the following result, we obtain the general characterization 
of a fluid limit.

\begin{lemma}
	\label{lem:sub}
	For almost all sample paths $\omega$ and sequence $r_k$, there exists a subsequence $r_{k_j}$ such that for all $c\in \mathcal C$ and $t\geq0$,
	\begin{equation}\label{eqfluid}
	\lim\limits_{j\rightarrow\infty} \frac{N_c^{IID,r_{k_j}}(r_{k_j}t)}{r_{k_j}}=n^{IID}_c(t) \ u.o.c\footnotemark\ \textrm{ and }
	\end{equation}
	\footnotetext{ where u.o.c. stands for uniformly on compact sets.}
	$$\lim\limits_{j\rightarrow\infty} \frac{T_c^{IID,r_{k_j}}(r_{k_j}t)}{r_{k_j}}=\tau^{IID}_c(t) \ \ u.o.c.,$$
	with $(n^{IID}_c(\cdot),\tau^{IID}_c(\cdot))$ continuous functions.
	In addition,  
	$$ n_c^{IID}(t)=n_c(0) + \frac{\lambda}{\binom{K}{d}}t-\mu\tau_c^{IID}(t),$$
	where $n^{IID}_c(t)\geq0$, $\tau^{IID}_c(0)=0$, $\tau^{IID}_c(t)\leq t$, and $\tau_c^{IID}(\cdot)$ are non-decreasing and Lipschitz continuous functions for all $c\in \mathcal C$.
\end{lemma}

The following lemma gives a partial characterization of the fluid process.

\begin{lemma}\label{lemmaIID} 
	The fluid limit $m_s^{IID}(t):= \sum_{c\in \mathcal C(s)} n_c^{IID}(t)$ satisfies:
	$$
	\frac{\mathrm{d} m^{IID}_s(t)}{\mathrm{d}t} \leq \lambda \frac{d}{K} - \mu d,  \ \ \mbox{if } \ \ m^{IID}_s(t) = \max_{l\in S}\{m^{IID}_{l}(t)\}>0.
	$$
\end{lemma}

In the case $\rho<1$, the drift in the above expression is strictly negative. That is, the maximum of the fluid process $\vec m(t)$ is strictly decreasing with constant rate. Hence, there is a finite time~$T$ when the fluid process is empty. From this, we can directly  conclude that the system is stable, see for instance~\cite{Robert03}.

\begin{proposition}\label{theorem1IID}
	Under either PS or ROS with i.i.d.\ copies, the system is stable  when  $\rho<1$.
\end{proposition}


{
	\begin{remark}[General scheduling policies]
		We believe that the above result holds for any non-preferential scheduling policy that treats all job types equally,  but we did not succeed in obtaining a unifying proof.  Our approach to prove Proposition~\ref{theorem1IID} can be readily extended to cover all policies whose fluid drift is \emph{(i)} continuous and \emph{(ii)} is equal or larger than $\mu d$ for the server(s) with the largest number of copies. Both PS and ROS satisfy this property, but not FCFS. Given the lack of generality of the class of policies that satisfy \emph{(i)}~and~\emph{(ii)}, we chose to restrict the presentation to PS and ROS.
	\end{remark} 
}

\begin{remark}[General service requirement distributions]
	In this paper we focus on exponential distributed service requirements.  The analysis of general service requirement distributions is  a very challenging problem and it will require a different proof technique.   For instance, FCFS with i.i.d.\ copies has been studied in~\cite{Raaijmakers2018} for a specific choice of highly variable service requirements. For an asymptotic regime, the authors show that the stability region \emph{increases without bound} as the service requirement becomes more variable and/or the number of redundant copies increases. This is explained by the fact that each job has $d$ independent copies, and hence, in the (unlikely) event that a copy has a relatively large size, the probability that this copy will be served will become very small as the number of redundant copies increases, or the sizes of the copies become more variable, since the completion of a small-sized copy will directly cancel this large copy.   Therefore, the combination of variable job sizes and redundancy, increases the stability region. 
	In Section~\ref{sec:numIID} we analyze the stability condition of the system under non-exponential service times under PS and ROS service policies. 
\end{remark}

\subsection{Priority policy}
\label{sec:priority}

Given Proposition~\ref{theorem1IID}, one might wonder whether any work-conserving policy would be maximum stable when copies are i.i.d. Indeed, whenever all servers have copies to serve, the total departure rate of jobs equals $K \mu$. This is however not enough to conclude for stability.  In Example~\ref{ex:counter} we give a counterexample.

\begin{example}
	\label{ex:counter}
	
	We consider the system with $K=3$ and $d=2$, hence there are three different types of jobs: $\mathcal C = \{\{1,2\},\{1,3\},\{2,3\} \}$. In server~1, FCFS is implemented. In server~2 and server~3, jobs of types $\{1,2\}$ and $\{1,3\}$ have preemptive priority over jobs of type $\{2,3\}$, respectively. Additionally, within a type, jobs are served in order of arrival.
	
	\begin{figure}[!ht]
		\begin{center}
			
			\includegraphics[scale=0.5]{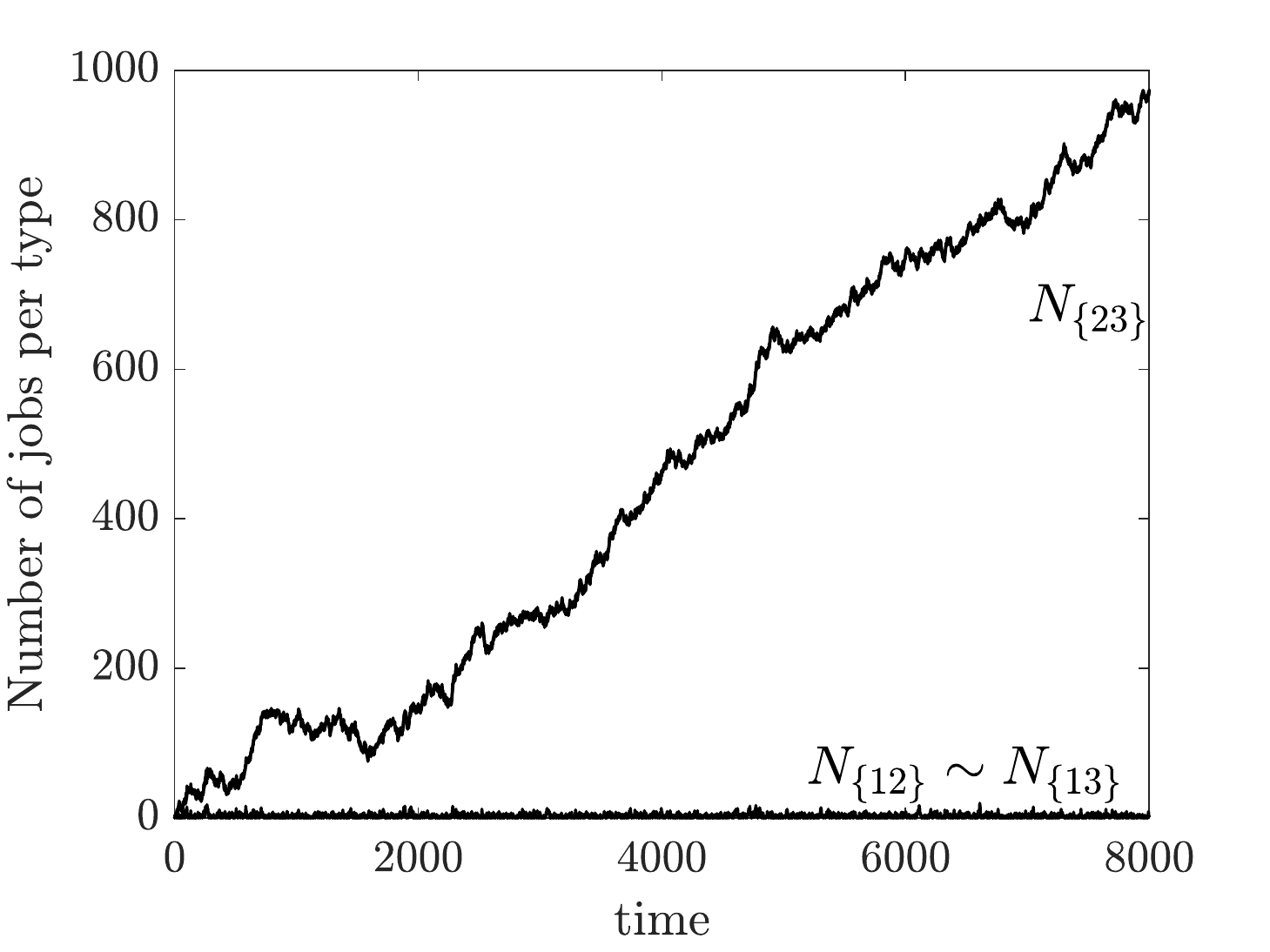}
			\caption{The trajectory of the number of jobs per type with time for the system with $\lambda=2.9$.}
			\label{fig:pri}
		\end{center}
		
	\end{figure}
	In Figure \ref{fig:pri} we have plotted the trajectory of the system when $\rho=0.96 <1$. 
	One observes that   the number of type-$\{2,3\}$  jobs in the system grows large, while  the number of type-$\{1,2\}$ and type-$\{1,3\}$ jobs stay close to 0.
	Hence, the system is clearly unstable, even though $\rho<1$. This is explained as follows: jobs of type-$\{1,2\}$ and jobs of type-$\{1,3\}$ are oblivious to the presence of jobs of type-$\{2,3\}$, due to the preemptive priority assumed in servers 2 and 3. Type-$\{1,2\}$ and type-$\{1,3\}$ jobs have only one server in common. Such a FCFS-redundancy system ($M$-model) has been analysed in~\cite{Gardner16}, where it was obtained that this system (and hence the number of type-$\{1,2\}$ jobs and type-$\{1,3\}$ jobs) is stable when $\rho=\frac{\lambda}{3\mu}< \frac{3}{2}$. 
	
	Type-$\{2,3\}$ jobs are served in server 2 (3) whenever there are no type-$\{1,2\}$ jobs (type-$\{1,3\}$ jobs) present in the system. 
	Note that type-$\{1,2\}$ and type-$\{1,3\}$ jobs behave independent from type-$\{2,3\}$ jobs. Assuming type-$\{1,2\}$ and type-$\{1,3\}$ are in steady state, the  drift of the number of type-$\{2,3\}$ jobs in the system is given by 
	\begin{eqnarray}
	&&\frac{\mathrm{d}}{\mathrm{d}t}{\mathbb E}^{\vec n}\big[N_{\{2,3\}}(t)\big]\Big|_{t=0} \nonumber\\
	&&\quad = \frac{\lambda}{3} - \mu P(N_{\{1,2\}}=0, N_{\{1,3\}}>0) - \mu P(N_{\{1,3\}}=0, N_{\{1,2\}}>0)-2\mu P(N_{\{1,2\}}=0, N_{\{1,3\}}=0).\nonumber
	\end{eqnarray}
	By  \cite{Gardner16}, we have that
	\begin{eqnarray}
	P(N_{\{1,2\}}=0,N_{\{1,3\}}>0)=P(N_{\{1,3\}}=0,N_{\{1,2\}}>0)=P(N_{\{1,3\}}=0,N_{\{1,2\}}=0)\left(\frac{2\mu}{2\mu-\lambda/3}-1 \right),\nonumber
	\end{eqnarray}
	where $P(N_{\{1,3\}}=0,N_{\{1,2\}}=0)=\left(\frac{(2\mu-\lambda/3)^2(3\mu - 2\lambda/3)}{4\mu^2(3\mu-2\lambda/3)+(\lambda/3)^2\mu}\right)$. 
	After some algebra, we have
	\begin{eqnarray}
	&&\frac{\mathrm{d}}{\mathrm{d}t}{\mathbb E}^{\vec n}\big[N_{\{2,3\}}(t)\big]\Big|_{t=0} = \frac{\lambda}{3} - 2\mu \left(\frac{(2\mu-\lambda/3)^2(3\mu - 2\lambda/3)}{4\mu^2(3\mu-2\lambda/3)+(\lambda/3)^2\mu}\right)\left(\frac{2\mu}{2\mu-\lambda/3} \right).\nonumber
	\end{eqnarray}
	It can be checked that the latter is strictly negative if and only if $\rho < 0.91$.
	From this one can conclude that the system is unstable when $\rho>0.91$, using fluid scaling techniques. We however omit the proof as it is out of the scope of this paper.

	
\end{example}

\section{FCFS service policy and identical copies}
\label{Sec:FCFSdes}

In this section we consider the redundancy-$d$ model when copies of a job are identical and when FCFS is employed. We  characterize the necessary and sufficient stability condition and show that the stability condition is reduced when adding redundant copies. This as opposed to the i.i.d.\ case, for which the stability condition remained fixed in~$d$. Appendix B contains the proofs of the results obtained in this section. 

\subsection{Characterization of stability condition}
\label{sec:expFCFS}

Under the FCFS service policy, jobs are served in order of arrival.
If the copies in service in the $K$ servers belong to $k$ different jobs, the departure rate of the system is equal to $k\mu$. The latter is strictly smaller than $K\mu$, even though $K$ servers are busy. This follows from the observation below.

O\scriptsize{BSERVATION}
\label{lem:fcfsd}
\normalsize 1. At every instant of time when the system is not empty, the job that is longest in the system will be running on $d$~servers.

Since  at every instant of time there is a subset of $d$ servers giving service to all the copies of the same job  and copies are identical,  the total output rate of this subset of $d$ servers is  reduced to $\mu$. 
Regarding the $K-d$ remaining servers, the order of arrivals of the jobs impacts the output rate of the remaining servers. As an example, when $K=4$ and $d=2$, the $K-d=2$ remaining servers have as total output rate either $\mu$ (if copies of the same job are first in line in both servers) or $2\mu$. In total, this would give as total output rate either $2\mu$ or $3\mu$.  In both cases, it is strictly less than $K\mu=4\mu$. 

From the above, it is clear that  the total departure rate is not order independent, that is, the total departure rate depends on the order of arrivals of the jobs that are in service. Note that in the case of i.i.d.\ copies, the order-independence property was key to obtain a product-form steady state distribution, see~\cite{Ayesta2019}.  For the case of identical copies (as considered here),  the \emph{lack} of the order-independence assumption prevents us from obtaining a closed-form expression for the stability condition. 
Instead, in the main result of this section we will characterize  the stability condition through the average departure rate in a corresponding system with infinite backlog, referred to as the \emph{saturated system}.
Formally, the saturated system is defined as follows. 
\begin{definition}
	\label{def:congested}
	\textbf{Saturated system:} There is an infinite backlog of  jobs waiting in the system, sampled uniformly over types. There are $K$ servers and the service policy within a server is FCFS. The $d$ copies corresponding to a job are identical.
\end{definition}
We denote the long-run time average number of distinct jobs served in the saturated system  by~$\bar\ell$. Hence, the total departure rate in a saturated system is $\bar\ell \mu$. Below we show that  $\lambda<\bar\ell \mu$, or equivalently, $\rho<\bar\ell/K$,  is a necessary and sufficient condition for the original FCFS system with identical copies to be stable. 
The characterization of the stability condition through a saturation system is reminiscent of the saturation rule obtained in \cite{baccelli95} to prove stability of a large class of queueing systems. We can however not use here their framework due to certain specifics of our model. Instead, in order to prove the stability condition, we resort to stochastic coupling, martingale arguments and fluid limits.  
The proof will be given in Section~\ref{sec:proofFCFS}.
\begin{proposition}\label{prop:FCFS_stab}
	Under FCFS and identical copies, the system is stable if   $\rho< \bar\ell/K$ and unstable if $\rho> \bar\ell/K$ .
\end{proposition}

We will prove in Section~\ref{Sec:Des}, that the stability condition under PS and identical copies is $\rho < 1/d$.  Note that  $\bar\ell\geq \lceil K/d\rceil$, since at least $\lceil K/d\rceil$ jobs are being served at a given time in the saturated system. This gives the following corollary.

\begin{corollary}
	\label{cor:FCFSPS}
	The stability region under FCFS, $\rho< \bar\ell/K$, is larger than under PS, $\rho<1/d$.
\end{corollary}

In the remainder of this section, we characterize~$\bar\ell$. In order to do so, we consider the Markovian state descriptor of the form $ \vec e=(O_{\ell^*},L_{\ell^*-1},\ldots,O_2,L_1,O_1)$. Here, $\ell^*$ denotes the number of jobs that receive service in state~$\vec e$ and $O_j$ denotes the type of the $j$-th job in service. Furthermore, there are $L_j$ jobs that arrived after job $O_j$ and cannot be served since they are waiting for servers that are busy  serving types $O_1,\ldots,O_j$. 
Note that the state descriptor $\vec e$ retains the order of the arriving jobs per type from right to left.

For a given state~$\vec e$, we let  $\ell^*(\vec e)$  denote the number of jobs in service, i.e., $\ell^*$.
Let $\bar E$ denote the state space of the saturated system. 
The mean number of jobs in service can formally be written as
\begin{equation}
\label{eq:averagesFCFS} 
\bar\ell  := \sum_{\vec e\in\bar E} \pi(\vec e)\ell^*(\vec e),
\end{equation}
with
$\pi(\vec e)$  the steady-state distribution of the saturated system.

\begin{figure*}[h]
	\caption{The table and figure show the values of $\bar\ell/K$ for different values of d and $K$.}
	\bigskip
	
	\begin{minipage}[c]{.57\textwidth}
		
		\centering
		\noindent
		\footnotesize
		\setlength{\tabcolsep}{2pt}
		\renewcommand{\arraystretch}{1.2}\begin{tabular}[b]{|c|c|c|c|c|c|c|c|c|}
			\hline
			$\bar\ell/K$	& $K=2$ & $K=3$ & $K=4$ & $K=5$ & $K=6$ & $K=7$ & $K=8$ & $K=9$  \\
			\hline
			$d=1$ & 1	  & 1	  & 1	  & 1     &   1	 & 1 & 1 & 1 \\
			\hline
			$d=2$ & 0.5   & 0.666  &	0.719 & 0.744 &  0.760 & 0.770 & 0.775 &   0.781 \\     
			\hline
			$d=3$ &       & 0.333     &	0.5	  & 0.547  &   0.573 & 0.589 & 0.600 &  0.608  \\
			\hline
			$d=4$ &       &       &	0.25	  &  0.4	  &  0.438 & 0.461 &	0.476&  0.488   \\
			\hline
			$d=5$ &       &       &		  &  0.2	  &   0.333   & 0.364 & 0.384&  0.398    \\
			\hline
			$d=6$ &       &       &		  &   	  &  0.166  & 0.285 & 0.310 & 0.328 \\
			\hline
			$d=7$ &       &       &		  &  	  &    &  0.142  &	 0.250 &  0.270  \\
			\hline
		\end{tabular}
		
	\end{minipage}
	\begin{minipage}[c]{.43\textwidth}
		\centering
		\includegraphics[width=0.75\textwidth]{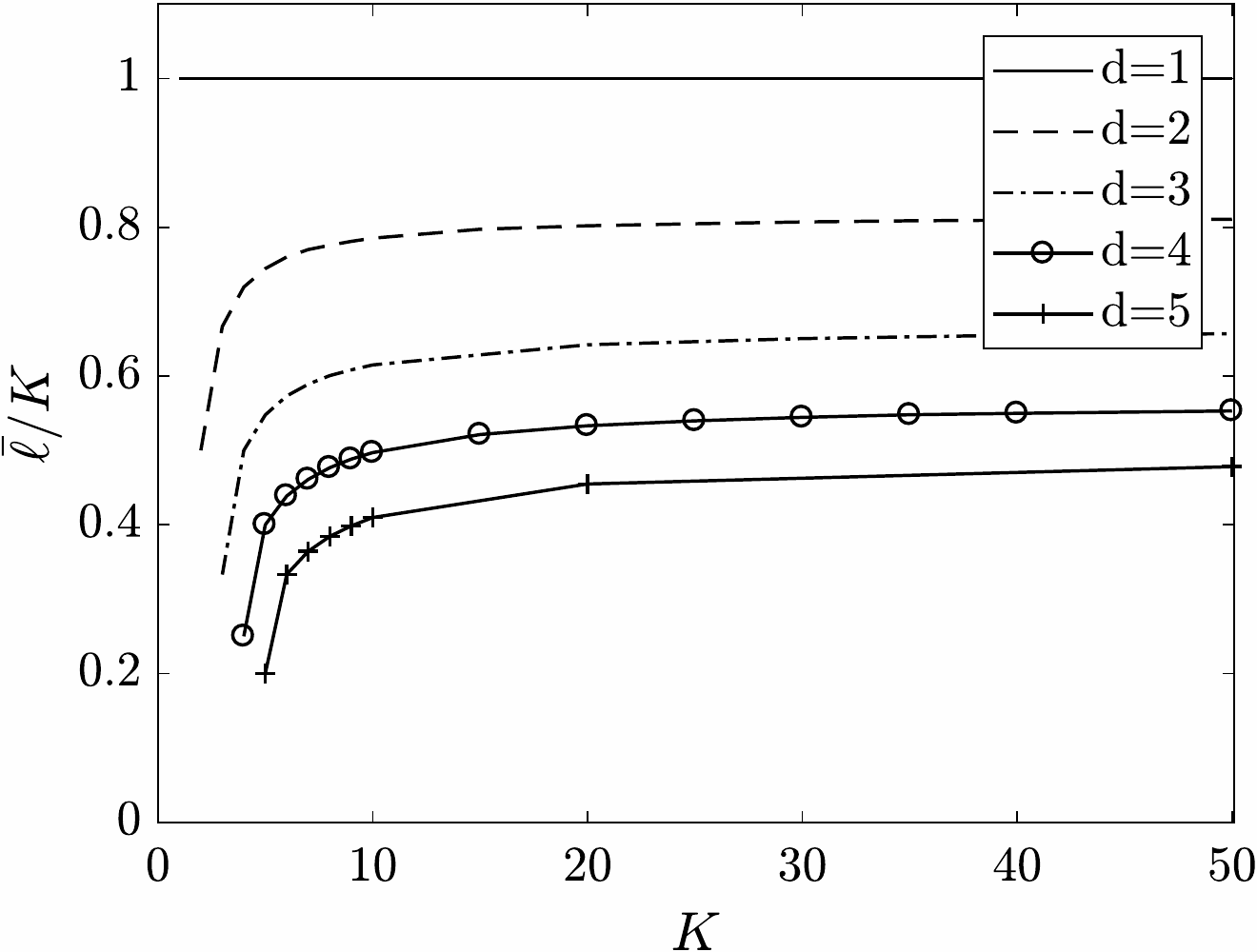}
	\end{minipage}
	\label{tab:fig}
\end{figure*}


In general, no closed-form expression is known for $\bar\ell$. 
In Appendix~B, we write a general expression for the balance equations of the saturated system and state them  explicitly for the case  $d=K-2$ (simplest non-trivial case, since then either two or three jobs are served in the saturated system). From this, we can obtain numerically the value of $\bar\ell$. 
When $d\in \{1, K-1,K\}$,  we can instead get closed-form expressions for $\bar \ell$. When $d=K-1$,  there are $d$ servers that process copies of one job, and the remaining $K-d=1$ server serves one additional job, hence, $\bar\ell=2$. 
When instead $d=1$, there is no redundancy and each server serves one job in the saturated system, i.e., $\bar\ell=K$. When $d=K$, the system behaves as a single server with capacity $\mu$, that is, $\bar\ell=1$.

In Figure~\ref{tab:fig}, we present  $\bar\ell/K$ for different values of $d$ and $K$: the table (\emph{left}) shows $\bar\ell/K$ for small values of $K$ and the figure (\emph{right}) plots the value of $\bar\ell/K$ as $K$ grows large. To  obtain the value of $\bar\ell$ for $d\neq 1,K-2, K-1,K$, we simulated the saturated system, rather than solving the balance equations. 
We  observe from Figure~\ref{tab:fig} that $\bar\ell/K$ (and hence the stability region) increases when the number of servers ($K$) grows large. We make this formal in the proposition below, which is proved using stochastic coupling arguments.

\begin{proposition}
	\label{prop:props}
	For the saturated system, it holds that
	$\bar\ell/K$ is increasing in $K$. 
\end{proposition}

It would be interesting to determine $\lim_{K\to\infty} \bar\ell/K$, as this would represent the stability condition in a mean-field setting. The values in the figure at Figure~\ref{tab:fig} seem to indicate that  $\lim\limits_{K\to\infty}\bar\ell/K=c$ with $c<1$, that is,  the stability region reduces as compared to $d=1$.  We observed that this value $c$ coincides with the value obtained by the numerical method developed in \cite{HH18}.

From  Figure~\ref{tab:fig} we further observe that $\bar\ell/K$ decreases when the number of redundant copies $(d)$ increases. Unfortunately, we did not succeed in finding a coupling argument to prove this property.

\subsection{Proof of stability condition}
\label{sec:proofFCFS}
In this section we show that $\rho<\bar\ell/K$ is both a necessary and sufficient stability condition, that is, we prove Proposition~\ref{prop:FCFS_stab}. 
The dependency on the order of arrivals of the total departure rate   makes exact analysis hard. In order to prove the  stability conditions, we formulate two auxiliary systems that we can compare sample-path wise to the original system. These systems will have the property that for a sufficiently large period of time, a saturated system is observed, and hence, have as average departure rate  $\bar\ell\mu$, which allows us to prove the stability condition.

\subsubsection{Necessary stability condition}

The auxiliary process $\tilde N^{(T)}(t)$ is defined as follows. At time $t=0$, we assume that $\tilde A_c(T)$ type-$c$ jobs arrive, $\forall c\in\mathcal C$.  During the interval $(0,T]$ there are no further arrivals. After time $t>T$, new type-$c$ jobs arrive according to the original Poisson process with rate $\lambda/\binom{K}{d}$.   In the $\tilde N^{(T)}$-system, each server serves according to FCFS. 

To compare the auxiliary process with the original FCFS system, we need to introduce some notation. 
The attained service of the copy of the $i$-th type-$c$ job in server~$s$, $a^{FCFS}_{cis}(t)$, will be compared to the attained service of the same copy in the $\tilde N^{(T)}$-system. For that, (with slight abuse of notation), we let $a^{\tilde N^{(T)}}_{cis}(t)$ denote the attained service of \emph{this same} copy, where we assume that in case this copy has already departed in the $\tilde N^{(T)}$-system, then $a^{\tilde N^{(T)}}_{cis}(t)$ is set equal to its service requirement $b_{ci}$. In the result below we show that sample-path wise, a job departs earlier in the $\tilde N^{(T)}$ system than in the original system. 
In particular, this implies that if the original FCFS model is stable, then the $\tilde N^{(T)}$-system is stable as well. 

\begin{lemma}\label{prop:FCFS_lb}
	Assume $N_c^{FCFS}(0) =\tilde N^{(T)}_c(0)$ and $a_{cis}^{FCFS}(0)= a^{\tilde N^{(T)}}_{cis}(0)$, for all $c,i, s$. Then, $\tilde N^{(T)}_c(t) \leq N^{FCFS}_c(t) + (\tilde A_c(T)-\tilde A_c(t))^+$ and $a_{cis}^{FCFS}(t)\leq a^{\tilde N^{(T)}}_{cis}(t)$, for all $i=1,\ldots,N^{FCFS}_c(t)$, $c\in \mathcal C$, $s\in S$.	
\end{lemma}

Let the random variable $\tau(T)>0$ denote the  moment that one of the servers becomes empty. In the time interval $[0,\tau(T)]$, the $\tilde N^{(T)}$-system will behave as a saturated system. We will prove that as $T$ grows large, $\tau(T)$ grows large, and due to the law of large numbers, the time-average number of jobs in service in  the interval $[0,\tau(T)]$ will be equal to $\bar\ell$,  as defined in~\eqref{eq:averagesFCFS}.  Since each job in service has a departure rate $\mu$, this allows us to prove that if the $\tilde N^{(T)}$-system is stable, then  $\lambda<\bar\ell \mu$. Together with Lemma~\ref{prop:FCFS_lb} this gives the following result. 

\begin{proposition}\label{prop:FCFS_n}
	Under FCFS and identical copies, the system is unstable if $\rho> \bar\ell/K$ .
\end{proposition}

\subsubsection{Sufficient stability condition} 
In order to prove that $\rho<\bar\ell/K$ is a sufficient stability condition, we define the process $\hat N(t)$ as follows. 
In the time interval $[0,|\hat N(0)|/\mu]$, only those jobs that were already present at time 0 are allowed to be served (according to FCFS). From time $t, t\geq |\hat N(0)|/\mu$  onwards, all jobs present in the system can be served. 

We first establish a sample-path comparison with the original FCFS system, which allows us to conclude for stability of the original process.   
We let $a_{cis}^{\hat N}(t)$ denote the attained service of the $i$-th type-$c$ job in the $\hat N$-system. The attained service of the $i$-th type-$c$ job in server~$s$ in the $\hat N$-system will be compared to the attained service of the same copy in the FCFS system. In order to do so, with a slight abuse of notation, we let $a_{cis}^{FCFS}(t)$ denote the attained service of \emph{this same} copy, where we assume that in case this copy has already departed in the FCFS-system, then it is set equal to its service requirement $b_{ci}$. 

\begin{lemma}\label{prop:FCFS_up}
	Assume $N_c^{FCFS}(0) =\hat N(0)$ and $a_{cis}^{FCFS}(0)= a^{\hat N }_{cis}(0)$, for all $c,i, s$. Then, $\hat N_c(t)\geq N^{FCFS}_c(t)$ and 
	$a_{cis}^{\hat N}(t)\leq a^{FCFS}_{cis}(t)$, for all $i=1,\ldots,\hat N_c(t), c\in \mathcal C,  s\in S$.
\end{lemma}

For the stochastic process~$\hat N(\cdot)$, we will see that the system is stable if $\rho<\bar\ell/K$. 
To do so, we will characterize the fluid limit. We will show that at the moment the auxiliary process can start serving jobs that were not present at time~0, the queue has built up, and during a considerable amount of time the system will behave as a saturated system. Hence,  the average number of occupied servers equals $\bar\ell$, which allows us to prove that $\hat N(t)$ is stable if $\rho<\bar\ell/K$. Together with Lemma~\ref{prop:FCFS_up} this gives the following result.  

\begin{proposition}\label{theo:FCFS_s}
	Under FCFS and identical copies, the system is stable if $\rho <\bar\ell/K$.
\end{proposition}


\section{PS service policy and identical copies}
\label{Sec:Des}
In this section we investigate the redundancy-$d$ model with identical copies when PS is employed in all the servers.
We will show that the system is stable if and only if $\rho<1/d$.
We note that this coincides with the stability condition of a system where all $d$ copies have to be fully served.

\begin{proposition}
	\label{prop:unstab}
	Under PS and identical copies, the system is stable if $\rho<\frac{1}{d}$ and unstable if $\rho>\frac{1}{d}$. 
\end{proposition}

Before proceeding to the intuition (Section~\ref{sec:phenomena})  and proof of Proposition~\ref{prop:unstab} (Section~\ref{sec:proof}), we first introduce a new notation. 
Under PS, the attained service of the copy of the $i$-th type-$c$ job in server~$s$ increases at speed $1/M^{PS}_{s}(t)$, that is, $\frac{\mathrm{d} a^{PS}_{cis}(t) }{\mathrm{d}t}= \frac{1}{M^{PS}_{s}(t)}$, $\forall c\in \mathcal{C}(s)$, $i=1,\ldots, N_c(t)$.  Note that a departure of a job is due to a departure in the server where it has the largest attained service. Denote by $s_{ci}^*(t)$ the server that contains the copy of the $i$-th type-$c$ job with the largest attained service, that is,  $s_{ci}^*(t):=argmax_{s\in c}\{a^{PS}_{cis}(t)\},$ for all $c\in \mathcal{C}$, $i=1,\ldots, N_c(t)$. The instantaneous departure rate of the $i$-th type-$c$ job under PS is hence 
$\frac{\mu}{M^{PS}_{s_{ci}^*(t)}(t)}.$
In particular, the number of type-$c$ jobs decreases at rate 
\begin{equation}
\label{eq:Msci}
\sum_{i=1}^{N^{PS}_c(t)} \ \frac{\mu}{M^{PS}_{s_{ci}^*(t)}(t)}.
\end{equation}

\subsection{Intuition behind stability condition and its proof}
\label{sec:phenomena}

\begin{figure*}[thb]
	\centering
	\begin{minipage}[b]{.45\textwidth}
		\includegraphics[width=1.1\textwidth]{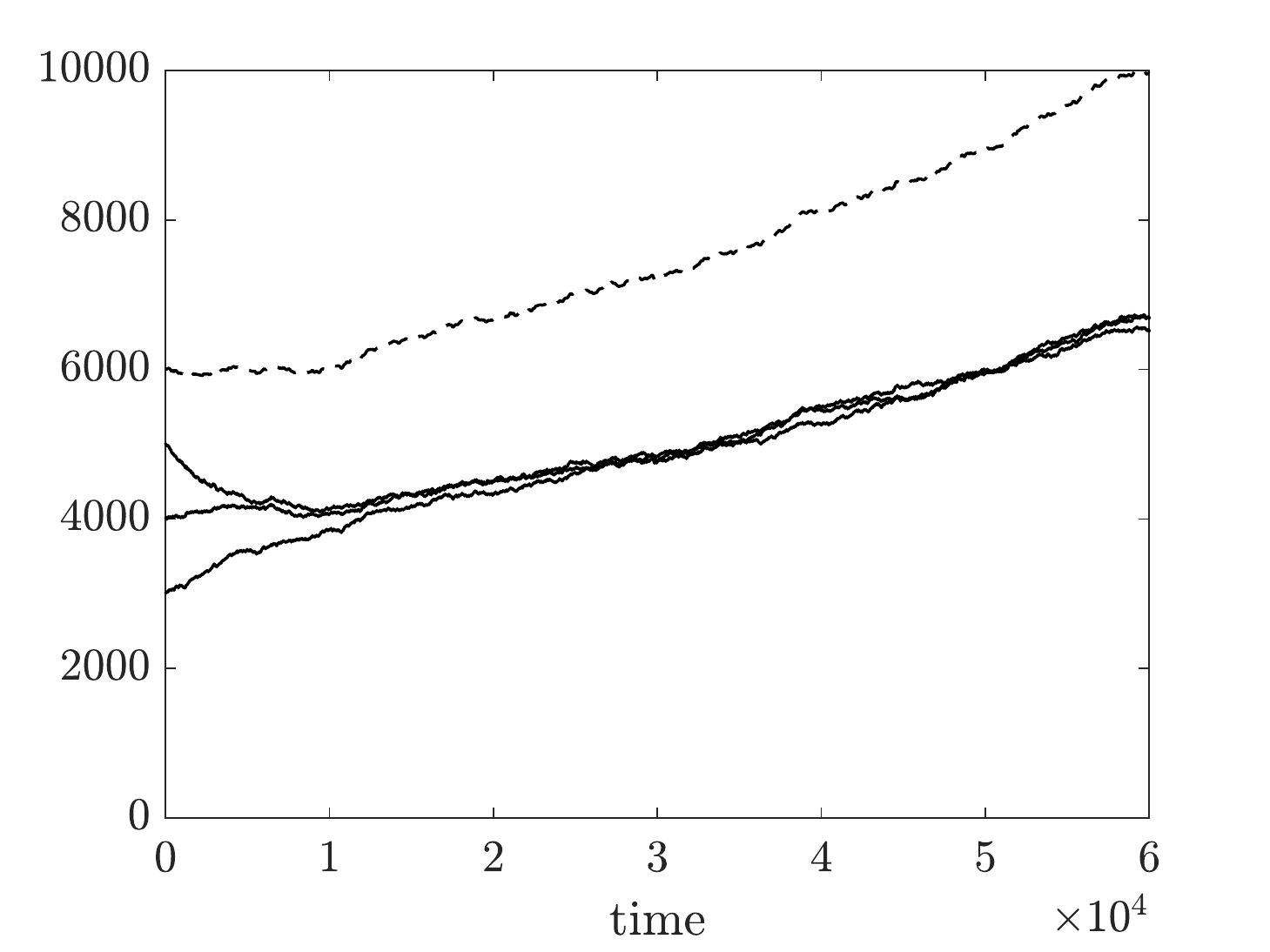}
	\end{minipage}
	\begin{minipage}[b]{.45\textwidth}
		\includegraphics[width=1\textwidth]{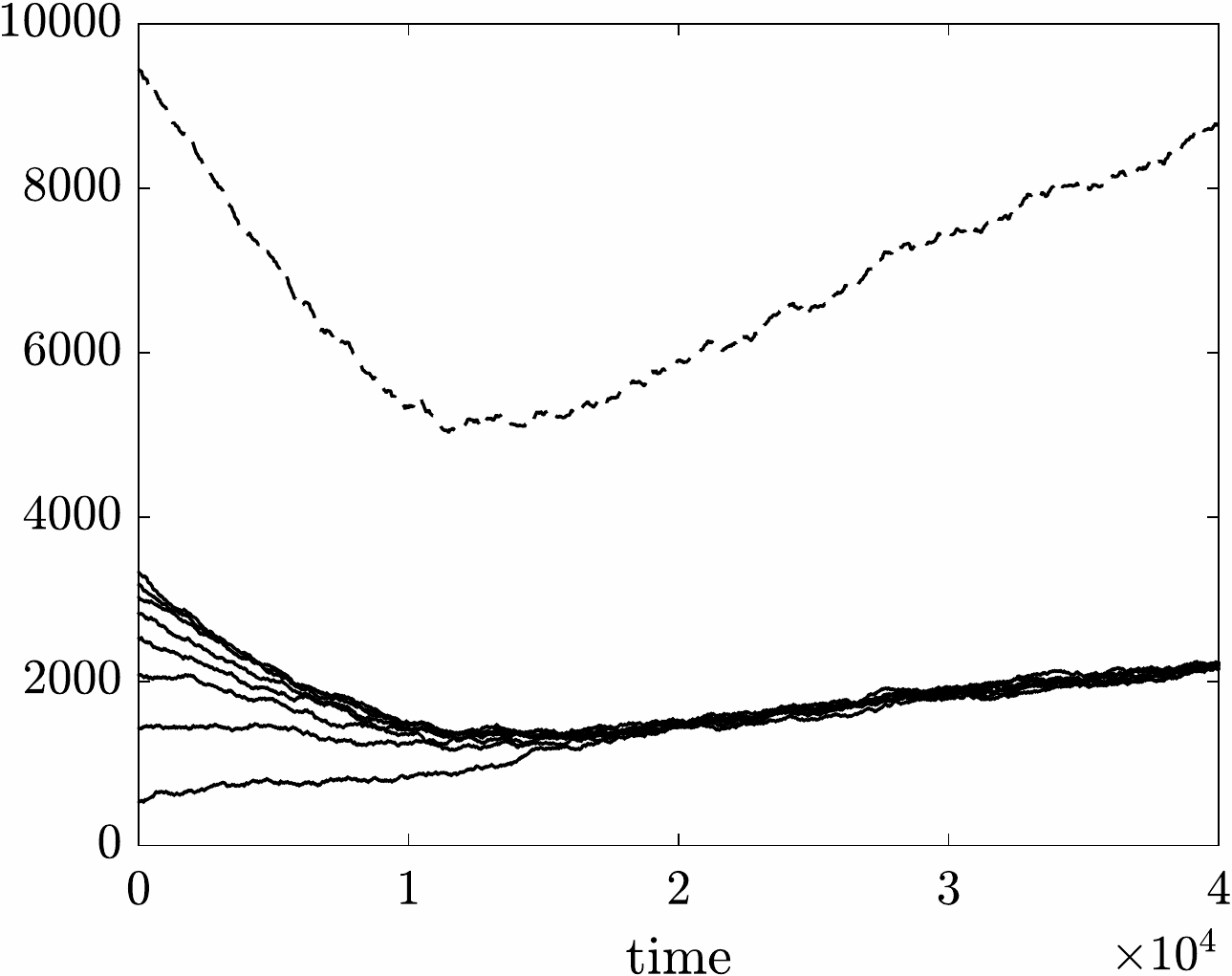}
	\end{minipage}
	\caption{The dashed line represents the total number of jobs in the system under PS with identical copies.  The other lines represent the number of copies in each of the servers. (left) $K=3$, $d=2$ and $\rho=0.53$, (right) $K=8$, $d=2$ and $\rho=0.52$.}
	\label{fig:int}
\end{figure*}

To illustrate why $\rho<1/d$ is the stability condition, we have plotted in Figure~\ref{fig:int} the trajectories of the number of copies in each of the servers,~$M^{PS}_s(t)$, for two settings, $K=3$ and $K=8$, with $d=2$. In both cases, we assume the load is such that $\rho>1/d$. We let the processes  start in a very large state, and plot the trajectories over a large time horizon.  

In Figure~\ref{fig:int}, we observe the following effect. When the processes $M^{PS}_s(t)$ are unbalanced (as is the case for $t<10^4$), the numbers of copies at the most loaded servers decrease. 
Consider one of the highly-loaded servers, referred to as server~$\tilde s$. Now, a copy in  service in server~$\tilde s$ will leave because either it has obtained full service in server~$\tilde s$, or a copy of the same job finished service in another \emph{less-loaded} server. 
The rate at which a copy is served at such a less-loaded server, is higher than that in the high-loaded server~$\tilde s$. 
Thus, the effective departure rate of copies from server~$\tilde s$ will be higher than $\mu$.
Since the arrival rate of new copies to a given server equals $\lambda \frac{d}{K}$, this explains why the number of copies in  server~$\tilde s$ (a higher loaded server) can go down, even though $\lambda d/K >\mu$ ($\rho>1/d$). 

Hence, during a certain time, the system experiences a "good phase" in which higher-loaded servers decrease and  the total queue length decreases as well. However, once the servers are more equally loaded, we observe that the total queue length starts to  build up. To explain this, consider the symmetric case, i.e., $M^{PS}_s(t)=m$, for all $s$. Then, each copy of a job receives in each server the same fraction of capacity.  Hence, the departure rate of copies from a server is $\mu$ (see Eq.~\eqref{eq:Msci}). Since $\lambda d/K>\mu$, the servers will build up from then on, and the total number of jobs will  diverge.
 
In order to prove the stability condition, the challenge is to prove instability. We note that the total number of jobs cannot be taken as Lyapunov function:  As we described above, inside some cone around the diagonal (symmetric states), the drift of the total number of copies in the system is strictly positive, while outside that cone,  the drift of the total number of jobs is decreasing. 
We further observe from   Figure~\ref{fig:int} that the drift of the server with the minimum number of copies is strictly positive, while the  drifts of the higher-loaded servers is first negative, until they join the minimum, from which point on they stay together and increase. 
This motivated us to study  the drift of the server with the minimum number of copies. Though it is a complicated (non-monotone) function for the stochastic process, one can show that for the fluid limit, the drift of the server with the minimum number of copies is strictly positive. 
So, even if at a short time scale, the minimum cannot be taken as Lyapunov function, the minimum at a fluid scale does go up if $\rho > 1/d$. This is exactly what is used  in order to prove unstability, see Lemma~\ref{lemma2}. 


	\subsection{Proof of stability condition}
\label{sec:proof}
Having identical copies makes exact analysis hard, as it requires to keep track of the attained service of the copies in each of the servers. 
In order to derive the necessary and sufficient stability condition, i.e. to prove Proposition~\ref{prop:unstab}, we describe two systems that lower and upper  bound the original PS system. These systems will have the property that the departure rate no longer depends on the attained service, which allows us to prove  necessary (sufficient) conditions for stability for the lower bound (upper bound), and hence also for the original system. The full proofs can be found in Appendix C.

\subsubsection{Necessary stability condition}
For the original system, the departure rate	of the number of type-$c$ jobs  depends on the attained service, see Equation~\eqref{eq:Msci}. More precisely, the departure rate of the $i$-th type-$c$ job equals 
\begin{equation}
\label{eq:Mpsf}
\frac{\mu}{M^{PS}_{s_{ci}^*(t)}(t)},
\end{equation}
where we recall that $s_{ci}^*(t)$ denotes  the server where a copy of this job has received most service so far. 
The lower-bound system is defined as follows: We replace~\eqref{eq:Mpsf} by 
$$
\frac{\mu}{M^{PS}_{s_{c}^{min}(\vec N(t))}(t)},
$$
where $s^{min}_{c}(\vec N(t)) := \arg\min_{s\in c}\{M_s(t)\}$ is the server with the least number of copies that contains a type-$c$ job at time $t$ (ties are broken at random). 
That is, in the lower-bound system, a type-$c$  job receives service from the 
server in the set $c$ with the minimum number of copies.
We note that since the lower-bound system does no longer depend on the attained service, it is more amenable to get the stability condition.

The lower-bound system is  described by $\{N_c^{LB}(t),c\in\mathcal{C}\}_{t\geq 0}$, living on the countable set $\mathbb{Z}_+^{\binom{K}{d}}$.  Here, $N_c^{LB}(t)$ denotes the number of type-$c$ jobs in the lower-bound system. 
The process $N_c^{LB}(t)$ increases by one at rate $\lambda/\binom{K}{d}$ (as is the case for the original process), and decreases by one at rate 
\begin{equation}\label{eq:depLB}
\mu\frac{N^{LB}_c(t)}{M^{LB}_{s^{min}_{c}(\vec N^{LB}(t))}(t)},
\end{equation}
where $M^{LB}_s=\sum_{c\in \mathcal C(s)} N^{LB}_c$.
Note that Equation (\ref{eq:depLB}) coincides with Equation~\eqref{eq:Msci}, where now $s^*_{ci}(t)$ is replaced by $s_c^{min}(\vec N(t))$ (because  for a given type, all jobs share the same server with the smallest number of copies). Below, we prove that this system gives a stochastic lower bound for the original system. 

\begin{lemma}\label{prop1} 
	Assume $N_c^{PS}(0)=N^{LB}_c(0)$, for all $c$. Then, 
	$N^{PS}_c(t)\geq_{\textrm{st}}N_c^{LB}(t)$, for all $c\in \mathcal C$ and $t\geq0$.
\end{lemma} 

In Lemma~\ref{lemma1} below, we  give an expression for  the departure rate from a server~$s$ in the lower-bound system. Before doing so, we need to introduce some notation.
For each server $s$, we define $D_s(\vec N^{LB}(t)):=\{l\in S  : M^{LB}_s(t)\geq M^{LB}_l(t)\}$.  
We denote by $\mathcal C_l^s(\vec N(t)):=  \{c\in \mathcal C(s) : l=s_c^{min}(\vec N(t))  \},$
the subset of types that are served in server~$s$ and for whom server~$l$ is the server with the minimum number of copies that serve type~$c$. 
Notice that $\mathcal C(s)$ is the disjoint union of the above elements, $\mathcal C(s)=\cup_{l\in D_s(\vec N(t))}\mathcal C_l^s(\vec N(t))$.	

\begin{lemma}\label{lemma1}
	For the lower-bound system, when in state $\vec N^{LB}(t)=\vec n^{LB}$, the number of copies in server~$s$, $M^{LB}_{s}(t)$, decreases by one at rate
	\begin{equation}
	\label{eq:lemma1}
	\mu\left( 1 +\sum_{l\in D_s(\vec n^{LB})} \frac{(M^{LB}_{s}(t)-M^{LB}_{l}(t))\sum_{c\in \mathcal C_l^s(\vec n)}N^{LB}_c(t)}{M^{LB}_{s}(t)M^{LB}_{l}(t)}\right).\end{equation}
\end{lemma}

In particular, 
from Equation (\ref{eq:lemma1}) we clearly see the improvement brought by redundancy. For the server with the minimum number of copies, Equation~\eqref{eq:lemma1} simplifies to $\mu$. This server is hence not receiving any help from the other (higher-loaded) servers.
However, servers that do not have the minimum number of copies, do benefit from redundant copies, as their service rate is $\mu$ plus some additional positive fractions.
This is due to the fact that in the lower-bound system, all types in server~$s$ that also have a copy in another server with less copies, will receive as effective service rate that what they would get in this latter server, and hence receive a higher capacity than what they would get in server~$s$. 


We study the fluid limit of the lower bound system in order to conclude the lower-bound system is transient when $\rho>1/d$. The fluid-scaling consists in studying the rescaled sequence of systems indexed by parameter $r$. For $r>0$, denote by $N_{c}^{LB,r}(t)$ the system where the initial state satisfies $N_{c}^{LB}(0)=rn_c(0)$, for all $c\in \mathcal C$. The associated number of copies per server is given by $M_s^{LB,r}(t)=\sum_{c\in \mathcal{C}(s)}N_c^{LB,r}(t)$, for all $s\in S$. For the fluid-scaled number of jobs per type we can write 
\begin{equation}
\label{eq:frelation}
\frac{N_c^{LB,r}(rt)}{r}= n_c(0) +   \frac{1}{r}\tilde A_c(rt)  - \frac{1}{r}\tilde S_{c}(T_c^{LB,r}(rt)),
\end{equation}
where $T^{LB,r}_c(t)$ is defined as the cumulative amount of capacity spent on serving type-$c$ jobs in server $s^{min}_c(\vec N^{LB,r}(\cdot))$ during the time interval $(0,t]$. 
The existence of   fluid limits can be proved  as before: The statement of Lemma~\ref{lem:sub} indeed directly translates to the process $\vec N^{LB,r}(t)$, and is therefore left out. In the following result, we obtain the general characterization of a fluid limit.

\begin{lemma}\label{lemma2} 
	The fluid limit $m_s^{LB}(t):= \sum_{c\in\mathcal C (s)} n_c^{LB}(t)$ satisfies:
	$$
	\frac{\mathrm{d} m^{LB}_s(t)}{\mathrm{d}t} = \lambda \frac{d}{K} - \mu,  \ \ \mbox{if } \ \ m^{LB}_s(t) = \min_{l\in S}\{m^{LB}_l(t)\}>0,
	$$
	and
	$$
	\frac{\mathrm{d} m^{LB}_s(t)}{\mathrm{d}t} \geq \lambda \frac{d}{K} - \mu,   \ \ \mbox{if } \ \ m^{LB}_s(t) = \min_{l\in S}\{m^{LB}_l(t)\}=0.
	$$
\end{lemma}

In case $\lambda d/K-\mu>0$, this partial characterization of the fluid limit implies the following. Consider servers whose amount of fluid   is the minimum, that is, consider servers belonging to the set $U(t):=\{s\in S: m^{LB}_s(t)\leq m^{LB}_{\tilde s}(t), \forall \tilde s\}$. By Lemma~\ref{lemma2}, the amount of fluid in these servers increases with  a strictly positive rate $\lambda d/K-\mu$.
Moreover, if at time $t_0>t$, some server~$\tilde s$ is added to this set, that is, $U(t_0)=U(t)\cup \{\tilde s\}$,  this server will increase as well from that moment on with the same rate  $\lambda d/K-\mu$.

This uniform divergence of the fluid limit, together with bounds on the macroscopic drifts, allows us to show instability of the stochastic process $\vec N^{LB}(t)$  via a usual transience criterion for Markov chains  whenever the fluid drift $\lambda d/K-\mu$ is strictly positive. Together with Lemma~\ref{prop1}, this allows us to prove the following result.

\begin{proposition}\label{theorem1}
	Under PS and identical copies, the system is unstable if  $\rho>1/d$.
\end{proposition}

\subsubsection{Sufficient stability condition}

For the original system, a job departs the system once a copy has received its service in one of the servers. We will now  upper bound this, by considering the same   system, but  where a job departs from the system only if \emph{all its copies} have completed service. 

For the UB-system, we let $N_c^{UB}(t)$ denote the number of type-$c$ jobs, and $A_c^{UB}(t)=(a^{UB}_{cis}(t))_{is}$, with $a^{UB}_{cis}(t)$ the attained service of the $i$-th type-$c$ job in server~$s$.  
Note that the number of copies in server~$s$ is given by $M^{UB}_s(t)=\sum_{c\in \mathcal C(s)}\sum_{i=1}^{N^{UB}_c(t)}\mathbf{1}_{(a^{UB}_{cis}(t)<b_{ci})}$, since a copy is only present in server~$s$ when $a^{UB}_{cis}(t)$ is strictly smaller than the service requirement $b_{ci}$.	
The $i$-th type-$c$ job in server $s$ is served at speed $1/M^{UB}_{s} (t)$, hence $\frac{\mathrm{d}a^{UB}_{cis}(t)}{\mathrm{d}t} = 1/M^{UB}_{s} (t)$. Now, the $i$-th type-$c$ job departs from the system once $a^{UB}_{cis}(t)=b_{ci}$ for all servers $s\in c$, that is,  when all the copies of a job are fully served. 

To compare the UB-system with the original PS system, we need to compare the attained service of the $i$-th arrived job in both systems. For that, we denote by $\alpha_{i,s}^{UB}(t)$ and $\alpha_{i,s}^{PS}(t)$ the attained service of the $i$-th arrived job in server~$s$, for UB and PS, respectively. With slight abuse of notation, we set $\alpha_{i,s}^{PS}(t)$ equal to $\beta_{i}$ (the service requirement of the $i$-th arrived job) for all servers~$s$, in case it has departed from the PS system.  

\begin{lemma}\label{lemmaup}
	Assume  $\alpha_{i,s}^{PS}(0)=\alpha_{i,s}^{UB}(0)$, for all $i=1,\ldots,$ and $s\in \mathcal{S}$. Then, 
	$\alpha^{UB}_{i,s}(t)\leq_{st} \alpha^{PS}_{i,s}(t)$	 for all $t\geq0$, $i=1,\ldots$, and   $s\in c$.
	In particular, $N_c^{UB}(t)\geq_{st}N^{PS}_c(t)$.
\end{lemma}

In the upper-bound system, all copies need to be served until a job departs. Hence, each queue receives copies at rate $\lambda d/K$ and copies depart at rate $\mu$, that is, its marginal distribution is that of an $M/M/1$ system with arrival rate $\lambda d/K$ and departure rate $\mu$. The latter is positive recurrent if and only if $\rho<1/d$.
Since $\vec N^{UB}(t)$ serves as an upper bound for our model (Lemma~\ref{lemmaup}), this implies that the original system is positive recurrent as well, as stated in the result below. 

\begin{proposition}~\label{lemma:ps1}
	Under PS and identical copies, the system is stable if $\rho<1/d$. 
\end{proposition}

\section{ROS service policy and identical copies}

\label{sec:ROS}

In this section we study ROS with identical copies and show that the stability condition is $\rho<1$. Appendix D contains the proofs of the results obtained in this section.

\subsection{Intuition behind stability condition and its proof}

Under ROS with identical copies, an idle server chooses uniformly at random a new copy from its queue and serves it until the copy finishes service, or one of its identical copies finishes service in another server.
Note that if $k$ servers are serving different jobs,  then the total departure rate of these $k$ servers is $\mu k$. If however these $k$ servers are serving a copy from the same  job, then these $k$ servers give together a total departure rate $\mu$ (since copies are identical),  hence capacity is wasted. 

From the above, we observe  that 
$
\mathbb{P}(\mbox{every copy in service belongs to a unique job})$
is an important measure to determine the stability condition under ROS. 
Note that this probability is strictly smaller than 1  when the queue length is small, hence capacity is wasted. However, as the queues grow large, this probability will converge to 1, showing that under the fluid scaling no capacity is wasted. This then allows us to conclude that the stability condition is not reduced when adding redundant copies, that is, $\rho<1$ is the stability condition.

\subsection{Proof of stability condition for ROS}
In order to prove the stability, we investigate the   fluid-scaled system. For $r>0$, denote by $N^{ROS,r}_c(t)$ the system where the initial state satisfies $\vec N^{ROS,r}(0)=r \vec n(0)$. Using routing arguments, we can write 
\begin{equation}
\label{eq:relation}
\frac{N^{ROS,r}_c(rt)}{r}= n_c(0) +  \frac{1}{r} \tilde A_c(rt)  -\frac{1}{r} \tilde S_{c}(T^{ROS,r}_c(rt)),
\end{equation}
where $T^{ROS,r}_c(t)$ is defined as  the cumulative amount of capacity spent on serving \emph{a first copy} of type-$c$ jobs in the interval~$(0,t]$. For a given job, we refer with ``first copy'' to that copy (out of the $d$) that was first to enter into service.  

The existence of the fluid limit can be proved. In fact, the statement of Lemma~\ref{lem:sub} and its proof directly carries over, and is therefore left out. 
The following lemma gives a partial characterization of the fluid process. For the proof see Appendix D.

\begin{lemma}\label{lemmaROS} 
	The fluid limit $m^{ROS}_s(t):= \sum_{c\in \mathcal C(s)} n^{ROS}_c(t)$ satisfies the following:
	$$
	\frac{\mathrm{d} m^{ROS}_s(t)}{\mathrm{d}t} \leq \lambda \frac{d}{K} - \mu d,  \ \ \mbox{if } \ \ m^{ROS}_s(t)=\max_{l\in S}\{m^{ROS}_l(t)\} > 0.
	$$
\end{lemma}

In case $\rho<1$, the drift in the above expression is strictly negative. That is, the maximum of the fluid process $\vec m(t)$ is strictly decreasing with constant rate. Hence, there is a finite time~$T$ when the fluid process is empty. From this, we can directly  conclude stability (same steps as in the proof of Proposition~\ref{theorem1IID}). 

\begin{proposition}
	\label{prop:ROSID}
	Under ROS with identical copies, the process  $\vec N^{ROS}(t)$ is ergodic  when $\rho<1$.
\end{proposition}
\begin{figure*}[t]
	\centering
	\begin{minipage}[b]{.45\textwidth}
		
		\includegraphics[width=1.\textwidth]{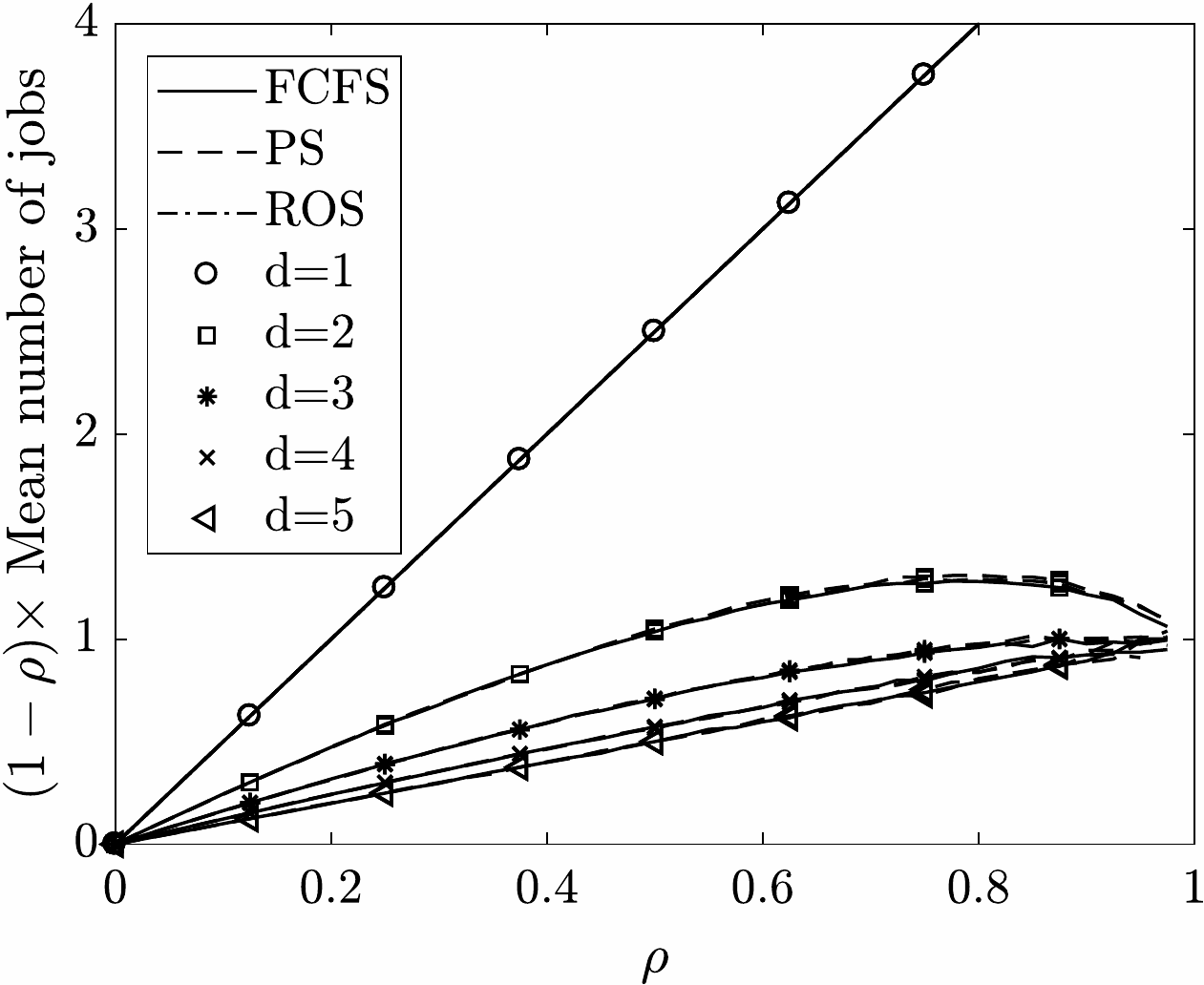}
	\end{minipage}
	\caption{Mean number of jobs for the homogeneous server system ($K=5$) with exponential service times and i.i.d.\  copies vs.\ the load for FCFS, PS and ROS service policies.}
	\label{fig:IID}
\end{figure*}

\begin{figure*}[t]
	\centering
	\begin{minipage}[b]{.45\textwidth}
		\includegraphics[width=1\textwidth]{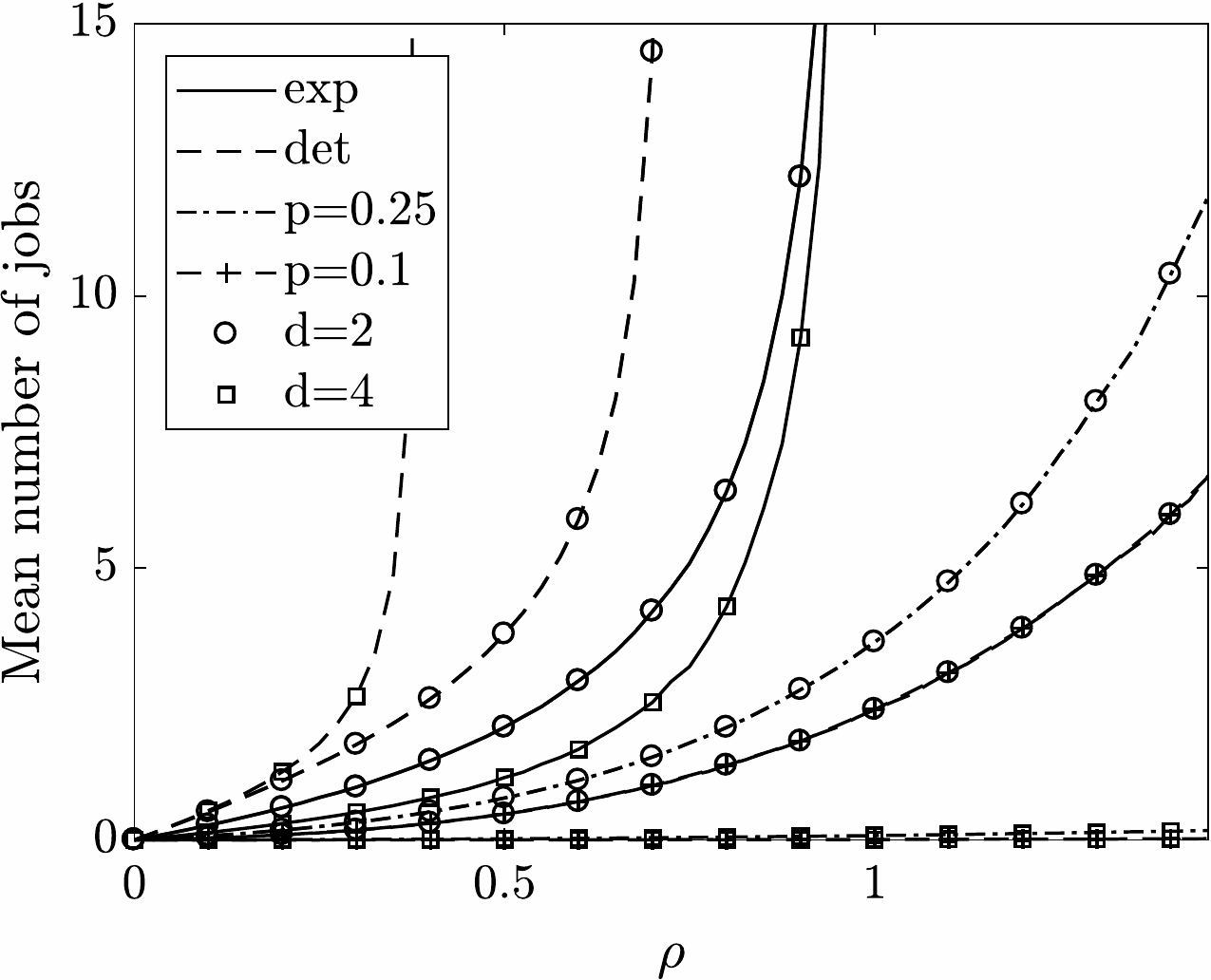}
	\end{minipage}
	\begin{minipage}[b]{.45\textwidth}
		\includegraphics[width=1\textwidth]{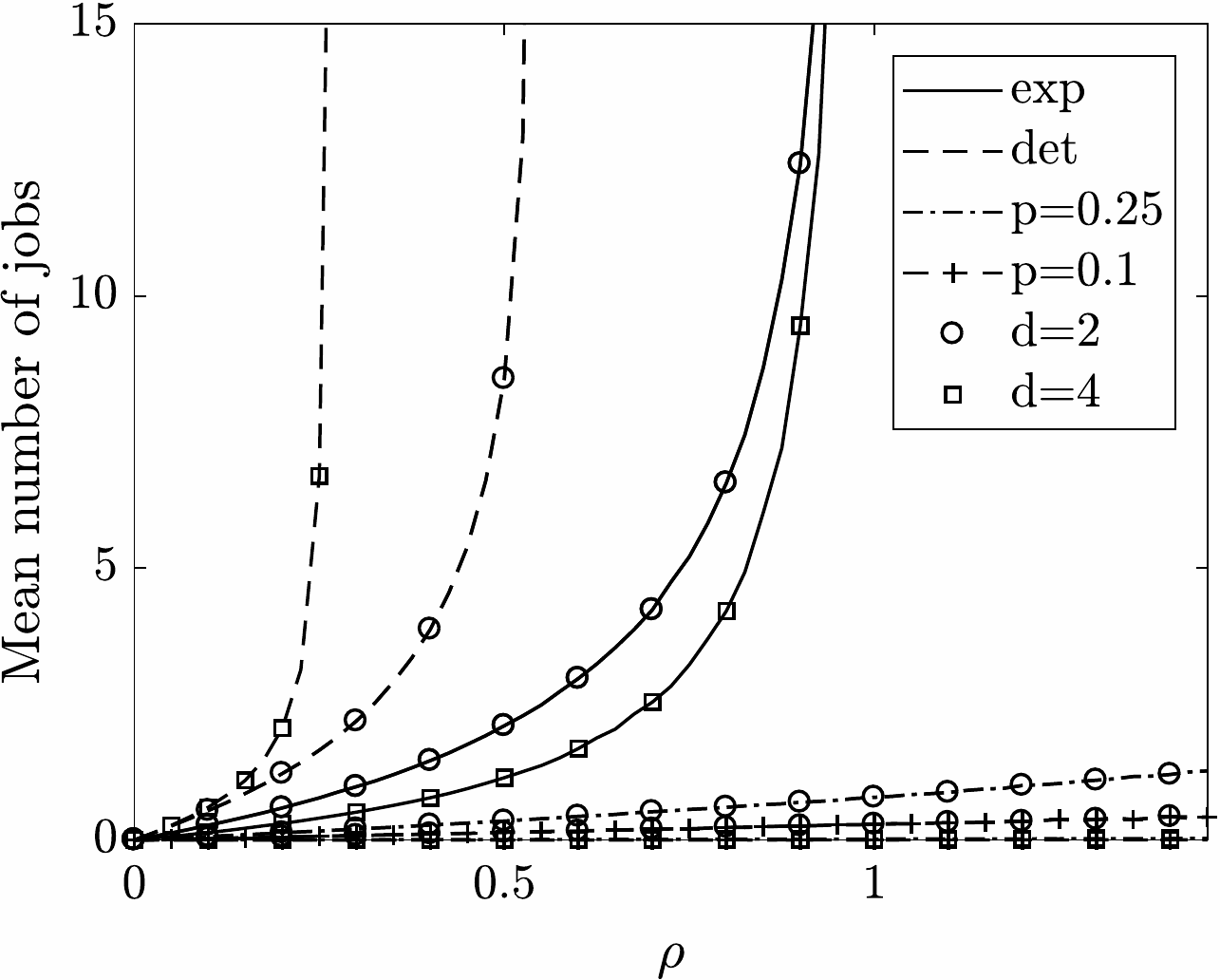}
	\end{minipage}     \\
	\begin{minipage}[b]{.45\textwidth}
		\includegraphics[width=1.1\textwidth]{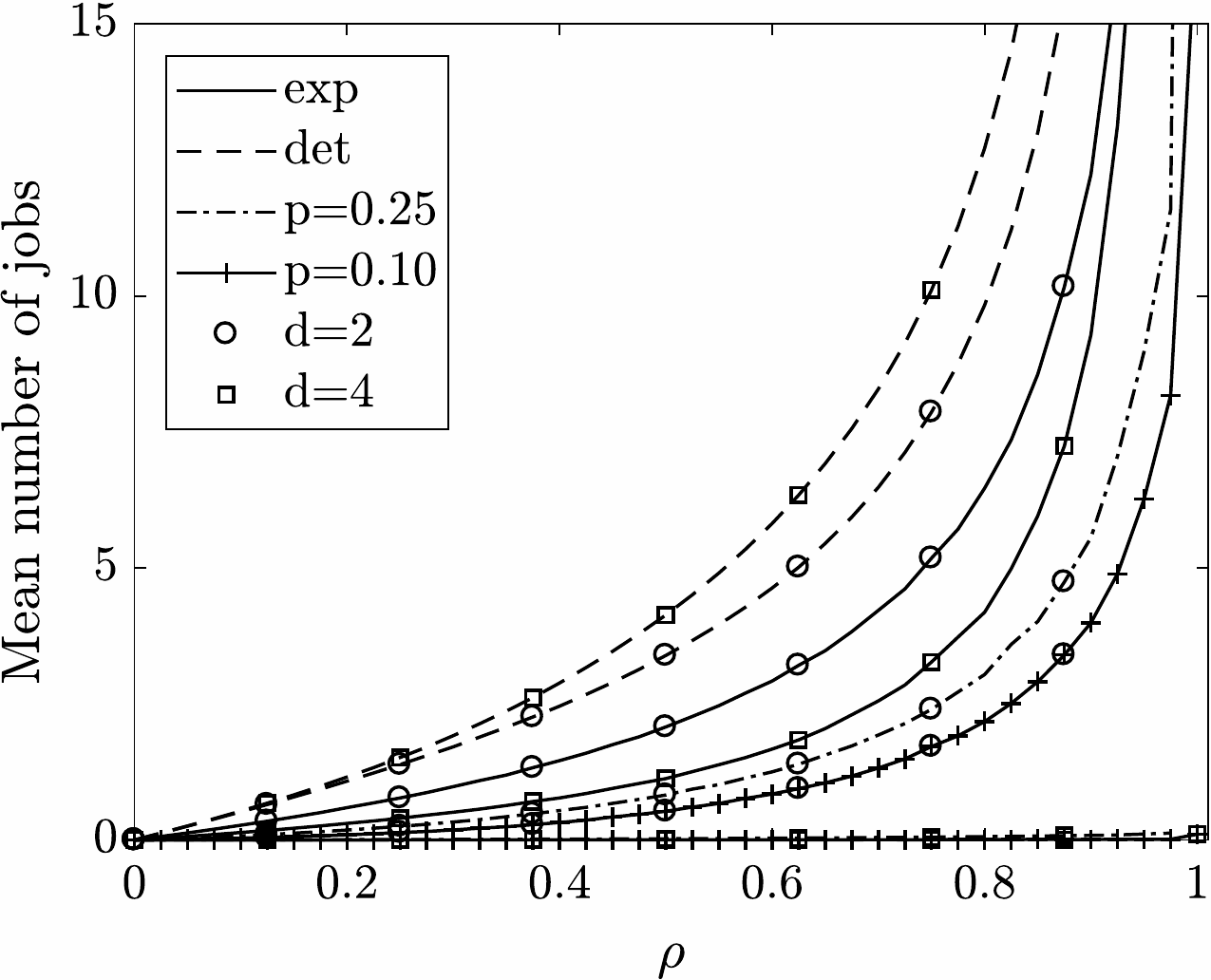}
	\end{minipage}        
	\caption{Mean number of jobs for the homogeneous server system for exponential, deterministic and degenerate hyperexponential ($p=0.25$ and $p=0.1$) service times (i.i.d.\ copies) vs.\ the load for {\em (top left)} FCFS, {\em (top right)} PS and {\em (bottom)} ROS.}
	\label{fig:comb_IID}
	
\end{figure*}

\section{Simulation analysis}
\label{sec:num}
We have implemented a simulator in order to assess numerically the impact of redundancy. We run these simulations for a sufficiently large number of busy periods ($10^6$), so that, the variance and confidence intervals of the mean number of jobs in the system are sufficiently small.

We simulate the system under the same assumptions as considered in the theoretical results, that is, exponential service times and homogeneous servers. In order to assess the impact of our modeling assumptions, we also simulate the queueing model with other service time distributions, such as deterministic services, or degenerate hyperexponential distributions. Under the latter distribution, with 
probability $p$ the service requirement is exponentially distributed with parameter $\mu p$, and is $0$ otherwise, hence, the mean service time equals $1/\mu$ (independent of $p$). The squared coefficient of variation however equals $\frac{2}{p}-1$, which increases as $p$ decreases. As a consequence, this distribution allows us to study the impact of the service time variability on the performance.  
We also consider a system with heterogeneous servers and present some preliminary analysis  in this setting.

Without loss of generality, throughout this section we assume that the mean service requirement of a copy equals 1.
In Section~\ref{sec:numIID} we present the numerics for i.i.d.\ copies and in Section~\ref{sec:numid} for identical copies. 

\subsection{IID copies}
\label{sec:numIID}
In this section we consider that copies are  i.i.d. copies. Under FCFS, PS and ROS, the system is stable whenever $\rho<1$. In Figure~\ref{fig:IID} we plot the mean number of jobs (scaled by $1-\rho$) under these policies for different values of $d$. 
For a given $d$, we observe that the plots under FCFS, PS and ROS are very similar. 
In addition, we observe that increasing the number of redundant i.i.d.\ copies, $d$,  improves the performance. 

In Figure~\ref{fig:comb_IID} we plot the mean number of jobs under FCFS, PS, and ROS for exponential, deterministic, and degenerate hyperexponential ($p=0.25$ and $p=0.1$) service time distributions. We assume $K=5$ servers and plot the performance for  $d=2$ and $d=4$ copies. 
For either FCFS, PS or ROS, we can draw similar qualitative observations: \emph{(i)}
For  variable service distributions, the performance improves as $d$ increases while it is the other way around for deterministic copies, and \emph{(ii)} for a given $d$,  the performance improves as the variability of the service time distribution increases. Only under FCFS and PS, with the  degenerate hyperexponential distribution the system remains stable even if $\rho>1$, while the stability region with deterministic service requirements seems to be reduced. 
The increase in the stability region when the service requirements become more variable was proved in an asymptotic regime by \cite{Raaijmakers2018} for FCFS. 
In general, this can be intuitively explained by noting that when copies of a job are i.i.d., the probability that a job departs due the completion of a rather large  copy will become small as the variability in the copies increases. 
For PS, the increase in performance due to variability of service sizes is even more profound than with FCFS (see also Figure~\ref{fig:comb_IID}), as only jobs that have a positive service time for their $d$ copies will enter service (which happens with probability $p^d$), while all other jobs are served instantaneously.

Under ROS, simulation results seem to indicate that the stability condition remains $\rho<1$ for any service time distributions. Since copies are being randomly chosen for service, the system does not seem to profit from the variability of the service times of the i.i.d. copies. For example, in the case of degenerate hyperexponential distribution, when $\rho$ is close to 1, with probability $p$ a job will start being served in a server where its copy has strictly positive service requirement. Due to ROS, in a high congested system, the probability that another copy (possible of size 0) of this job will receive service in another server will be close to zero. Hence, with probability $p$ a job needs exponential service time with parameter $\mu p$ and with probability $1-p$ its service time equals zero. On average it needs $1/\mu$, which explains why $\rho<1$ can be the stability condition.
Further note that for $d=4$, it seems that the system remains stable when $\rho= 1$. This is however not the case. For $\rho$ close to 1, the  mean number of jobs in the system is close to zero, which can be explained as follows: A job has a strictly positive service requirement with probability $p^d$. When $d=4$ and   $p\in \{0.25, 0.1\}$, it holds that $p^4<10^{-2}$. Hence, it is very likely that for a very long time, zero jobs are present in the system (as all arriving jobs spend zero time in the system). This explains why for $\rho<1$, the mean number of jobs stays very close to zero.  We however believe that the stability condition is $\rho<1$.


\begin{figure*}[t]
	\centering
	\begin{minipage}[b]{.45\textwidth}
		\includegraphics[width=1.\textwidth]{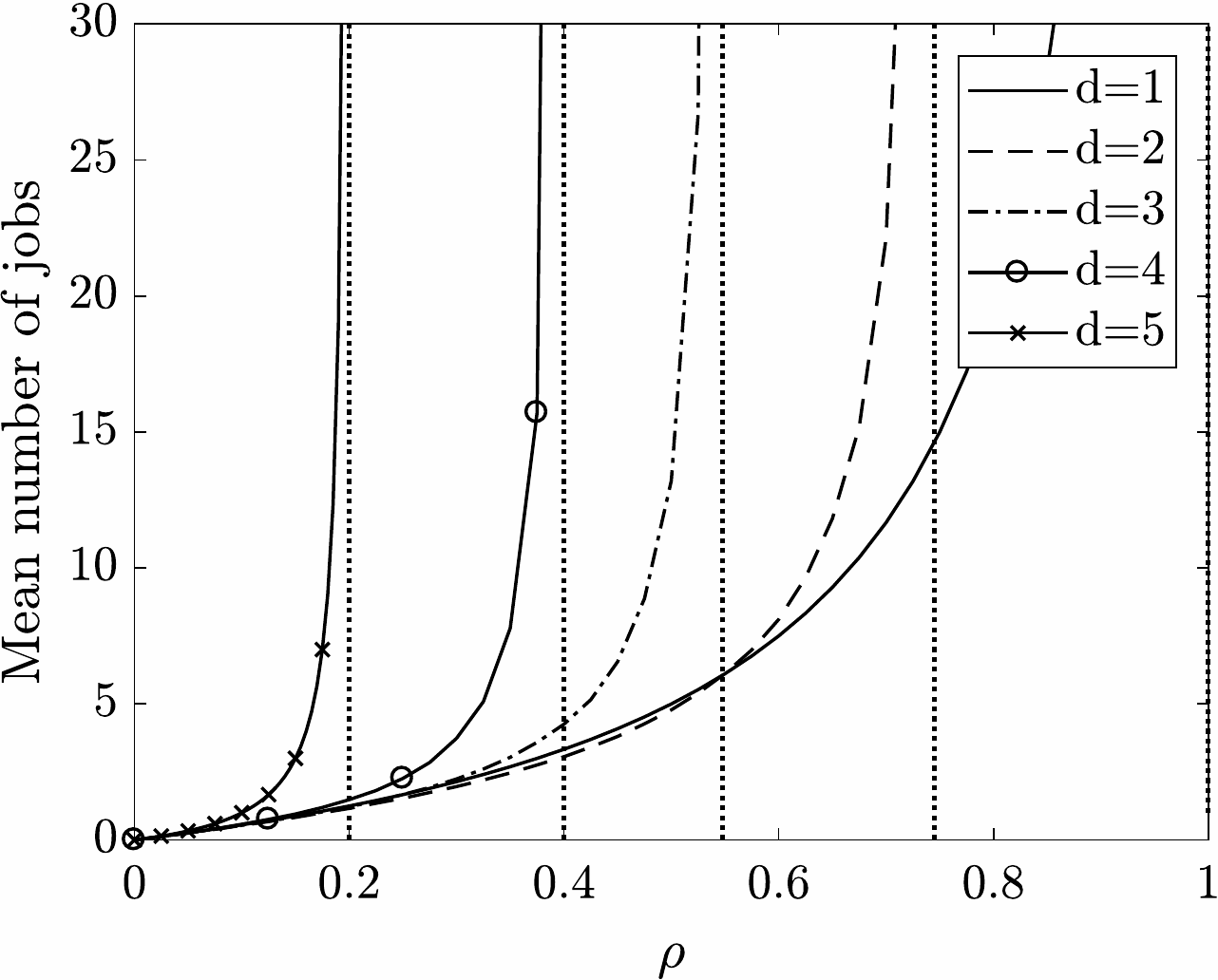}
	\end{minipage}
	\begin{minipage}[b]{.45\textwidth}
		\includegraphics[width=1\textwidth]{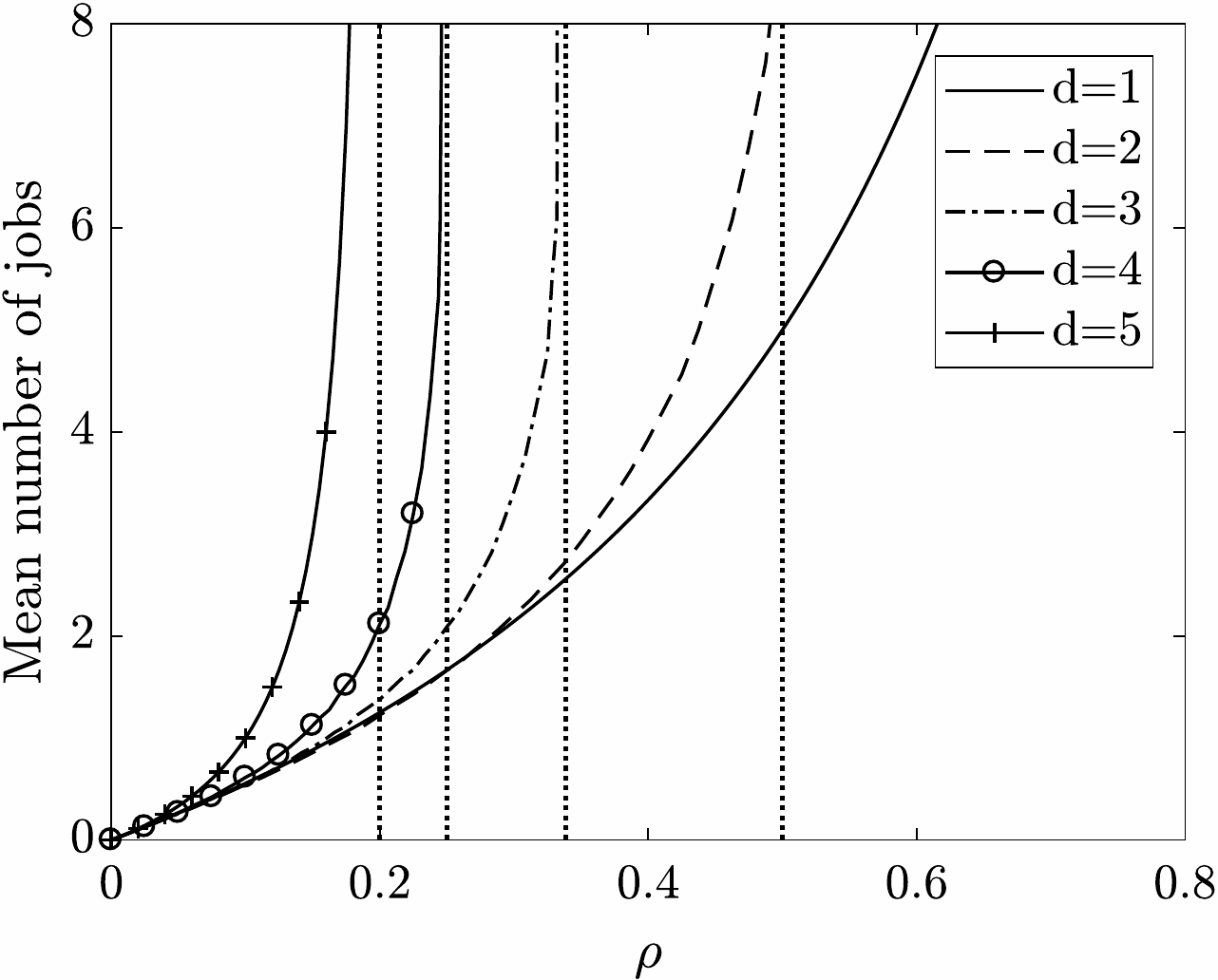}
	\end{minipage}\\
	
	\begin{minipage}[b]{.45\textwidth}
		\includegraphics[width=1.\textwidth]{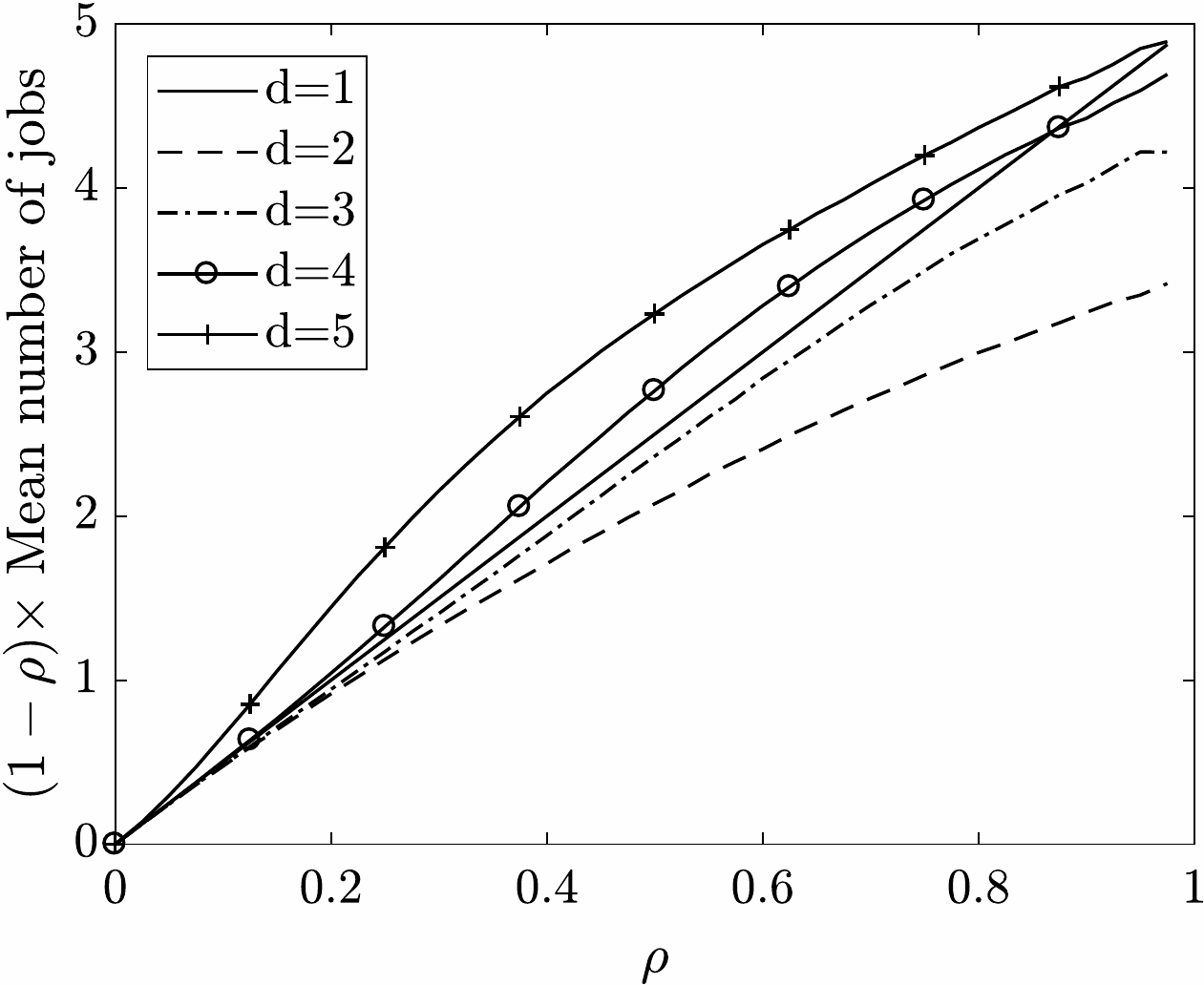}
	\end{minipage}

	\caption{Mean number of jobs for the homogeneous server system ($K=5$) with exponential service times and identical copies vs.\ the load:   {\em (top left)} FCFS, {\em (top right)} PS and {\em (bottom)} ROS.}
	\label{fig:PS}
\end{figure*}
\subsection{Identical copies}
\label{sec:numid}
In this section we consider jobs with identical copies. We have proved that the stability condition strongly depends on the employed scheduling policy and on the number of copies $d$. In Section~\ref{sec:numexp} we evaluate the system for different  values of $d$ and observe that the stability region reduces as $d$ grows large. In Section~\ref{sec:LT} we characterize the performance and its dependence on~$d$ for a light-load regime, and observe that when the load is \emph{small enough}, redundancy can improve the performance. In Section~\ref{sec:nonexp}, we numerically  study the impact of different service time distributions on the performance. Finally, in Section~\ref{sec:heterogeneous},  we present a preliminary analysis and numerics for a heterogeneous servers setting. 

\subsubsection{Exponential service times}
\label{sec:numexp}
In Figure~\ref{fig:PS} we plot the mean number of jobs under FCFS ({\em top left}),  PS ({\em top right}) and  ROS ({\em bottom}) respectively, for different values of $d$. The vertical lines in Figure~\ref{fig:PS} correspond to the stability regions (for different values of~$d$)  as derived in Proposition~\ref{prop:unstab}, Proposition~\ref{prop:FCFS_stab} (where $\bar\ell/K$ has been obtained via simulation of the saturated system) and Proposition~\ref{prop:ROSID}. Indeed, we observe that the mean numbers of jobs under FCFS, PS and ROS have an asymptote at the point $\bar\ell/K$, $1/d$ and 1, respectively. 

An interesting observation we can draw from Figure~\ref{fig:PS} is that for every $d$, the stability region under FCFS is larger than under PS,  as proved in~Corollary~\ref{cor:FCFSPS}.

Under ROS, we observe that the best performance of the system is achieved when $d=2$, for any load $\rho$. Hence, adding redundant copies ($d=2$) when the scheduling policy is ROS \emph{does} improve the performance of the system for any value of $\rho$. This as opposed to FCFS and PS, where the best performance for high load is obtained in $d=1$.

\begin{figure*}[t]
	\centering
	\begin{minipage}[b]{.45\textwidth}
		\includegraphics[width=1\textwidth]{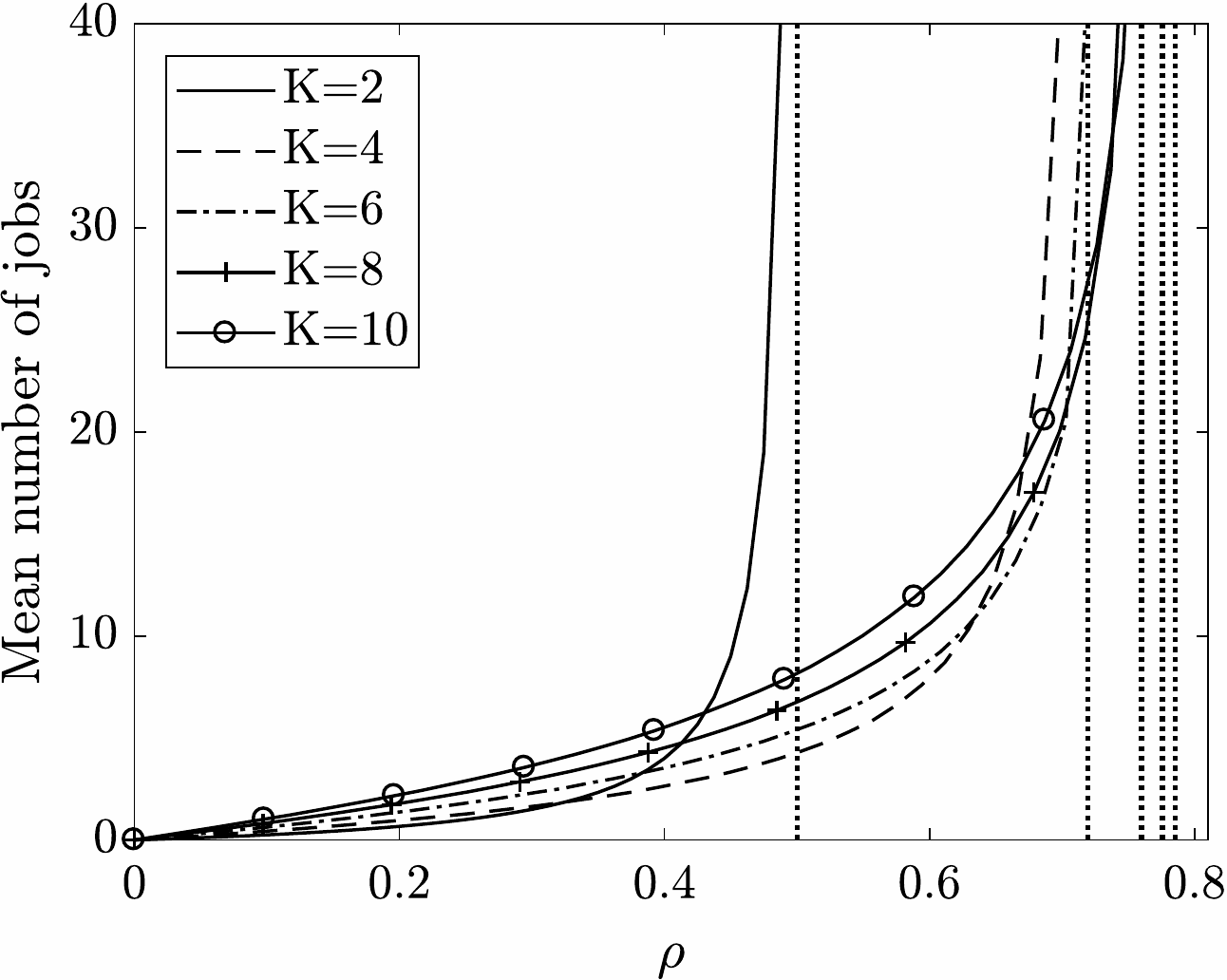}
	\end{minipage}
	\begin{minipage}[b]{.45\textwidth}
		
		\includegraphics[width=1\textwidth]{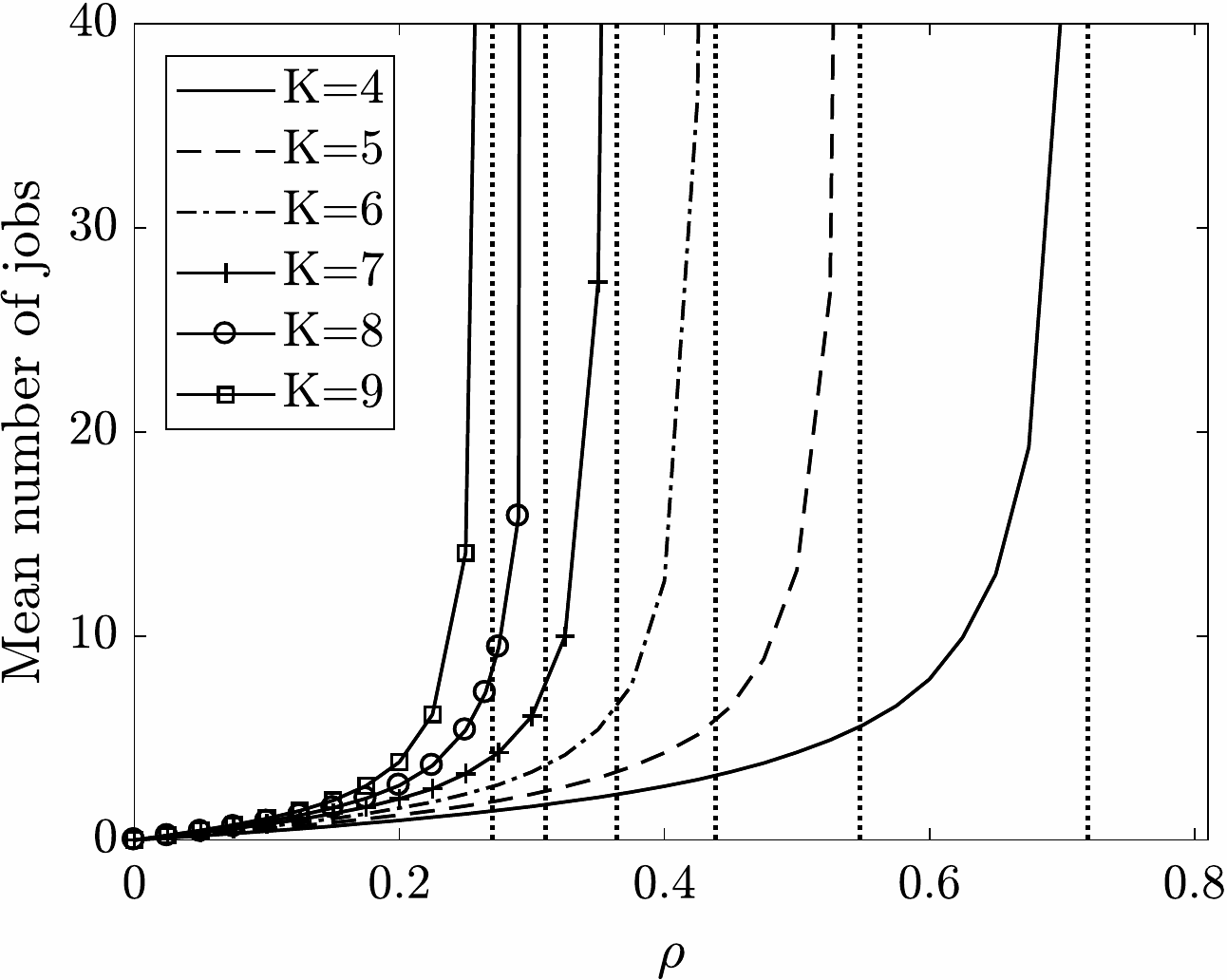}
	\end{minipage} 
	\caption{Mean number of jobs for the homogeneous server system under FCFS  with exponential service times and identical copies vs.\ the load: ({\em left}) $d=2$ and $K=2,\ldots,9$, ({\em right})  $d=K-2$ and $K=4,\ldots,10$.}
	\label{fig:FCFS_Kd2}
\end{figure*}

In Figure~\ref{fig:FCFS_Kd2} we focus on FCFS. As before, the vertical lines correspond to the stability region,  $\rho<\bar\ell/K$.
In the figure on the left, we fix the number of copies to  $d=2$ and plot the performance for several values of the number of servers $K$.  We note that the stability region increases in $K$ (as also proved in Proposition~\ref{prop:props}) and that it converges as $K$ grows large to a constant value. 
In the figure on the right, we instead set $d=K-2$, so that the number of copies increases with $K$. Now, we observe that the stability condition reduces as the number of servers $K$ increases. Hence, the negative impact due to having one more redundant copy, is more important than the benefit of having one more server.  This is in agreement with the case $d=K-1$, for which the the stability region ($\rho<\frac{2}{K}$) also decreases in $K$. 

\subsubsection{Light-traffic approximation}
\label{sec:LT}
In this section we consider the steady-state performance for extremely low traffic load, i.e., the so-called light-traffic regime pioneered in \cite{RS88a}, see also \cite{Walrand1990519}. 
The light-traffic approximation corresponds to the first-order asymptotic expansion of the system as $\lambda \to 0$.
More precisely, as $\lambda \to 0$ we seek to write $E(\lvert \vec N^P(\infty) \rvert) =   \bar{N}^{LT,P}(\lambda) + o(\lambda^2)$, for a given service policy~$P$. We defer the details of the light-traffic analysis to Appendix E, and we give here the main result of the approach in which we characterize $\bar{N}^{LT,FCFS}(\lambda)$, $\bar{N}^{LT,ROS}(\lambda)$ and $\bar{N}^{LT,PS}(\lambda)$.

\begin{figure*}[t]
	\centering
	\begin{minipage}[b]{.45\textwidth}
		\includegraphics[width=1\textwidth]{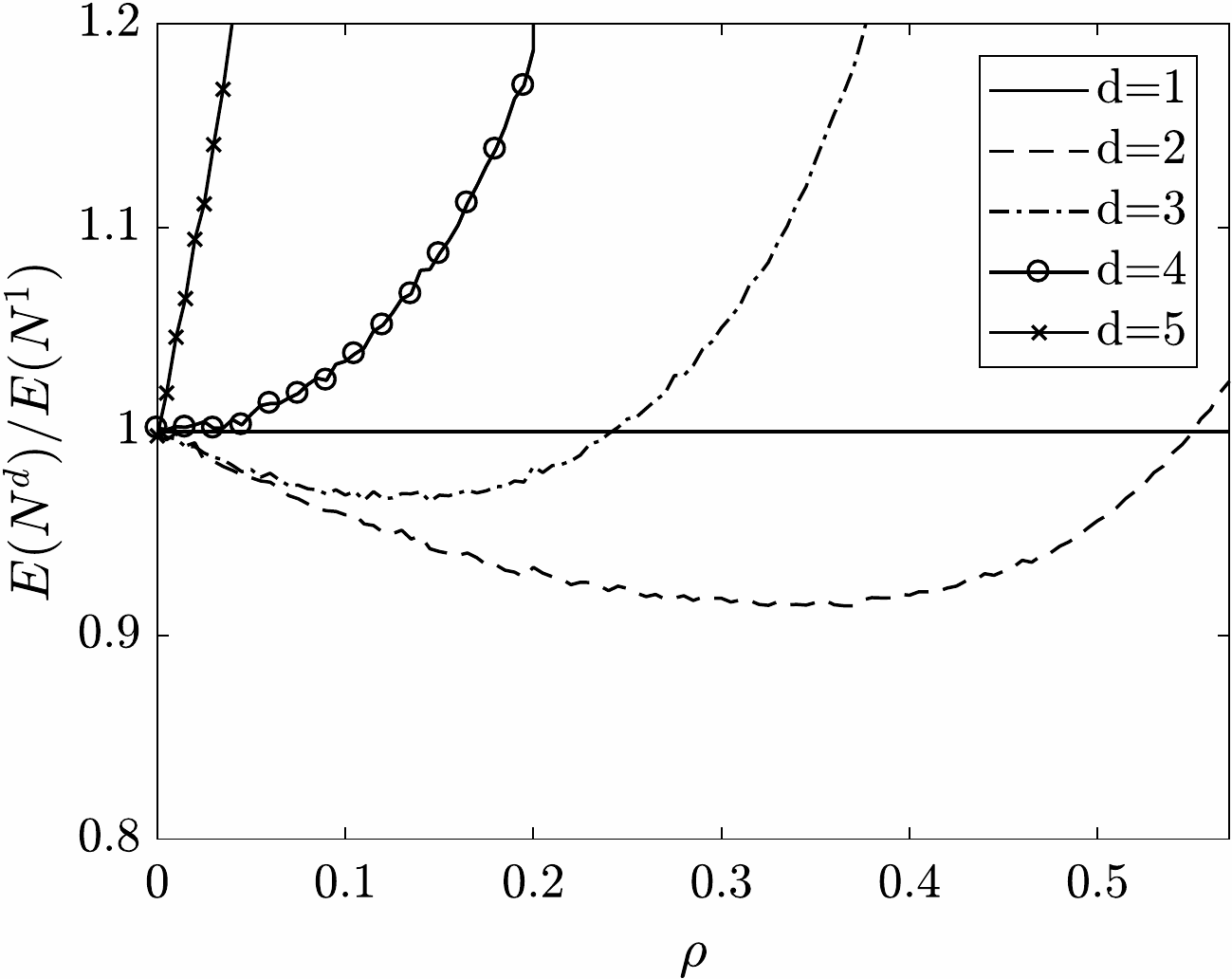}
	\end{minipage}
	\begin{minipage}[b]{.45\textwidth}
		\includegraphics[width=1.05\textwidth]{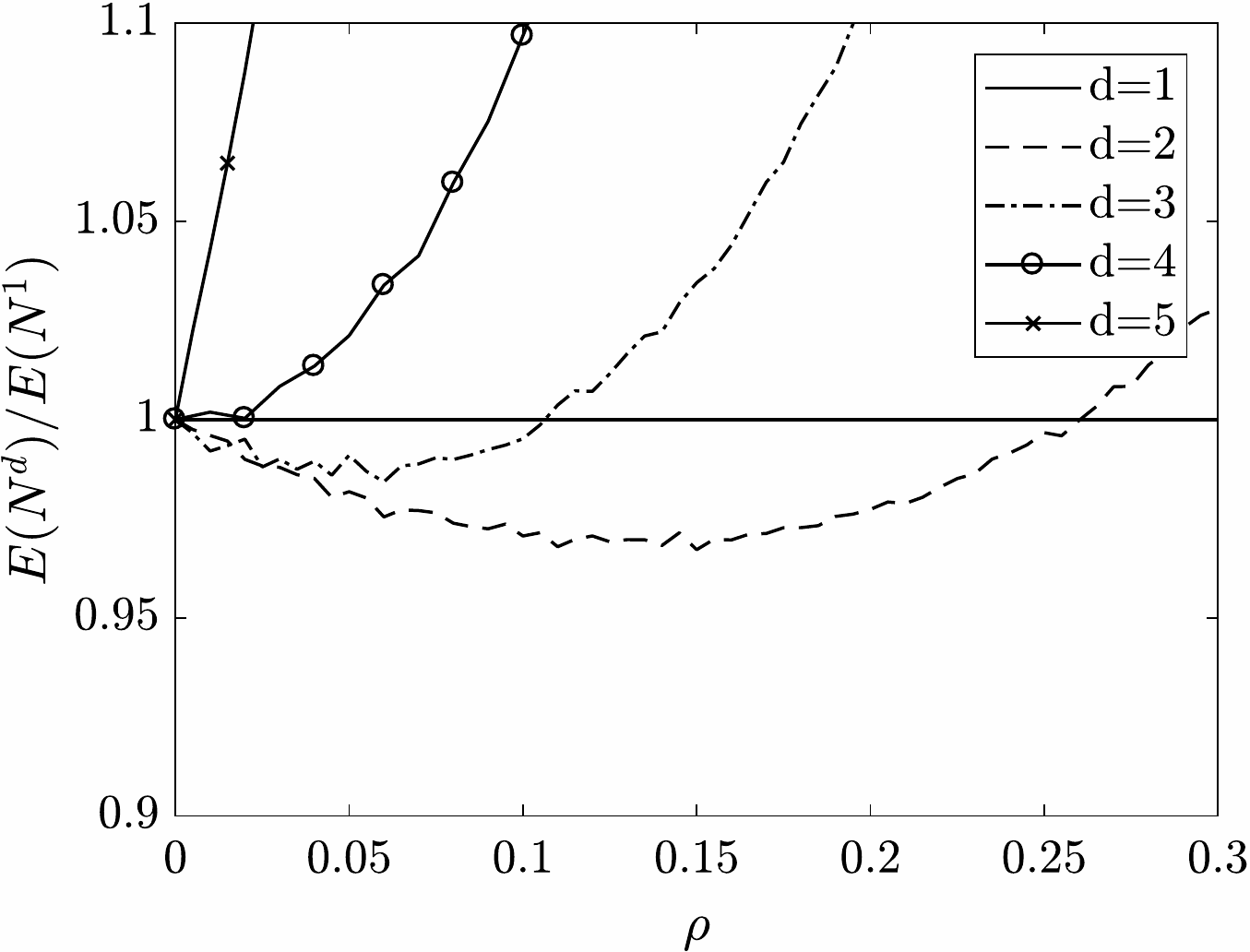}
	\end{minipage}              
	\caption{Ratio of the  mean delay with $d$ identical copies and the mean delay with no redundant copies ($d=1$), as a function of $\rho$. For the homogeneous server system ($K=5$) with exponential service times and identical copies: {\em (left)} FCFS and {\em (right)} PS. }
	\label{fig:PS-FCFS-LT}
\end{figure*}

\begin{lemma}\label{lem:LT_2} The leading term  of the light-traffic approximation  for FCFS, ROS and PS with identical copies is given by $\bar{N}^{LT,FCFS}(\lambda)  = \bar{N}^{LT,ROS}(\lambda) = \frac{\lambda}{\mu} + \frac{3\lambda^2}{2\mu^2}\frac{1}{\binom{K}{d}}$, and $\bar{N}^{LT,PS}(\lambda) = \frac{\lambda}{\mu} + \frac{\lambda^2}{\mu^2}\frac{1}{\binom{K}{d}}$, respectively.
\end{lemma}

We note that for all three policies, the light-traffic term is minimized in $d^*:=\lfloor K/2\rfloor$. To explain this, we note that at very low loads, an arriving  job will find at most one other job present. In particular this implies that this new arrival will wait for service if and only if it is of the same type as the job already present in the system. The probability of being of the same type is equal to $1/\binom{K}{d}$, which is minimized by setting $d$ equal to $d^*$.

For ROS, we saw in Figure~\ref{fig:PS} that $d=2$ indeed minimizes the mean number of jobs.
In Figure~\ref{fig:PS-FCFS-LT} (left) and (right) we consider   FCFS and PS, respectively, for low load.  We plot the ratio of the  mean total number of jobs for the system with $d$ identical copies with that of a system with no redundant copies ($d=1$). If the ratio is below 1, this implies that redundancy (for the particular value of $d$) improves the performance. As predicted in Lemma~\ref{lem:LT_2}, redundancy reduces the mean delay for $\rho$ small enough, and the best performance is obtained in $d=2$.
Recall that for sufficiently large load, the minimum delay is obtained with $d=1$, see Figure~\ref{fig:PS}.

\subsubsection{Non-exponentially distributed service requirements}
\label{sec:nonexp}

\begin{figure*}[h]
	\centering
	\begin{minipage}[b]{.45\textwidth}
		\includegraphics[width=1.015\textwidth]{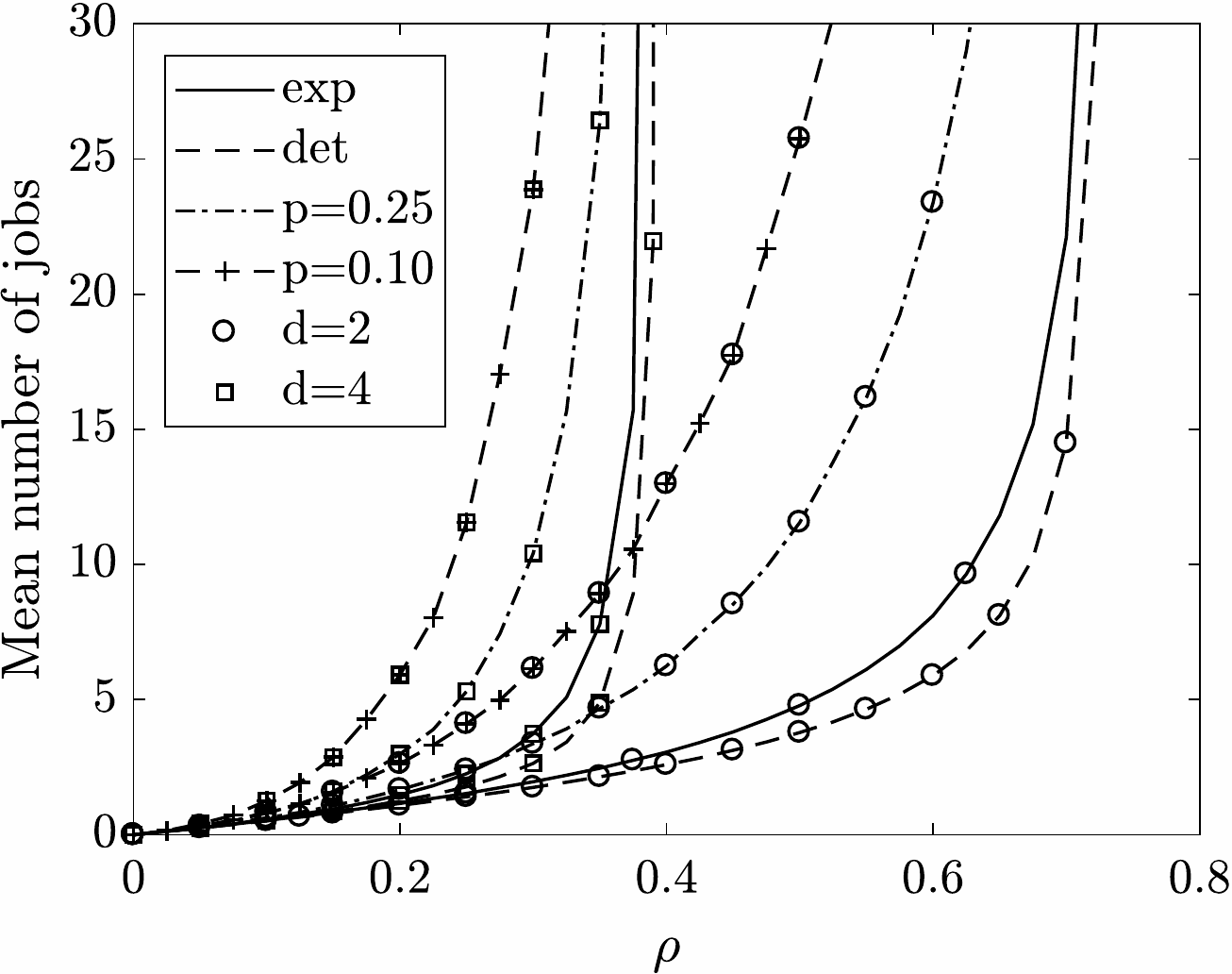}
	\end{minipage}
	\begin{minipage}[b]{.45\textwidth}
		\includegraphics[width=1.\textwidth]{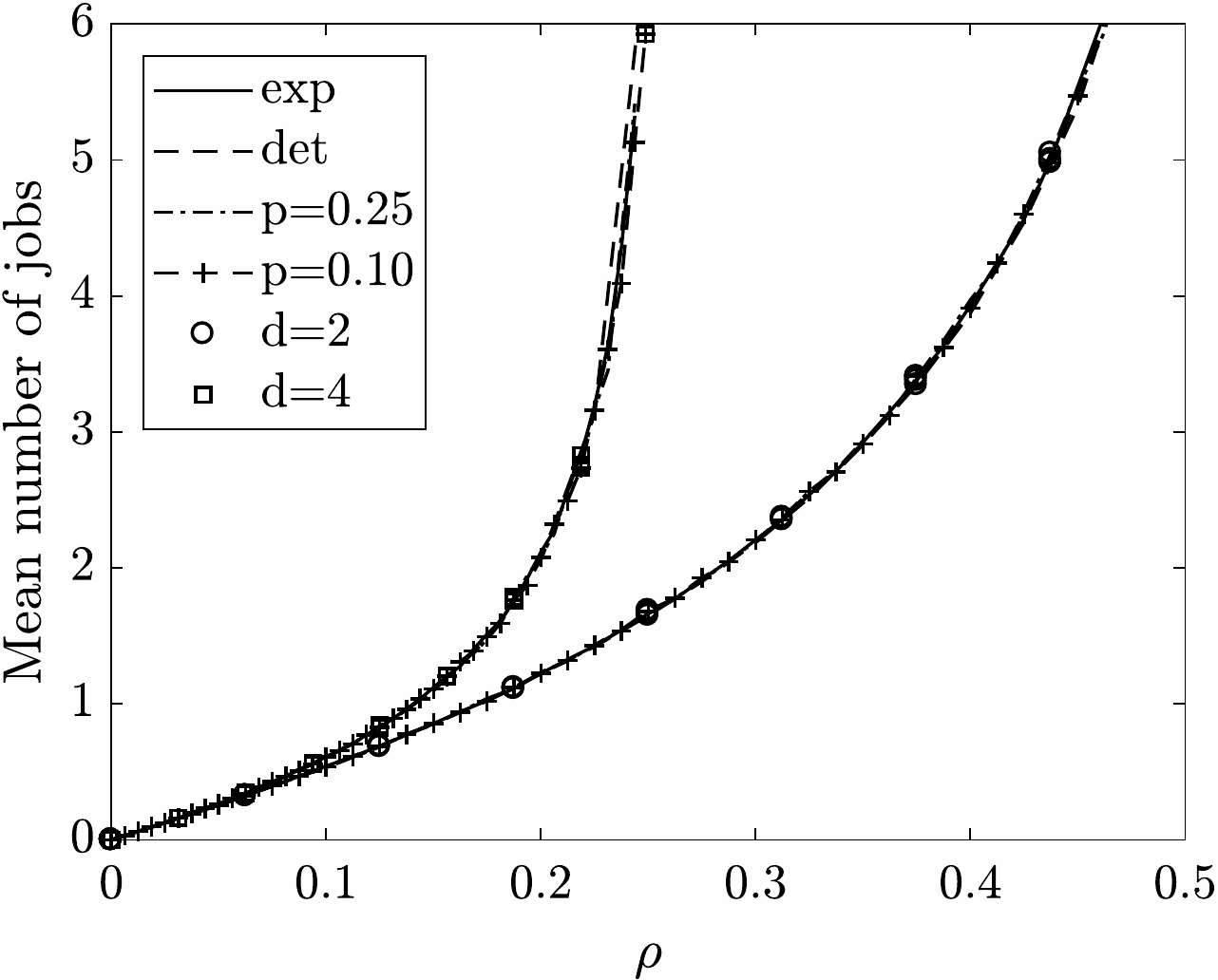}
		
	\end{minipage}    \\
	
	\begin{minipage}[b]{.45\textwidth}
		\includegraphics[width=1.\textwidth]{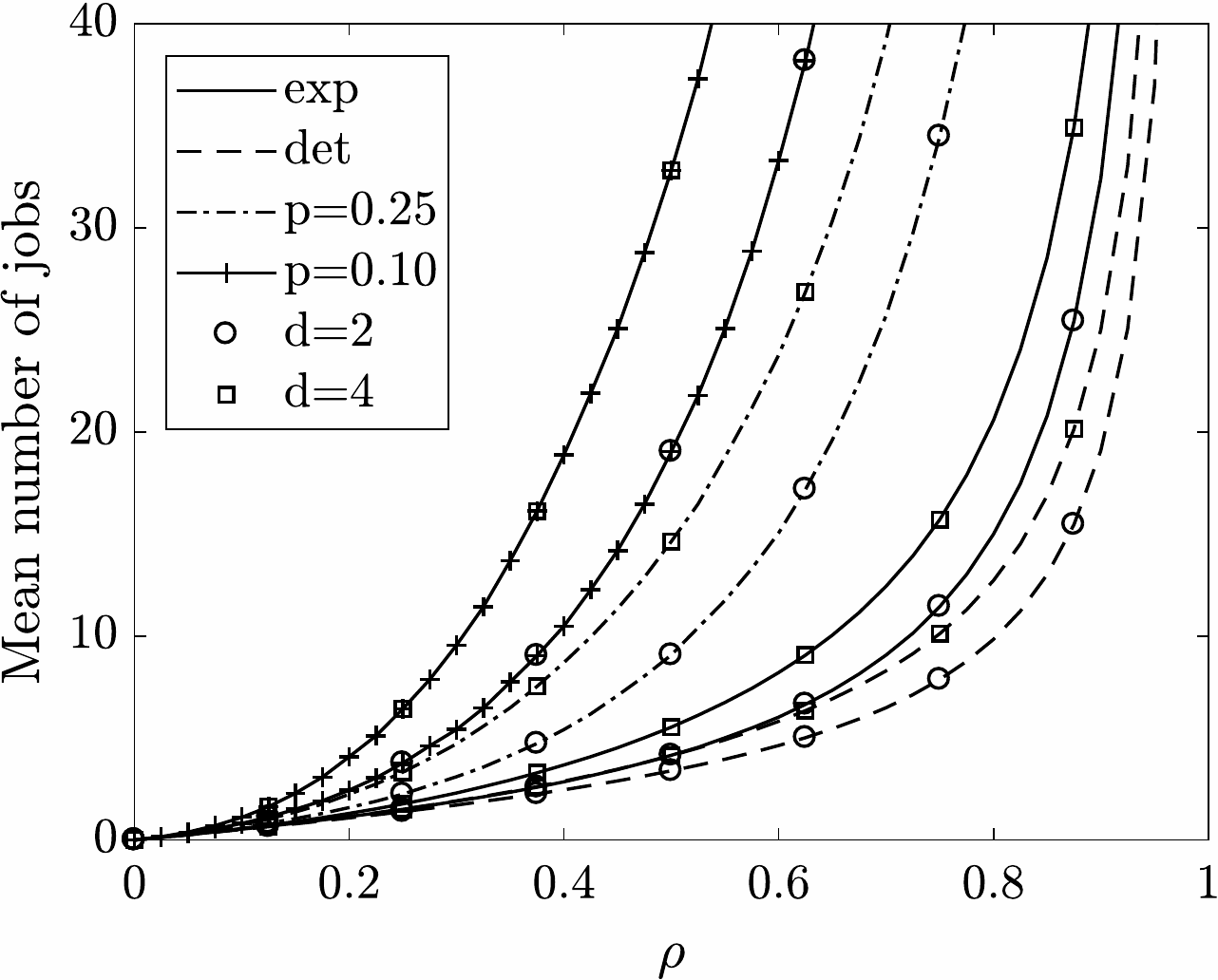}
		
	\end{minipage}            
	\caption{Mean number of jobs for the homogeneous server system ($K=5$) and  exponential, deterministic and degenerate hyperexponential ($p=0.25$ and $p=0.1$) service times (identical copies) vs.\ the load: {\em (top left)} FCFS, {\em (top right)} PS and {\em (bottom)} ROS.}
	\label{fig:PS-FCFS-gen}
\end{figure*}

In Figure~\ref{fig:PS-FCFS-gen} we  compare the mean number of jobs for   exponential, deterministic, and degenerate hyperexponential service time distributions. We consider $K=5$ servers and $d=2$ and $d=4$ identical copies. 

We observe that  for  FCFS, PS and ROS,   the performance degrades as $d$ increases. This is in contrast to the i.i.d. case, where we observed the opposite effect. 
This is due to the fact that with identical copies, capacity is wasted on serving the exact same copy, while with i.i.d.\ copies, the system benefits from the difference in the requirement per copy.  

For FCFS and ROS, 
we observe that, unlike in the i.i.d. case, the performance of the system degrades as the variability of the service time increases. 
In particular, for a given $d$, the best performance is obtained with deterministic service times. Moreover, for the degenerate hyperexponential service distribution, the  performance deteriorates as $p$ decreases. 
For FCFS, these observations are in agreement with the results obtained by \cite{HH18} for the mean field analysis. 

From the numerics, it seems that for deterministic or degenerate hyperexponential service requirements, ROS is more stable than when FCFS and PS are implemented.

Another interesting observation is that for PS the performance seems to be nearly insensitive to the service time distribution (beyond its mean value). When degenerate hyperexponential service times are considered, it is trivial that the performance coincides with that of exponential service requirements. This can be explained as follows: With identical copies, only a fraction $p$ of the arrivals have a non-zero service requirement, and this is exponentially distributed with mean $1/(p\mu)$. Thus, the system with degenerate hyperexponential service requirements with parameter $p$ is equivalent to the system with arrival rate $\lambda$ and exponentially distributed service requirements with mean $1/\mu$ where time is parametrized with parameter $p$.
\subsubsection{Heterogeneous server capacities}
\label{sec:heterogeneous}
In this section we investigate the stability region of the previously analysed systems PS and FCFS for \emph{heterogeneous} servers. 
We take exactly the same model as before, that is, a type-$c$ job arrives at rate $\lambda/\binom{K}{d} $ and sends $d$ identical copies (exponentially distributed with parameter $\mu$) to 
$d$ servers chosen at random.
However, now, instead of having homogeneous servers, we   assume that server~$s$ has capacity $\nu_s$, for $s\in S$. 
This is a rather simple heterogeneous model, as the arrival rates of the different types are taken uniformly, but in spite of this, it provides  interesting insights. In fact,  redundancy might have the complete opposite effect when the capacity of the servers is sufficiently spread out. 

When $d=1$, there is no redundancy and each server receives arrivals at rate $\lambda/K$.   For any work-conserving policy implemented in the servers, the latter system is stable if and only if  $\rho< \nu_{min}$, where $\nu_{min}=\min_{s\in S} \{\nu_s\}$. 

When $d=K$, each job sends identical copies to all $K$ servers. Assume one starts with an empty system at time~0. It can easily be seen that in each server the queue length is the same, and a departure of a job is always due to a copy finishing its service requirement in the server with the highest capacity. This holds for both FCFS and PS. Hence, the system behaves as a single server with capacity $\mu \nu_{max}$, where $\nu_{max}=\max_{s\in S} \{\nu_s\}$. Therefore, under FCFS or PS with $d=K$, the system is stable if and only if $\lambda<\mu \nu_{max}$, i.e.,\ $\rho<\frac{\nu_{max}}{K}$. 
From this, we observe that adding $d=K$ identical copies to the system reduces the stability region if and only if $\nu_{max}<K\nu_{min}$. Hence, when the difference between the smallest and largest capacity is not that large, redundancy reduces the stability region, as we saw for the homogeneous case. However, when the difference is sufficiently large, adding $K$ redundant copies to the system can in fact improve the system.  

\begin{figure*}[h]
	\centering
	\begin{minipage}[b]{.45\textwidth}
		\includegraphics[width=1.\textwidth]{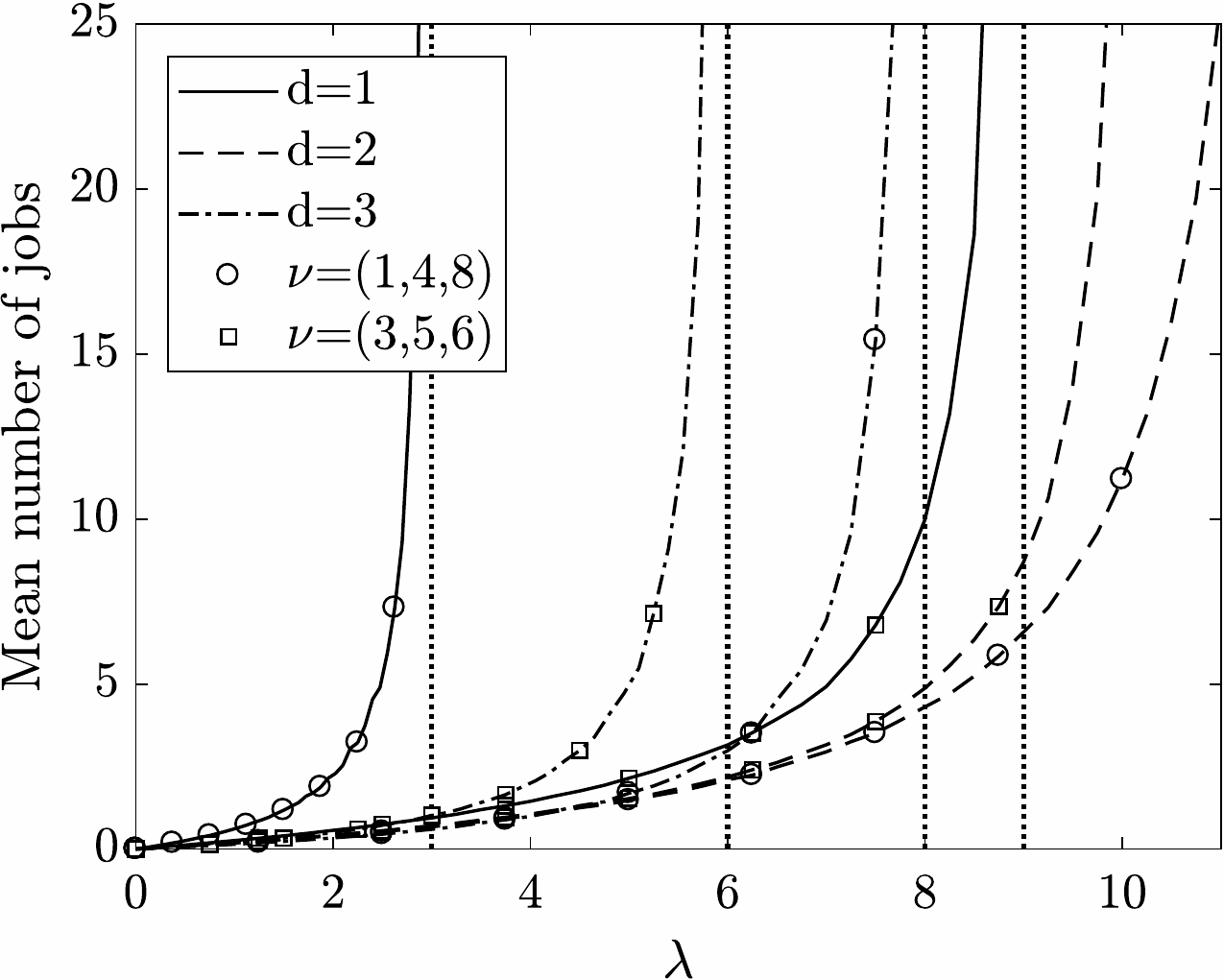}
	\end{minipage}
	\begin{minipage}[b]{.45\textwidth}
		\includegraphics[width=1.\textwidth]{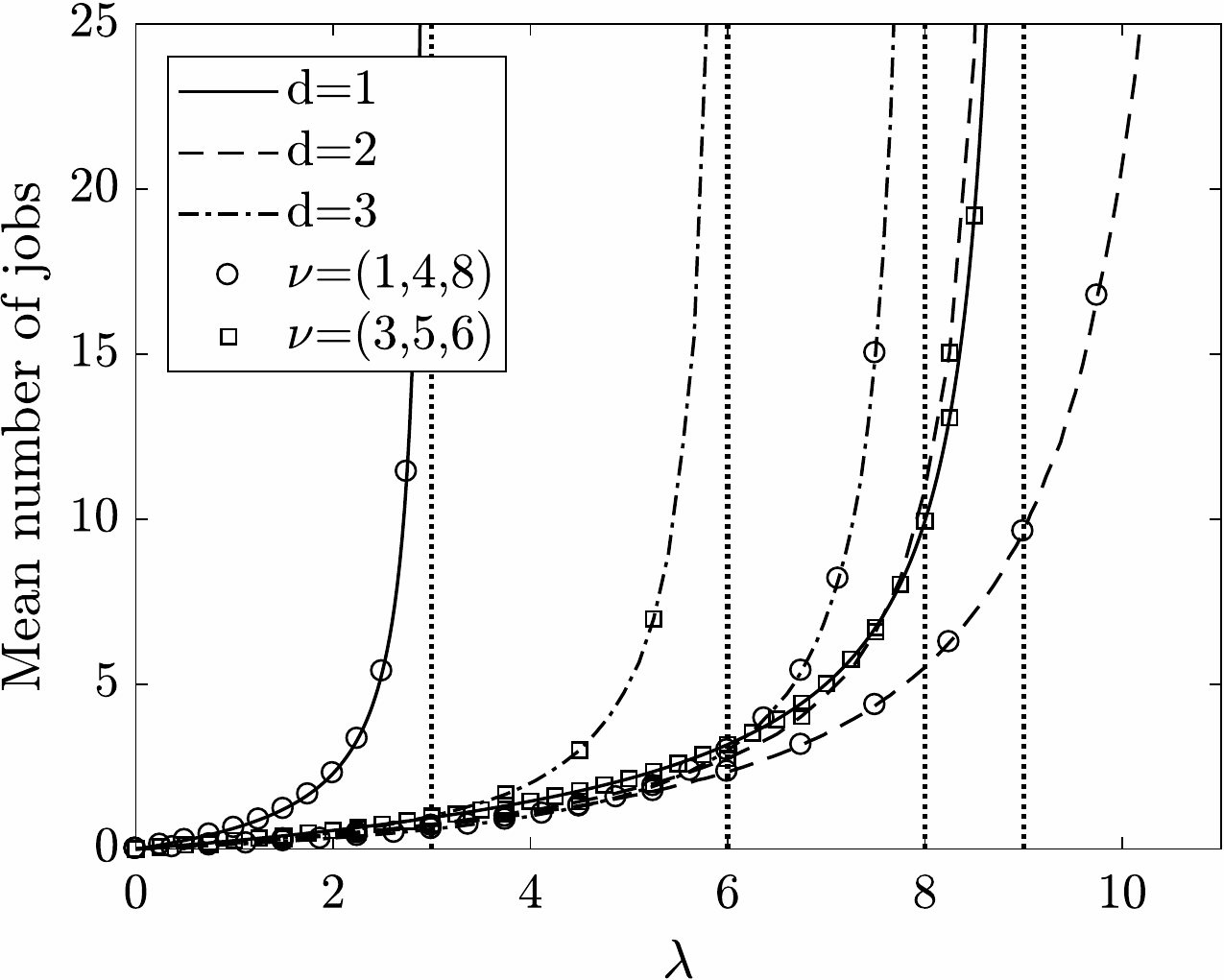}
	\end{minipage}     
	
	\caption{Mean number of jobs with heterogeneous servers ($K=3$ and $\nu=(1,4,8)$ and $\nu=(3,5,6)$) with exponential service times (identical copies) vs.\ the arrival rate ($\lambda$): {\em (left)} FCFS and {\em (right)} PS. }
	\label{fig:heterPS-FCFS}
\end{figure*}

To see the impact of  $d<K$ identical copies, we performed simulations. 
In Figure~\ref{fig:heterPS-FCFS} we plot the mean number of jobs under FCFS and PS, respectively, for $K=3$ and two different setting of parameters: $(\nu_1,\nu_2,\nu_3)=(3,5,6) $ (lines with~$\square$) and $(\nu_1,\nu_2,\nu_3)= (1,4,8)$ (lines with $\circ$). Note that when $(\nu_1,\nu_2,\nu_3)=(3,5,6)$, the stability condition for $d=1$ is $\lambda<9$, and for $d=K$ is $\lambda<6$. As in the case of homogeneous servers, the stability condition is larger in the system with no redundancy than when $d=K$. However, when the server capacities are $(\nu_1,\nu_2,\nu_3)= (1,4,8)$ (lines with $\circ$), we observe the opposite effect: the stability condition for $d=1$ is $\lambda<3$, and for $d=K$ is $\lambda<8$.  Since $\nu_{max} < K\nu_{min}$, having $d=3$ redundant copies improves the stability of the system.

From Figure~\ref{fig:heterPS-FCFS} we observe that for FCFS, the performance of the system is improved when $d=2$ copies are considered. In addition, the stability region under $d=2$ seems to be larger than that under $d=1$. Similarly for PS service policy we observe that under different capacity parameters, the stability condition of the system could be improved when $d=2$ copies are considered.

Even if these results are not conclusive for the study of the stability condition under heterogeneous server capacities and identical copies, these results are insightful for a future understanding of the impact of redundancy in heterogeneous server systems.
As a general conclusion, we see that heterogeneity in server's speeds has a profound impact on the stability condition with identical copies. 


\section{Conclusion}
\label{sec:conclusion}
In recent years, redundancy has emerged as a promising technique to reduce the response time of jobs in data centers, and researchers have obtained encouraging results showing that indeed, redundancy could help improving the performance. Due to mathematical tractability, a large body of the literature has assumed that redundant copies are exponentially distributed and independent among each other,  and that the scheduling discipline in servers is FCFS. Under these assumptions, one of the main conclusions from literature is that redundancy does not impact the stability region, that is, the amount of work that the system can handle remains unchanged.

However, we believe our analysis serves as an indication that redundancy needs to be implemented with care, in order to prevent an unnecessary degradation of the performance. Indeed, the main takeaway message from our work is that the stability region strongly depends on the modeling assumptions (in some cases in a somehow unexpected manner), for instance on \emph{(i)} the scheduling discipline deployed in servers, \emph{(ii)} the correlation structure between copies, \emph{(iii)} service requirement distribution, and \emph{(iv)} the variability of server speeds.
The fact that the stability condition  depends strongly on the employed scheduling discipline was initially unexpected for us, given that in the $M/M/1$  case, size-unaware policies like PS, ROS and FCFS have the same steady-state distribution with exponential service times. 

An important question for future work is to characterize a set of disciplines that, just like ROS, do not reduce the stability region under the identical copies assumption,  and  to analyze their performance as a function of the redundancy degree~$d$. 
Other interesting questions that emerge from our work are to investigate the stability of redundancy when combined with size based scheduling policies like Shortest-Remaining-Processing-Time and Least Attained Service. 
All our theoretical results are restricted to the case of exponentially distributed service times. To characterize the stability condition to general service times  is a very challenging problem and it will require a different proof technique. For instance, to extend the maximum stability result of alpha-fair bandwidth allocations in networks from exponential assumptions (proved in \cite{BM01}) to general distributions (proved in \cite{PFTA12}) took a decade. 

We conclude by noting that in our study we have considered a rather basic model of redundancy.
In practice, one might expect servers to be heterogeneous in terms of speeds and scheduling discipline,
and the service time distribution of copies to show other correlation structures. 
However, we believe that our analysis provides sufficient ground to conclude that redundancy needs to be implemented with care, in order to prevent an unnecessary degradation of the performance.

\section*{Acknowledgments}
The authors are thankful to Rhonda Righter for useful comments that helped improve the presentation of the paper. The Ph.D. project of E. Anton is funded by the French ``Agence Nationale de la Recherche (ANR)" through the project ANR-15-CE25-0004 (ANR JCJC RACON). This work (in particular research visits of E. Anton and M. Jonckheere) was partially funded by a STIC AMSUD GENE project. Urtzi Ayesta has received funding from the Department of Education of the Basque Government through the Consolidated Research Group MATHMODE (IT1294-19).

\begin{appendix}
\section*{Appendix}

\subsection*{A: Proofs of Section~\ref{sec:IID}}

\subsubsection*{Proof of Lemma~\ref{lem:sub}:} 
From the law of large numbers, we obtain that almost surely, 
\begin{equation}
\label{eq:urtzi}
\lim\limits_{r\rightarrow\infty} \frac{1}{r}  \tilde A_c(rt)  = \frac{\lambda}{\binom{K}{d}} t \textrm{  and  } \lim\limits_{r\rightarrow\infty} \frac{1}{r}\tilde S_c(s)ds = \mu t.
\end{equation}

The cumulative amount of capacity spent on serving type-$c$ jobs in server~$s$,  $T_{s,c}^{IID,r}(t)$ increases at rate $N^{IID}_c(t)/M^{IID}_{s}(t)\leq 1$. Hence,  $\frac{1}{r}T^{IID,r}_{s,c}(rt)-\frac{1}{r}T^{IID,r}_{s,c}(ru)\leq t-u$ for every $t\geq u$, i.e.,  $T_{s,c}^{IID,r}(rt)/r$ is Lipschitz continuous. Therefore, by the Arzela-Ascoli theorem we obtain that for almost all sample path $\omega$ and any sequence $r_{k}$, there exists a subsequence $r_{k_j}$ such that $\lim\limits_{j\rightarrow\infty} \frac{T^{IID,r_{k_j}}_c(r_{k_j}t)}{r_{k_j}}=\tau^{IID}_c(t),$ u.o.c.. Together with~\eqref{eq:frelationiid} and~\eqref{eq:urtzi}, we obtain Equation (\ref{eqfluid}).

\subsubsection*{Proof of Lemma~\ref{lemmaIID}:}

For ease of notation, we removed the superscript IID throughout the proof. 
Let $f(\vec n)= (f_c(\vec n),c\in\mathcal C)$, with $f_c(\vec n): \mathbb{R}_+^{|\mathcal C|} \to \mathbb{R}^{|\mathcal C|}$, denote  the drift vector field of $\vec N(t)$ when starting in state   $\vec N(0)= \vec n$, i.e., $f(\vec n)={d \over dt} \mathbb{E}^{\vec n}\big[ \vec N(t) \big] \Big|_{t=0}$.  
We can deduce from the results of~\cite[Proposition 5]{GG12} that the fluid limit $\vec n(t)$ satisfies 
\begin{equation}
\label{eq:fluid}
\frac{\mathrm{d} \vec n(t)}{\mathrm{d}t}\in F(\vec n(t)),
\end{equation}
where
\begin{equation}
F(\vec n):= conv\left( acc_{r\rightarrow\infty} \ f(r \vec n^r) \ \ with \ \ \lim\limits_{r\rightarrow\infty} \vec n^r = \vec n \right).
\label{eq:F}
\end{equation}
Here,  $ acc_{r\rightarrow\infty} \ x^r $ denotes the set of accumulation points of the sequence $x^r$ when $r$ goes to infinity and $conv(A)$ is the convex hull of set $A$.
An illustration of how $F$ is constructed  is
available in Figure~1 in \cite[Section 2]{GG12}.

Using Equations~\eqref{eq:fluid} and~\eqref{eq:F}, we  can  partly characterize the fluid process $m_s(t)=\sum_{c\in \mathcal C(s)} n_c(t)$. Denote by  $\tilde f_s(\vec n) = \sum_{c\in \mathcal C(s)} f_c(\vec n)$ the one-step drift of $M_s(t)$. The arrival rate of copies to server~$s$ equals $\lambda \binom{K-1}{d-1}/\binom{K}{d}=\lambda d/K$.
Recall from~\eqref{eq:tdr} that the  total departure rate of copies from server~$s$ equals $\mu
\left( \sum_{c\in \mathcal C(s)}\sum_{l\in c} \frac{n_c}{m_l} \right)$.
Hence, in state $\vec n$, the drift of server~$s$ is equal to
\begin{align}
\tilde f_{s}(\vec n) & = \lambda \frac{d}{K}-\mu
\left( \sum_{c\in \mathcal C(s)}\sum_{l\in c} \frac{n_c}{m_l} \right).  
\label{eq:fluidIID}
\end{align}

Let $G_2(\vec n) := \{s\in S : m_s\geq m_l, \forall l\}.$  Note that if $s\in G_2(\vec n)$, then $ \left( \sum_{c\in \mathcal C(s)}  \sum_{l\in c} \frac{n_c}{m_l} \right) \geq  \left( \sum_{c\in \mathcal C(s)}   \sum_{l\in c} \frac{n_c}{m_s} \right) = \sum_{c\in \mathcal C(s)} d\frac{n_c}{m_s} = d$. 

Now, let $\lim_{r\to\infty} \vec n^r=\vec n$, and $\vec n\neq \vec 0$. Then, for $s\in G_2(\vec n)$,    
$$
\lim_{r\to\infty} \tilde f_s(r \vec n^r) = \lambda d/K - \mu\left( \sum_{c\in \mathcal C(s)}  \sum_{l\in c} \frac{n_c}{m_l} \right)  \leq 
\lambda d/K-\mu d.
$$ 
Together with  \eqref{eq:fluid},~\eqref{eq:F} and
$\sum_{c\in \mathcal C(s)}\frac{\mathrm{d} n_c(t)}{\mathrm{d}t}=\frac{\mathrm{d} m_s(t)}{\mathrm{d}t}$, this concludes the proof.

\subsubsection*{Proof of Proposition~\ref{theorem1IID}:} 
Define $m_{max}^{IID}(t):=\max_{s\in S}\{ m_s^{IID}(t)\}$ and fix $T=m_{max}^{IID}(0)/d(\mu-\frac{\lambda}{K})$. 	From Lemma~\ref{lemmaIID}, we know that  at time $T$, $m_{max}(T)=0$.  
Since for any $s\in S$, $m_s^{IID}(t)\leq m_{max}^{IID}(t)$,  we conclude that at time $T$   the fluid system is empty. From~\cite[Corollary 9.8]{Robert03}, we then conclude that the process is ergodic.

\subsection*{B: Comments and proofs of Section~\ref{Sec:FCFSdes}}

\subsection*{Balance equations of the saturated system}
\label{Sec:ApFCFS}
For $\vec e_1, \vec e_2 \in \overline E$, we denote by $q(\vec e_1 ,\vec e_2)$ the transition probability from state $\vec e_1$ to state $\vec e_2$. 
Recall that in state $\vec e=(O_{\ell^*},\ldots,O_2,L_1,O_1) \in\bar E$, exactly $\ell^*(\vec e):=\ell^*$ jobs are being served, each of them with departure rate $\mu$. 
Hence, the balance equations of the saturated system are given by 
$$  \mu \ell^*(\vec e_1) \pi(\vec e_1 )=\sum_{\vec e_2\in \bar E } q(\vec e_2, \vec e_1) \pi(\vec  e_2),$$ 
with $\pi(\vec e)$ the steady-state distribution.

We will write down the balance equations in the case $d=K-2$. 
In that case, at any moment in time there 
is one job that is served in $d=K-2$ servers. In the remaining two servers, either one job, or two jobs are served. Hence, the states of the saturated system are of the form $\vec e=  (O_3, L_2, O_2, L_1, O_1)$ and $\vec e= (O_2, L_1, O_1)$. 
We denote by 
$\mathcal C(O_1):=\{ c\in\mathcal C :  c\cup O_1 = \{1,\ldots, K\} \}$  the subset of types that together with type $O_1$ make all servers busy. Hence, if the system is in state $\vec e=(c, L_1, O_1)$, $c\in \mathcal{C}(O_1)$, the total departure rate is $2\mu$.  
We denote by  $\bar{\mathcal C}(O_1) := \mathcal C -  O_1\cup \mathcal C(O_1) $ the subset of types that together with $O_1$ do not use all servers. For $O_1, O_2\in \mathcal C$, we denote by $\mathcal{C}(O_1, O_2):=  \{ c\in\mathcal{C}: c\cup O_1 \cup O_2 =\{1,\ldots, K\} \}$ the subset of types that together with $O_1$ and $O_2$ make all servers busy.

The balance equations are given by: 
\begin{eqnarray*} 
	2\mu\pi(O_2,L_1,O_1)&=& \mu\pi(O_2,L_1+1,O_1)+ \mu\sum_{c\in\mathcal C(O_1)} \sum_{j=0}^{L_1} (\frac{1}{\vert\mathcal C\vert})^{j+1}\pi(c,L_1-j,O_1)\\
	&&+\mu\sum_{c\in\mathcal C(O_1)}(\frac{1}{\vert\mathcal C\vert})^{L_1+1}\pi(O_1,0,c)
	+ \mu\sum_{c\in\bar{\mathcal{C}}(O_1)} \left(\frac{1}{\binom{K-1}{d}}\right)^{L_1} \pi(O_2,L_1,O_1,0,c)\\
	&&+\mu \sum_{c\in\bar{\mathcal C}(O_1)} \sum_{j=0}^{L_1}  \left(\frac{1}{\binom{K-1}{d}}\right)^{j} \pi(O_2,j, c,L_1-j,O_1),
\end{eqnarray*}
with $L_1\geq 0$, $O_1\in\mathcal C$, and $O_2\in\mathcal C(O_1)$. The term $(1/\binom{K-1}{d})^{j}$  in the fourth and fifth term on the right represents the probability that all $j$ waiting jobs are of type~$O_1$ (types $O_1$ and $c$ occupy $K-1$ servers, hence $\binom{K-1}{d}$ is the number of possible types that can compose $L_1$). 
For a $3\mu$ departure rate configuration state we have
$$
\begin{array}{l}
3\mu\pi(O_3,L_2,O_2,L_1,O_1)= \mu\pi(O_3,L_2,O_2,L_1+1,O_1)\\
+\mu\sum_{c\in\bar{\mathcal{C}}(O_1) \cap \mathcal{C}(O_1, O_3)}\sum_{j=0}^{L_1}\left(\frac{1}{\binom{K-1}{d}}\right)^{j+1} \pi(O_3,L_2+j+1,c,L_1-j,O_1)\\
+\mu\sum_{c\in {\mathcal{C}}(O_1,O_2)}
\sum_{j=0}^{L_2}\left(\frac{\binom{K-1}{d}}{\vert\mathcal C\vert}\right)^{j}\frac{1}{\vert\mathcal C\vert} (c,L_2-j,O_2,L_1,O_1)\\
+\mu\sum_{c\in \bar{\mathcal{C}}(O_1) \cap \mathcal{C}(O_1, O_3) }
\left(\frac{1}{\binom{K-1}{d}}\right)^{L_1+1}\pi(O_3,L_1+L_2+1,O_1,0,c) \\
+\mu \sum_{c\in\bar{\mathcal C}(O_1)\cap{\mathcal{C}}(O_1,O_2)}
\left(\frac{\binom{K-1}{d}}{\vert\mathcal C\vert}\right)^{L_2}\frac{1}{\vert\mathcal C\vert}\left(\sum_{j=0}^{L_1} \left(\frac{1}{\binom{K-1}{d}}\right)^{j}\pi(O_2,j, c,L_1-j,O_1)\right.\\
\left.
\quad + \left(\frac{1}{\binom{K-1}{d}}\right)^{L_1}\pi(O_2,L_1,O_1,0,c)\right)\\
+ \mu \sum_{c\in\mathcal C(O_1)}\left(\frac{\binom{K-1}{d}}{\mathcal\vert C\vert}\right)^{L_2}\left(\frac{1}{\vert\mathcal C\vert}\right)^2\left(\sum_{j=0}^{L_1} \left((\frac{1}{\vert\mathcal C\vert})^{j}\pi(c,L_1-j,O_1)\right)+(\frac{1}{\vert\mathcal C\vert})^{L_1} \pi(O_1,0,c)\right),
\end{array}$$  
with $O_1\in \mathcal{C}$, $O_2\in\bar{\mathcal{C}}(O_1)$, $O_3\in\mathcal C(O_1,O_2)$ and $L_1,L_2\geq0$.    
Note that on the right-hand-side, the term $\frac{\binom{K-1}{d}}{\mathcal\vert C\vert}$ is the probability that an arriving job is of type~$c\in O_1\cup O_2$ (the number of types $c$ with $c\in O_1\cup O_2$ is equal to $\binom{K-1}{d}$). 

\subsubsection*{Some properties of $\bar \ell$} 
In the proof, we will make use of the following properties for the saturated system (as defined in Definition~\ref{def:congested}). 
Recall that the average total departure rate of the saturated system is given by~$\bar \ell \mu$, where $\bar \ell$ is defined in~\eqref{eq:averagesFCFS} as the average number of jobs in service. 
For the saturated system, recall that $\ell^*(\vec e)$ denotes the number of jobs that is served in state~$\vec e$. Hence
\begin{equation}
\label{eq:barmu}
\lim_{t\rightarrow\infty}
\frac{1}{t}\int_0^{t} \ell^*(\vec e(u)) \textrm du=\bar \ell, \mbox{ \ almost surely}.
\end{equation} 
For the saturated system, let $\ell^*_c (\vec e)$ equal 1 if a type-$c$ job is served in state~$\vec e$ and 0 otherwise. 
Note that $\sum_{c\in \mathcal{C}} \ell_c^* (\vec e) = \ell^*(\vec e)$.  
Hence, using the system symmetry, (or more precisely, the exchangeability of the server contents),
and together with~\eqref{eq:barmu},  we obtain that 
\begin{equation}
\label{eq:sc}
\lim_{t\to\infty}\frac{1}{t}\int_0^{t}  \ell^*_c(\vec e(u)) \textrm du =\frac{\bar \ell }{\binom{K}{d}}, \ \mbox{almost surely.}
\end{equation}

\subsubsection*{Proof of Proposition~\ref{prop:props}}

We are given a saturated system with $K$ servers, with a central queue where jobs wait in order of arrival. The system starts serving at time~0. 
Let $c^K(i)$ denote the type of the $i$-th job at time~0 in this central queue. Let $\alpha^K_{is}(t)$ denote the attained service of this job at time~$t$ in server~$s\in c^K(i)$. Once the job~$i$ departs, the attained service $\alpha_{is}^K(t)$ is set equal to $\beta_i$, the service requirement of job~$i$. 
Let  $D_c^{K}(t)$ denote the number of departed type-$c$ jobs in the interval $(0,t]$ and $D_s^K(t):= \sum_{c\in \mathcal{C}^K(s)} D_c^{K}(t)$ the number of departed jobs from server~$s$, with $\mathcal{C}^K$ the set of types with $K$ servers. 
We  will prove that  
\begin{equation} 
\label{eq:D}
D_s^{K}(t)\geq_{st} D_s^{K-1}(t), \mbox{ with $s$ an arbitrary server in each of the systems.}
\end{equation}

Before proving this, we first show how \eqref{eq:D} implies that $\bar \ell/K$ is increasing in $K$, as stated in Proposition~\ref{prop:props}.
From~\eqref{eq:D} we have 
$\lim_{t\to\infty}\frac{1}{t}  D_s^{K}(t)\geq \lim_{t\to\infty}\frac{1}{t}  D_s^{K-1}(t)$, that is, the long-run departure rate from server~$s$ is increasing in the number of servers. Note that $\mu \sum_{c\in \mathcal{C}(s)}\ell^*_c(\vec e(t))$ is the instantaneous departure rate from server~$s$, where $\ell^*_c(\vec e(t))$ equals 1 if a type-$c$ job is served, and equals 0 otherwise. From~\eqref{eq:sc}, we have that the long-run departure rate from server~$s$ can equivalently be written as $\lim_{t\to\infty}  \frac{1}{t} \mu \sum_{c\in \mathcal{C}(s)}  \int_0^t\ell^*_c(\vec e(t)) \mathrm{d}u = \frac{\bar \ell \mu}{\binom{K}{d}}\binom{K-1}{d-1} = \frac{d \bar \ell \mu}{K}.$
Since the long-run departure rate is increasing in the number of servers, this implies that $\frac{\bar \ell}{K}$ is increasing in $K$ and proves the statement of Proposition~\ref{prop:props}. 

We are left with proving~\eqref{eq:D}. In order to do so,  we will couple the system with $K$ servers to a system with $K-1$ servers  as follows.  We consider the  central queue  associated to the saturated system with $K$ servers,   which corresponds to an infinite backlog of jobs (at time~0) ordered according to arrival (from $-\infty$). In the $K-1$ server model,   server~$K$ is removed. We couple the $K-1$ server model to the $K$ server model, by creating the central queue for the $K-1$ system as follows. 
For each $i$-th job in the central queue that has a copy in server~$K$, i.e., $K\in c^K(i)$, we choose uniformly at random another server among the remaining $K\backslash c^K(i)$ servers, denoted by $s^{K-1}(i)$. Hence, for any job with $K\in c^K(i)$, we set its type in the $K-1$ system as $c^{K-1}(i)= (c^K(i)\backslash K)\cup s^{K-1}(i)$. For all  jobs with $K\notin c^K(i)$, we set  $c^{K-1}(i)= c^K(i)$.
Below we show that   for all $t\geq 0$, 
\begin{equation} 
\alpha^K_{is}(t)\geq \alpha^{K-1}_{is}(t), \forall  i=1,\ldots \mbox{ and } \forall s\in c^K(i)\backslash K.
\label{eq:alpha1}
\end{equation}
and 
\begin{equation}
\alpha^K_{iK}(t)\geq \alpha^{K-1}_{is^{K-1}(i)}(t).
\label{eq:alpha2}
\end{equation}
From~\eqref{eq:alpha1} and \eqref{eq:alpha2} we obtain that \eqref{eq:D} holds: If a job~$i$ departs from a server~$s$ in the $K-1$ system, then (i) either also $s\in c^{K}(i)$, in which case this job   has  departed at a time $u\leq t$ in the $K$ system (from~\eqref{eq:alpha1}), (ii) or $s\neq c^{K}_i$, which implies that the type of the job is different in the  $K$ system and $K-1$ system, hence $s=s^{K-1}(i)$. Then,   from~\eqref{eq:alpha2} it follows that this job has  departed  at a time $u\leq t$  in the $K$ system. 
To conclude, in both cases, job~$i$ has already departed in the $K$ system before it departs in the $K-1$ system, hence,~\eqref{eq:D} holds. 

The result in~\eqref{eq:alpha1} and \eqref{eq:alpha2} will be proved by induction. It holds at time~0. Now assume that for all $u\leq t$ it holds that 
$\alpha^K_{is}(u)\geq \alpha^{K-1}_{is}(u), \forall  i=1,\ldots \mbox{ and } \forall s\in c^K(i)\backslash K.
$ 
and
$\alpha^K_{iK}(u)\geq \alpha^{K-1}_{is^{K-1}(i)}(u).
$
We prove that this remains true at time~$t^+$. 

In order for the inequality~\eqref{eq:alpha1} to no longer be valid at time~$t^+$, it  needs to hold that either \eqref{eq:alpha1} or \eqref{eq:alpha2} hold with strict equality. We first assume the first case, that is, $\alpha^K_{is}(t) = \alpha^{K-1}_{is}(t),$ for some $i$ and $s\in c^K(i)\backslash K.$
If $\alpha^K_{is}(t) = \alpha^{K-1}_{is}(t)=0$ and in the $K-1$ system it holds that $\alpha^{K-1}_{is}(t^+)>0$,  then this implies that one of the following occurs:
\begin{itemize}
	\item[(1)] in the $K$ system, the server~$s$ is serving the $i_1$-th job, with $i_1<i$, while in the $K-1$ system, server~$s$ starts serving job~$i$ at time $t^+$. However, since $\alpha^K_{i_1\tilde s}(t)\geq  \alpha^{K-1}_{i_1 \tilde s}(t),$ for all $\tilde s\in c^K(i)\backslash K$, this implies that job $i_1$ should not have a copy in server~$s$ in the $K-1$ system, since otherwise, job~$i_1$ was also still in service in the $K-1$ system.  However, due to the construction of the coupling and since $s\neq K$,   such a job $i_1$ does not exist.
	\item[(2)] in the $K-1$ system, this job~$i$ finishes its service in server~$\tilde s$, that is, $\alpha_{i\tilde s}^{K-1}(t^+)=\beta_i$ and hence $\alpha_{i  s}^{K-1}(t^+)=\beta_i$. But since \eqref{eq:alpha1} and \eqref{eq:alpha2} hold at time~$t$, this job is then also finished in the $K$ system, and hence also $\alpha^K_{is}(t^+)=\beta_i$. 
\end{itemize}

Now assume $\alpha^K_{is}(t) = \alpha^{K-1}_{is}(t)>0$ for some $i$ and $s\in c^K(i)\backslash K.$
Then, both jobs are in service in server~$s$, hence the inequality remains valid, unless the job departs in the $K-1$ system (and hence $\alpha^{K-1}_{is}(t^+)=\beta_i$), but not in the $K$ system. This can however not happen, since $\alpha^K_{j\tilde s}(t)\geq \alpha^{K-1}_{j\tilde s}(t), \forall j\tilde s\in c^K(j)\backslash K.
$ and $\alpha^K_{jK}(t)\geq \alpha^{K-1}_{js^{K-1}(j)}(t).
$
Hence, the inequality remains valid at time~$t^+$. 

To prove that    $\alpha^K_{is}(t) = \alpha^{K-1}_{is}(t),
$ implies $\alpha^K_{is}(t^+) = \alpha^{K-1}_{is}(t^+)$ follows exactly the same steps and is therefore left out.

\subsubsection*{Proof of Lemma~\ref{prop:FCFS_lb}.}
Both systems are coupled as follows: 
At time $t=0$, $N_c^{FCFS}(0)=0$ and $\tilde N_c^{(T)}(0)=\tilde A_c(T)$, where $\tilde A_c(t)$ is the arrival process of type-$c$ jobs.
During the time interval $[0,T]$, we couple the original system and its modified version by using the same arrivals and service times in the FCFS systems, as those that arrived in the $\tilde N^{(T)}$-system at time~0.

The result will be proved by induction. It holds at time $0$. Now assume that for all $u\leq t$ it holds that $\tilde N^{(T)}_c(u)\leq N^{FCFS}_c(u)+ (\tilde A_c(T)-\tilde A_c(u))^+$ and  $a_{cis}^{FCFS}(u)\leq a_{cis}^{\tilde N^{(T)}}(u)$, for all $i=1,\ldots,N_c^{FCFS}(t)$, $c\in \mathcal C$, $s\in S$. We prove that this remains true at time $t^+$. 

For that, assume there is a $c$ such that  $\tilde N^{(T)}_c(t) = N^{FCFS}_c(t)+ (\tilde A_c(T)-\tilde A_c(t))^+$. If $t<T$, only in the  FCFS system we can have an arrival, in which case $N^{FCFS}_c(t^+)= N^{FCFS}_c(t) +1$ and $(\tilde A_c(T)-\tilde A_c(t^+))^+ = (\tilde A_c(T)-\tilde A_c(t))^+ - 1$. Hence, the inequality remains valid. 
If $t\geq T$, then an arrival in the FCFS system is coupled to an arrival in the $\tilde N^{(T)}$-system, hence $\tilde N^{(T)}_c(t^+) = N^{FCFS}_c(t^+)$ (and note that $(\tilde A_c(T)-\tilde A_c(t))^+=0$). Now, assume the $i$-th type-$c$ job departs in the FCFS system (which can cause a violation of the inequality). Since $a_{ci\tilde s}^{FCFS}(t)\leq a_{ci\tilde s}^{\tilde N^{(T)}}(t)$, for all $\tilde s$, it holds that the same job departs in the $\tilde N^{(T)}$-system.  Hence, in all cases, the inequality $\tilde N^{(T)}_c(t^+) \leq N^{FCFS}_c(t^+)+(\tilde A_c(T)-\tilde A_c(t^+))^+$ remains valid at time~$t^+$. 

Now assume there exists  a $c,i,s$ such that $a_{cis}^{FCFS}(t)=a_{cis}^{\tilde N^{(T)}}(t)$.  
First assume $a_{cis}^{FCFS}(t)=a_{cis}^{\tilde N^{(T)}}(t)>0$.   Because of FCFS, in both systems this copy has entered service in server~$s$ at the same instant of time. Hence, it cannot happen that  $a_{cis}^{FCFS}(t^+)>a_{cis}^{\tilde N^{(T)}}(t^+)$.
If instead $a_{cis}^{FCFS}(t)=a_{cis}^{\tilde N^{(T)}}(t)=0$, the $i$-th type-$c$ copy in server~$s$ is waiting in the queue in both systems. We need to prove that if this copy would enter  service in server~$s$ at time $t^+$ in the FCFS system, it also enters service in the $\tilde N^{(T)}$-system in server~$s$.  From the FCFS discipline and $a_{\tilde cjs}^{FCFS}(t)\leq a_{\tilde cjs}^{\tilde N^{(T)}}(t)$, for all $\tilde c, j\leq i$, this follows directly. 

\subsubsection*{Proof of   Proposition~\ref{prop:FCFS_n}}

We will prove that if the $\tilde N^{(T)}$-system is stable and $T$ is sufficiently large, then $\lambda\leq \bar \ell \mu$. From Lemma~\ref{prop:FCFS_lb}, it follows that stability of the FCFS system, implies stability of the $\tilde N^{(T)}$-system, and hence $\lambda\leq \bar \ell \mu$, which would conclude the proof.

We will now prove that if the $\tilde N^{(T)}$-system is stable, then $\lambda \leq \bar \ell \mu$.  We  define  the random variable $\tau(T)$ as the first moment in  time  a servers gets empty in the $\tilde N^{(T)}$-system, i.e., 
$\tau(T):= \min\{u: \tilde M_s^{(T)}(u)=0, \mbox{for some server~$s$}\}$.  
Up till time $\tau(T)$, the $\tilde N^{(T)}$-system is stochastically equivalent to the saturated system. 
Hence, using the Markovian description of the process $M_s^{\tilde N^{(T)}}$ and  Dynkin's formula, we have that there exists a martingale $(Z_s(t))_{t \ge 0}$ such that
\begin{equation}
\label{eq:ast}
\frac{M^{\tilde N^{(T)}}_s(\tau(T))}{\tau(T)} = \frac{d\lambda}{K}\frac{ (T+(\tau(T)-T)^+)}{\tau(T)} -  \frac{1}{\tau(T)}\int_0^{\tau(T)} \mu \sum_{c\in C(s)}\ell_c(\vec e(u)) \textrm du + \frac{Z_s(\tau(T))}{\tau(T)},
\end{equation}
Since the increasing process associated to $Z_s$ is  bounded in mean by  $C t$, with $C>0$,  
it follows that $\sup_t E\big( {Z_s(t) \over t} \big)^2 \le {C t \over t^2}$, which in turn implies that $ \frac{Z_s(\tau(T))}{\tau(T)} \to 0$ in $L^2$ and hence the convergence holds almost surely.

Since by the law of large numbers for the Poisson process, $\liminf_{T} {\tau(T) \over T} \ge c >0,$ almost surely, together with~\eqref{eq:sc}, it follows that 
\begin{equation}
\label{eq:final}
\lim_{T\to \infty} \frac{M_s^{\tilde N^{(T)}}(\tau(T))}{\tau(T)} = \frac{d\lambda}{K} \frac{ T+(\tau(T)-T)^+}{\tau(T)} - \frac{d}{K}\bar \ell \mu, \ \mbox{ almost surely. }
\end{equation}

By assumption, the $\tilde N^{(T)}$-system is stable, hence,   $\mathbb{E}(\tau(T))<\infty$.    We can therefore focus on a sample path realization~$\omega$ such that $\tau(T)<\infty$.
From~\eqref{eq:final}, it follows that for each $\epsilon>0$, we can find $T$ such that    $\frac{M^{\tilde N^{(T)}}_s(\tau(T))}{\tau(T)} \geq \frac{d\lambda}{K}\frac{ T+(\tau(T)-T)^+}{\tau(T)} - \frac{d}{K}\bar \ell \mu   -\epsilon.$
Let  server~$\tilde s$ be such that  $M_{\tilde s}^{\tilde N^{(T)}}(\tau(T))=0$. If $\tau(T)\le T$, then, $0=\frac{M_{\tilde s}^{\tilde N^{(T)}}(\tau(T))}{\tau(T)} \geq \frac{d}{K}(\frac{\lambda T}{\tau(T)} - \bar \ell \mu )-\epsilon$. Since  $\tau(T)\le T$, this implies $\lambda\leq \bar\ell \mu+\frac{K}{d}\epsilon$. On the other hand, if $\tau(T)>T$, then, $0=\frac{M_{\tilde s}^{\tilde N^{(T)}}(\tau(T))}{\tau(T)} \geq \frac{d}{K}(\lambda -\bar \ell \mu )-\epsilon,$ i.e., $\lambda \leq  \bar \ell\mu + \frac{K}{d}\epsilon$. Since this holds for any $\epsilon$,  we conclude that $\lambda\leq \bar \ell \mu$.

\subsubsection*{Proof of Lemma~\ref{prop:FCFS_up}}
We couple both systems as follows: at time zero, both systems start in the same initial state, that is, $\hat N_c(0)=N_c^{FCFS}(0)$ and  $a_{cis}^{\hat N}(0) = a^{FCFS}_{cis}(0)$, for all $c,i,s$. Arrivals and their service requirements are coupled. 


The result will be proved by induction. It holds at time 0. Now assume that for all $u\leq t$ it holds that  
$\hat N_c(u)\geq N_c^{FCFS}(u)$ 
and  $a_{cis}^{\hat N}(u)\leq a_{cis}^{FCFS}(u)$ for all $i=1,\ldots,N_c^{FCFS}(t)$, $c\in \mathcal C$, $s\in S$. Below, we prove that this remains true at time~$t^+$. 

Assume there is a $c$ such that $\hat N_c(t)=N_c^{FCFS}(t)$. The inequality can be violated at time $t^+$ if there is a type-$c$ departure in $\hat N$, but not in FCFS. However, note that if the head-of-the-line type-$c$ job in the $\hat N$-system departs, since $a_{c1\tilde s}^{\hat N}(t)\leq a_{c1\tilde s}^{FCFS}(t)$, this job would also depart from the FCFS system. Hence,  $\hat N_c(t^+)= N_c^{FCFS}(t^+)$  at time~$t^+$. 

Now assume there exists a $c,i,s$ such that $a_{cis}^{\hat N}(t)=a_{cis}^{FCFS}(t)$.
First assume $a_{cis}^{\hat N}(t)=a_{cis}^{FCFS}(t)>0$.   Because of FCFS, in both systems this copy has entered service in server~$s$ at the same instant of time. Hence, it cannot happen that  $a_{cis}^{\hat N}(t^+)>a_{cis}^{FCFS}(t^+)$.
If instead $a_{cis}^{\hat N}(t)=a_{cis}^{FCFS}(t)=0$, the $i$-th type-$c$ copy in server~$s$ is waiting in the queue in both systems. We need to prove that if this copy would enter  service in server~$s$ at time $t^+$ in the $\hat N$-system, it also enters service in the FCFS-system in server~$s$.  From the FCFS discipline and $a_{\tilde cjs}^{\hat N}(t)\leq a_{\tilde cjs}^{FCFS}(t)$, for all $\tilde c, j\leq i$, this follows directly. 

\subsubsection*{Proof of Proposition~\ref{theo:FCFS_s}}

From Lemma~\ref{prop:FCFS_up} we have that $\hat N(t)$ is an upper bound for  the original FCFS system. Hence, it will be enough to prove stability of the process $\hat N(t)$. 

In order to prove stability of $\hat N(t)$, we study the fluid-scaled system. That is, for each $r$, we study $\hat N^r(t)$,  with $\hat N^r(0) =r n(0)$. 
Define $T_0=\frac{|n(0)|}{ \mu}$. 
By definition of $\hat N^r(t)$,   in the interval $[0,rT_0]$, only those jobs present at time 0 are served (according to FCFS). From time $rT_0=\frac{|N(0)|}{ \mu}$  onwards, all jobs can be served. 

We write $\hat N^r(t)= \hat N^r_A(t)+ \hat N^r_B(t)$, where $\hat N^r_A(t)$ denotes the number of old  jobs, that is, the number of jobs present at time $t$ among those that were already present at time $t=0$.
We let  $\hat N^r_B(t)=\hat N^r(t) - \hat N^r_A(t)$  denote the number of new  jobs present at time $t$. Similarly, we let  $\hat M_{s,B}(t)$  denote the number of new jobs that have a copy in server~$s$. 

We now show that for any fluid limit $\hat n(t)$ of $\hat N^r(r t)$ it holds that it is zero at some time smaller than or equal to $T_1:=T_0 \frac{\lambda}{\lambda-\bar \ell \mu}=|n(0)| \frac{\lambda}{\mu (\lambda-\bar \ell \mu)}$, that is $|\hat n(\tilde T_1)|=0$.

In the interval $[0,r T_0]$, the system serves only the jobs present at time $0$.  
Let $\hat\ell^r_A(t)$ denote the number of such jobs in service at time $t$ in the $\hat N^r$-system. 
Hence, using the Markovian description of the process $\vert \hat N^r_A(r t)\vert$ and  Dynkin's formula, we have that there exists a martingale $(Z(t))_{t \ge 0}$ such that
\begin{equation}
\label{eq:ast}
\frac{\vert \hat N^r_A(r t)\vert}{r} = \vert\vec n(0)\vert   -\mu \frac{1}{r}\int_0^{r t}\hat\ell_A(u)\mathrm{d}u
+ \frac{Z(r t)}{r},
\end{equation}
for $t\in [0,T_0]$. 
Since the increasing process associated to $Z$ is  bounded in mean by  $C t$, with $C>0$,  
it follows that $\sup_t E\big( {Z(r t) \over r} \big)^2 < \infty$, which in turn implies that $ \frac{Z(r t)}{r} \to 0$ almost surely.		
Further note that $\hat\ell_A(u)\geq 1$ whenever $\vec n\neq \vec 0$. Together, this gives that 
$$\lim_{r\to \infty} \frac{\vert \hat N^r_A(r t)\vert }{r} =\max(0, \vert\vec n(0)\vert   -\mu t), \mbox{for } t\in [0,T_0].$$  
Hence, for any fluid limit $\hat n_A(t)$, we have  $ \vert \hat n_{A}(t)\vert \leq  \max(0,\vert n(0)\vert - \mu t),$
so that  $|\hat n_A(T_0)|=0$.

In the remainder of the proof, we study fluid limits of the process $\hat N_B(rt)/r$.
We define the random variable $T_1^r:=\inf\{t>T_0: \mbox{there is an } s \mbox{ s.t.} \hat M_{s,B}(r(T_0+t))=0 \}$ as the first moment after time~$T_0$ that one of the servers gets empty. 
By the law of large numbers, $\liminf_{r} {T_1^r} \ge c >0$, almost surely. Hence,
without loss of generality, we can focus on sample paths, such that the latter is the case.
For a given sample path, let $$\tilde T_1 :=\liminf_{r\to\infty} T_1^r.$$  Note that  $\tilde T_1 \ge c.$
We consider henceforth the subsequence $r_j$ of any given sequence $r$, such that 
\begin{equation*}
T_1^{r_j}>\tilde T_1 -\epsilon, \ \forall r_j.
\end{equation*} 
In particular, this implies that all servers are working on new jobs during the interval $[r_j T_0,r_j T_0+r_j (\tilde T_1-\epsilon)]$, for any $r_j$. Also note that all jobs $\hat N_B(T_1)$ are ``freshly'' sampled, and hence the system behaves as a saturated system during this time frame.

Using the Markovian description of the process $ {\hat M}^r_B(r t)$ and  Dynkin's formula, we have that there exists a martingale $(Z_s(t))_{t \ge 0}$ such that 
\begin{eqnarray}
\label{eq:3o}
&&\lim_{j\to \infty} \frac{\vert \hat M^{r_j}_{s,B}(r_j (T_0+ \tilde T_1-\epsilon))\vert }{r_j}\nonumber\\
&&\quad =\lambda \frac{d}{K} (T_0+\tilde T_1-\epsilon)  - 
\frac{\mu}{r_j} \int_{r_jT_0}^{r_j(T_0+\tilde T_1-\epsilon)} \sum_{c\in \mathcal{C}(s)}\ell_c^*(\vec e(u))\mathrm{d}u + \frac{Z_s(r_j (T_0+ \tilde T_1-\epsilon))}{r_j},
\end{eqnarray}
where we recall that $\ell_c^*(\vec e(u))$ equals 1 if a type-$c$ job is in service, and equals zero otherwise.
Since the increasing process associated to $Z_s$ is  bounded in mean by  $C t$, with $C>0$,    
it follows that $\sup_t E\big( {Z_s(r t) \over r} \big)^2 \le {C t r \over r^2}= {Ct \over r}$, 
which in turn implies that $ \frac{Z_s(r t)}{r} \to 0$ almost surely.		
Now together with~\eqref{eq:sc}, 
we conclude that for this sample path 
\begin{eqnarray*}
	\lim_{j\to \infty} \frac{\vert \hat M^{r_j}_{s,B}(r_j (T_0+ \tilde T_1-\epsilon))\vert }{r_j} 
	=\lambda \frac{d}{K} (T_0+\tilde T_1-\epsilon)  - \bar \ell \mu \frac{d}{K} (\tilde T_1-\epsilon)= \lambda \frac{d}{K} T_0+ (\lambda- \bar \ell \mu)\frac{d}{K} (\tilde T_1-\epsilon).
\end{eqnarray*}

Hence,   the corresponding fluid limit satisfies
\begin{equation}
\label{eq:msb}
\hat m_{s,B}(T_0+\tilde T_1-\epsilon)=\lambda \frac{d}{K} T_0+ (\lambda- \bar \ell \mu)\frac{d}{K} (\tilde T_1-\epsilon).
\end{equation}

In case $\tilde T_1>T_1$, we set $\epsilon=\tilde T_1-T_1$, and since $T_1=T_0 \frac{\lambda}{\lambda-\bar \ell \mu},$ one obtains from~\eqref{eq:msb} that $\hat m_{s,B}(T_1)=0$ for all $s$. 
Now assume $\tilde T_1\leq T_1$. 
By definition of $T_1^r$, it holds that $\prod_{s}\hat M^r_{s,B}(r(T_0+ T_1^r)) =0$ and hence $\prod_{s}\hat m_{s,B}(T_0+\tilde T_1) =0$. From~\eqref{eq:msb} one has $\hat m_{s,B}(T_0+ \tilde T_1-\epsilon) = \hat m_{\tilde s, B}(T_0+ \tilde T_1-\epsilon)$, for any $s,\tilde s$ and any $\epsilon>0$. This, together with the fact that a fluid limit $\hat n_B(t)$ is a continuous function and $\prod_{s}\hat m_{s,B}(T_0+ \tilde T_1) =0$, it follows that $\hat m_{s,B}( T_0+\tilde T_1) =0$ for all $s$. 

We conclude that at time $T_0+\tilde T_1$, for any fluid limit $\hat n(\cdot )$ of $\hat N^r(\cdot)$, it holds that $\hat n(T_0+\tilde T_1)=0$. 
From~\cite[Corollary 9.8]{Robert03}, we conclude that the process $\hat N(t)$ is ergodic.


\subsection*{C: Proofs of Section~\ref{Sec:Des}}

\subsubsection*{Proof of Lemma \ref{prop1}:} 

We couple the two systems as follows: at time zero, start in the same initial state. Arrivals are coupled in both systems. Below it will become clear how the departures are coupled under both systems. 

Assume that at time $t\geq0$, $N^{PS}_c(t)\geq N_c^{LB}(t)$ for all $c\in \mathcal C$. We prove that this remains valid at time $t^+$. We only need to analyse states such that $N^{PS}_c(t)=N^{LB}_c(t)=n_c$, for some $c\in \mathcal C$. Under this situation, note that $M^{PS}_{s_{ci}^*(t)}(t) \geq M^{PS}_{s_{c}^{min}(\vec N^{PS}(t))}(t)\geq M^{LB}_{s_{c}^{min}(\vec N^{PS}(t))}(t)\geq M^{LB}_{s_{c}^{min}(\vec N^{LB}(t))}(t)$ for all $i=1,\ldots, N^{PS}_c(t)$. 
Hence, the departure rate of type-$c$ jobs in the PS system, $\mu \sum_{i=1}^{N^{PS}_c(t)} \ \frac{1}{M^{PS}_{s_{ci}^*(t)}(t)}$, is smaller than or equal to that in the LB-system, $\mu \ \frac{N^{LB}_c(t)}{M^{LB}_{s^{min}_{c}(\vec N^{LB}(t))}(t)}$.  
We can therefore couple the systems such that if there is a type-$c$ departure in the original PS model, then also a type-$c$ departure occurs in the LB-system. 
Since arrivals are coupled in both systems, it follows directly that at time $t^+$, 
$N^{PS}_c(t^+)\geq N_c^{LB}(t^+)$.

\subsubsection*{Proof of Lemma \ref{lemma1}:}

For simplicity in notation, we remove the superscript $LB$ throughout the proof. From~\eqref{eq:depLB}, we have that the departure rate of $M_s(t)$  is given by
\begin{equation}\label{eq5}
\displaystyle \sum_{c\in \mathcal C(s)}  \frac{N_c(t)}{M_{s_{c}^{min}(\vec n)}(t)}.\end{equation}
Recall that  $c\in \mathcal C_l^s(\vec n)$ if server $l$ is the server with the minimum number of copies that serves a type-$c$ job.  Hence, if $c\in \mathcal C_l^s(\vec n)$, then $s^{min}_c(\vec n)=l$. Since $\mathcal C(s)=\cup_{l\in D_s(\vec n)}\mathcal C_l^s(\vec n)$, Equation~\eqref{eq5} can be written as
\begin{eqnarray}\label{eq4}
&&\sum_{l\in D_s(\vec n)} \frac{ \sum_{c\in \mathcal C_l^s(\vec n)} N_c(t)}{M_{l}(t)}.\end{eqnarray}
Using that $\sum_{c\in \mathcal C_s^s(\vec n)} N_c(t)$ can equivalently be written as $M_{s}(t) - \sum_{l\in D_s(\vec n), l\neq s}\sum_{c\in \mathcal C_l^s(\vec n)} N_c(t)$, we obtain that (\ref{eq4}) is equal to 
\begin{eqnarray}
\sum_{l\in D_s(\vec n), l\neq s} \frac{\sum_{c\in \mathcal C_l^s(\vec n)}N_c(t)}{M_{l}(t)}+1- \sum_{l\in D_s(\vec n), l\neq s} \frac{\sum_{c\in \mathcal C_l^s(\vec n)} N_c(t)}{M_{s}(t)} 	 =1+ \sum_{l\in D_s(\vec n)} \frac{(M_s(t)-M_l(t)) \sum_{c\in \mathcal C_l^s(\vec n)}N_c(t)}{M_s(t)M_{l}(t)}.\nonumber\end{eqnarray}

\subsubsection*{Proof of Lemma \ref{lemma2}:}
For ease of notation, we removed the superscript PS throughout the proof. 

Let $f(\vec n)= (f_c(\vec n),c\in\mathcal C)$, with $f_c(\vec n): \mathbb{R}_+^{|\mathcal C|} \to \mathbb{R}^{|\mathcal C|}$, denote  the drift vector field of $\vec N(t)$ when starting in state $\vec N(0)= \vec n$, i.e. $f(\vec n)={d \over dt} \mathbb{E}^{\vec n}\big[ \vec N(t) \big]_{t=0}$. Recall that the fluid limit can be characterized as in
Equations \eqref{eq:fluid} and \eqref{eq:F}. We want to partly characterize the fluid process $m_s(t)=\sum_{c\in \mathcal C(s)} n_c(t)$. We denote by  $\tilde f_s(\vec n) = \sum_{c\in \mathcal C(s)} f_c(\vec n)$ the drift of $M_s(t)$. 

From~Lemma~\ref{lemma1}, we can write the drift of $M_s(\cdot)$, starting in state $\vec N(0)=\vec n$, as
\begin{align}
\tilde f_{s}(\vec n)	& = \lambda \frac{d}{K}-\mu \mathbf{1}_{(m_s>0)}
- \mu\mathbf{1}_{(m_s>0)}\left(\sum_{l\in D_s(\vec n)} \frac{(m_{s}-m_{l})\sum_{c\in \mathcal C_{l}^s(\vec n)}n_c}{m_{s}m_{l}}
\right),\label{eq:f}
\end{align}  
where  $D_s(\vec n)=\{l\in S: m_s\geq m_l\}$ is the set of servers that have less than or equal number of copies, compared to server~$s$, in state~$\vec n$.

Let $G_1(\vec n) := \{s\in S : m_s\leq m_l, \forall l\}$. If $s\in G_1(\vec n)$ and $\lim_{r\to\infty } \sum_{c\in \mathcal C(s)} n^r_c = \lim_{r\to\infty } m^r_s>0$, it follows from~\eqref{eq:f} that $$\lim_{r\to\infty} \tilde f_s(r \vec n^r) = \lambda d/K-\mu.$$ 
If instead $s\in G_1(\vec n)$ and  $\lim_{r\to\infty } m^r_s=0$, then 
\begin{eqnarray*}
	&&conv\left( acc_{r\rightarrow\infty} \ \tilde f_s(r \vec n^r) \ \ with \ \ \lim\limits_{r\rightarrow\infty} \vec n^r = \vec n \right)= conv(\lambda d/K-\mu, \lambda d/K ).
\end{eqnarray*}
Combining~\eqref{eq:fluid} and $\sum_{c\in \mathcal C(s)}\frac{\mathrm{d} n_c(t)}{\mathrm{d}t}=\frac{\mathrm{d} m_s(t)}{\mathrm{d}t}$, we conclude the proof. 

\subsubsection*{Proof of Proposition \ref{theorem1}:} 
From Lemma~\ref{prop1} we have that if the lower-bound system is unstable, then also the original PS system. Hence, to prove Proposition~\ref{theorem1}, it will be enough to prove  that $\vec N^{LB}(t)$ is unstable if $\rho>1/d$. This is done in the remainder of the proof. 

For ease of notation, we remove the superscript $LB$ throughout the proof. 
To prove that the system is transient, below we will show that there is  a subsequence of $t$ such that the system $\vec N(t)$ converges towards $+\infty$. 

Define $m_{min}(t):= \min_{s\in S} \{m_s(t)\}$ and fix $T= (|\vec n|+\delta)/(\lambda d/K-\mu)$, for some $\delta>0$.
From Lemma~\ref{lemma2}, we know that at time $T$, $m_{min}(T) \geq |\vec n | +\delta$, when $\vec n(0)=\vec n$. Hence, as well, 
\begin{equation}
\label{eq:dd}
\vert  \vec{n}(T)\vert \geq m_{min}(T) \geq |\vec n| +\delta.
\end{equation}

For almost all sample paths, and  any subsequence $r_k$ of $r$, there exists a further subsequence $r_{k_j}$ 
such that $\lim_{j\to\infty} \frac{|\vec N^{r_{k_j}}(r_{k_j} T)|}{r_{k_j}} = |\vec n(T)|\geq |\vec n|+\delta$, with $\vec N^r(0)=r \vec n$ and $\vec n(t)$ a fluid limit (the inequality follows from~\eqref{eq:dd}).
Hence, when considering the liminf subsequence, this gives, for all $\vec n$,
$$\liminf_{r\to\infty} \Big| \frac{\vec N^{r} (rT)}{r} \Big| \ge |\vec{n} | + \delta,$$
where $\vec N^r(0)=r \vec n$. 
From Fatou's lemma, this implies
$$ \liminf_{r\to\infty} \mathbb{E} \Big| \frac{\vec N^{r} (rT)}{r} \Big|  \ge |\vec {n} | + \delta.$$ 
Hence, there exists $r_0(\vec {n}) \ge 1$, such that 
$\vec N^{r_0(\vec n)}(0)= r_0(\vec n) \vec n$ and 
\begin{equation}\mathbb{E} \Big| \vec N^{r_0(\vec {n}) } (r_0(\vec n) T)  \Big|  \ge r_0(\vec {n}) \left(|\vec {n}| + \delta - \epsilon\right),
\label{eq:fld}
\end{equation}
for some $\epsilon$, with $0<\epsilon<\delta$.  
Now, for any $\vec n$,  define the discrete time stochastic process $(\vec Y_l, \vec Z_l)$, $l\geq 0$:
\begin{align*}
&\vec {Z}_0=\vec n,\\
& \vec {Y}_{l+1}= \vec {N}^{ r_0(\vec Z_l)}(r_0(\vec {Z_l} ) T ), \mbox{ where } \vec N^{r_0(\vec Z_l)}(0)= r_0(\vec Z_l)\vec Z_l,   \\
& \vec {Z}_{l+1} = {\vec {Y}_{l+1} \vec r_0(\vec {Z_l})}, l \ge 0.
\end{align*}

Observe that:
\begin{enumerate}
	\item $(\vec {Y_l},\vec {Z_l})$ is Markov, since $\vec {N}$ is a Markov process.
	\item It follows from~\eqref{eq:fld} that $\mathbb{E}\Big(|\vec {Z}_{l+1}| \Big| \vec {Z}_l \Big) - |\vec {Z}_l| \ge  \delta- \epsilon>0$, $l\geq 0$.
	\item  
	Using Dynkins formula for the continuous time process $\vec N(t)$, we see that
	\begin{align*}
	&\mathbb{E}\Big( \mathbb{E} |\vec {Z}_{l+1}|- |\vec {Z}_l| \Big)
	= \mathbb{E}\Big( |\vec{Z}_l|   +  {1 \over r_0(\vec{Z}_l)} \int_0^{r_0(\vec{Z}_l)T} a(\vec{N}_s) ds - |\vec {Z}_l| \Big)=  \mathbb{E}\Big( {1 \over r_0(\vec{Z}_l)} \int_0^{r_0(\vec{Z}_l)T} a(\vec{N}_s) ds \Big),
	\end{align*}
	where $a(\cdot)$ is the drift of the norm function.
	Note that given the model (bounded rates of arrival and departures) this drift is a bounded function (say by $\gamma$), which implies that
	$$\mathbb{E}\Big| |\vec {Z}_{l+1}|- |\vec {Z}_l| \Big| \le \mathbb{E}\Big( \mathbb{E} \Big({\gamma r_0(\vec{Z}_l) T \over r_0(\vec{Z}_l)}\Big| \vec{Z}_l \Big)\Big)= \gamma T < \infty.$$

\end{enumerate}
Using a classical transience criterion for Markov chains, (see for instance Proposition 8.9 in \cite{Robert03}), we obtain that $Z_l$ is transient. This in turn directly implies that $\vec N(t)$ converges along one subsequence of $t$ towards $+\infty$, which implies that it is transient.
\subsubsection*{Proof of Lemma \ref{lemmaup}:}

We couple both systems as follows: at time zero, we start in the same initial state. Arrivals and their service requirements are coupled in both systems. 

The result will be proved by induction. It holds at time~0. Now assume that for all $u\leq t$ it holds that 	  $\alpha^{UB}_{i,s}(u)\leq \alpha^{PS}_{i,s}(u)$	 for all $i=1, \ldots,$ and  $s\in S$. Below we prove that this remains true at time~$t^+$.

Let $c(i)$ denote the type of the $i$-th arrived job and let $\tilde A(t)$ denote  the number of arrivals until time~$t$.
Assume there is an $i\leq \tilde A(t)$ and $s\in c(i)$ such that $\alpha^{UB}_{i,s}(t)= \alpha^{PS}_{i,s}(t)$. 
Note that for all $c$, 
$$N_c^{UB}(t)=\sum_{j=1}^{\tilde A(t)} \mathbf{1}_{\{c = c(j)\}}\mathbf{1}_{\{\exists \tilde s\in c(j), \mbox{ s.t. } \alpha_{j,\tilde s}^{UB}(t)< \beta_j\}} \ \mbox{ and } \ N_c^{PS}(t)=\sum_{j=1}^{\tilde A(t)} \mathbf{1}_{\{c=c(j)\}} \mathbf{1}_{\{\forall \tilde s\in c(j), \alpha_{j,\tilde s}^{PS}(t)< \beta_j\}}.$$
Since $\alpha^{UB}_{j,\tilde s}(t)\leq \alpha^{PS}_{j,\tilde s}(t)$, for all $j, \tilde s$, it follows that  $N^{UB}_{c}(t) \geq N^{PS}_{c}(t)$, for all $c$,  hence   $M^{UB}_{\tilde s}(t) \geq M^{PS}_{\tilde s}(t)$ for all servers $\tilde s$. In particular, this implies that   $\frac{\mathrm{d} \alpha^{UB}_{i,s}(t)}{\mathrm{d}t} = \frac{1}{M^{UB}_s(t)}\leq \frac{1}{M^{PS}_s(t)} = \frac{\mathrm{d}a^{PS}_{i,s}(t)}{\mathrm{d}t}$, for $s\in c(i)$, which together with $\alpha^{UB}_{i,s}(t)= \alpha^{PS}_{i,s}(t)$ gives  that $\alpha^{UB}_{i,s}(t^+)\leq  \alpha^{PS}_{i,s}(t^+)$ holds at time~$t^+$.


\subsection*{D: Proofs of Section~\ref{sec:ROS}}

\subsubsection*{Proof of Lemma~\ref{lemmaROS}:}

For ease of notation, we remove the superscript $ROS$ throughout the proof. 

Assume at time~0 we are in state $\vec N(0)=\vec N$. We first will write the probability that a given server~$s$ is serving a copy that is not in service in any other server. We denote this probability by $P_s(\vec N)$. 
In order to derive that, we   consider $P_s(\vec N|c)$ defined as the probability that server~$s$ is serving a type-$c$ job, $s\in c$, and this job is not in service in any other server. 

Let~$-\tilde T_s<0$ denote the time that server~$s$ started working on the copy which it is serving at time~0. 
When the server becomes
idle, it chooses a copy uniformly at random. 
Hence, the probability that a copy from a type-$c$ job is being served in server~$s$ is  given by
$\frac{N_c(-\tilde T_s)}{M_s(-\tilde T_s)}$. 
Using the law of total probability, we have
\begin{eqnarray}
\label{mo2}
P_s(\vec N)= \sum_{c\in \mathcal{C}(s)}  \frac{N_c(-\tilde T^r_s)}{M_s(-\tilde T^r_s)} P_s(\vec N|c).
\end{eqnarray}
To calculate $P_s(\vec N|c)$, note that $\frac{N_c(-\tilde T^r_l)-1}{M_l(-\tilde T^r_l)}$ is the probability that server~$l$ is \emph{not} serving the type-$c$ copy that is now in service in server~$s$, with $l,s\in c$. Hence,
\begin{eqnarray}
\label{mo}
P_s(\vec N|c) = \Pi_{l\in c, l\neq s} \frac{M_l(-\tilde T^r_l)-1}{M_l(-\tilde T^r_l)}, \  s\in c.
\end{eqnarray}

We now characterize the fluid limits, which we recall can be characterized as in Equations~\eqref{eq:fluid} and~\eqref{eq:F}. 
Let $f(\vec n)= (f_c(\vec n) , c\in \mathcal{C})$, with $f_c(\vec n) : \mathbb{R}_+^{|\mathcal C|} \to \mathbb{R}^{|\mathcal C|}$, denote  the drift vector field of $\vec N(t)$ when starting in state $\vec N(0)= \vec n$, i.e., $f(\vec n)={d \over dt} \mathbb{E}^{\vec n}\big[ \vec N(t) \big]_{t=0}$. 
Hence, we study the fluid drift in points $r \vec n^r$, where $\lim_{r\to\infty} \vec n^r=\vec n$. 
That is,      $\vec N(0)=r\vec n^r$. 

Since the transition rates $\mu$ and $\lambda$ are of order $O(1)$, it follows directly that $\tilde T_s^r$ and $\vec N(-\tilde T_s^r)- \vec N(0)$ are of order $O(1)$ as well, so that   
\begin{equation}
\label{eq:roschoose}
\lim_{r\to\infty}\frac{N_c(-\tilde T_s^r)}{M_s(-\tilde T_s^r)}= \lim_{r\to\infty}\frac{N_c(0)}{M_s(0)}= \frac{n_c(0)}{m_s(0)} \mbox{ and } \lim_{r\to\infty}\frac{M_l(-\tilde T_l^r)-1}{M_l(-\tilde T_l^r)}=1.
\end{equation}
It hence follows  from~\eqref{mo2} and~\eqref{mo} that 
\begin{equation}
\label{eq:1}
\lim_{r\to\infty} P_s(r\vec n^r) =  1.
\end{equation}

We denote by  $\tilde f_s(\vec n) = \sum_{c\in \mathcal C(s)} f_c(\vec n)$ the one-step drift of $M_s(t)$.
When starting in state $\vec N(0)=r \vec n^r$, the latter is in the limit equal to 
\begin{align}  
&\lim_{r\to\infty}	\tilde f_{s}( r \vec n^r)
=  \lambda \frac{d}{K} -\mu
\left( \sum_{c\in \mathcal C(s)}\sum_{l\in c}  \frac{ n_c}{ m_l}\lim_{r\to\infty} \left( P_l(r \vec n)\right)  -  \lim_{r\to\infty} (g_{c,l,s}( r \vec n^r) (1- P_l(r \vec n)))\right)\label{eq:new}
\end{align}  
with $g_{c,l,s}=O(1)$.   
Note that the first term multiplied by $\mu$ in~\eqref{eq:new} represents departures  of type-$c$ jobs, $c\in \mathcal{C}(s)$,   who were served in one unique server. Here $\frac{ n_c}{ m_l}$ represents the probability 
(in the limit) that a copy from type~$c$ is being served in server~$s$, see~\eqref{eq:roschoose}.
The second term multiplied by~$\mu$ in~\eqref{eq:new} represents departures due to a type-$c$ job that is being served in more than one server simultaneously.  	
Together with~\eqref{eq:1}, we obtain
\begin{align}  
&\lim_{r\to\infty}	\tilde f_{s}( r \vec n^r) =  \lambda \frac{d}{K}-\mu
\sum_{c\in \mathcal C(s)}\sum_{l\in c}  \frac{ n_c}{ m_l}, 	\label{eq:driftROS}	
\end{align}

Now, note that~\eqref{eq:driftROS} is equal to~\eqref{eq:fluidIID} (the fluid drift for the PS  model with i.i.d.\ copies). Hence, the proof now follows   as in the proof of Lemma~\ref{lemmaIID}.


\subsection*{E: Proof of Section~\ref{sec:LT}}
\subsubsection*{Light-traffic approximation}

In the light-traffic regime, the number of jobs in the system will be very small. In particular, in our light-traffic approximation, we will assume that at most two jobs will be in the system. Then, the main idea consists in calculating the mean sojourn time  of a tagged job conditioned on its service requirement~$b$,  its type~$c$, and on having at most one other job  present in the system upon its arrival. Unconditioning on the service requirement~$b$, one then obtains the light-traffic approximation for the  unconditional mean sojourn time, denoted by $\bar D^{LT,P}(\lambda)$, where $P$ denotes the scheduling discipline used in the servers.
Using Little's law on the light-traffic approximation, i.e. $\bar{N}^{LT,P}(\lambda) =\lambda \bar{D}^{LT,P}(\lambda)$, one obtains the  result for the mean number of jobs in the system as presented in Lemma~\ref{lem:LT_2}.  


For a given arrival rate~$\lambda>0$, let $\bar{D}^P(\lambda,b)$ denote the mean sojourn time for the tagged  job conditioned on its size being $b$. We let $c$ be the type of the tagged job. Using the ideas as presented in~\cite{Walrand1990519}, we can write   $\bar{D}^P(\lambda,b)=\bar D^{LT,P}(\lambda,b)+o(\lambda^{n+1})$, as $\lambda\to 0$,  where   
\begin{equation}
\label{eq:LTs} \bar D^{LT,P}(\lambda,b) := \bar D^{(0)}(0,b) + \lambda \bar D^{(1)}(0,b)+\ldots+\frac{\lambda^n}{n!} \bar D^{(n)}(0,b),
\end{equation}
is referred to as the light-traffic approximation of order $n$.  
Here, the $i$-th term, $\bar D^{(i)}(0,b)$, denotes the mean sojourn time when in addition to the tagged job, $i$ other jobs arrive to the system in the interval $(-\infty,\infty)$. We note that for the ease of notation we drop the dependency on $P$ of $\bar D^{(n)}(0,b)$.  

We calculate the light-traffic approximation of order $1$, that is, in~\eqref{eq:LTs} we set $n=1$. Hence, we will calculate the sojourn time of the tagged job conditioned on, at most, having one other job present in the system. 
Let $\tilde A(t_0,t_1)$ denote the number of arrivals in the time interval $[t_0,t_1)$ in addition to the tagged job who is assumed to arrive at time $0$.  
The zeroth and first light-traffic derivatives satisfy, see \cite{Walrand1990519}:
$$\bar{D}^{(0)}(0,b):= \mathbb{E}\left(\bar{D}(0,b)\vert \tilde A(-\infty,\infty)=0\right)$$
and 	
\begin{align}
&\bar{D}^{(1)}(0,b):= \int_{-\infty}^{\infty}\left(\mathbb{E}\left(   \bar{D}(0,b)\vert \tilde A(-\infty,\infty)=1,\tau=t\right)\right. \left.- \mathbb{E}\left( \bar{D}(0,b)\vert \tilde A(-\infty,\infty)=0\right)\right)\textrm dt,\nonumber
\end{align} 
where $\tau$ is the arrival time of the other job. 
For any work-conserving policy, it readily follows that $\bar{D}^{(0)}(0,b)= b$, since only the tagged job is present, and all copies of this job are equal to~$b$. 

When in addition to the tagged job, another job is present, the delay of the tagged job will depend on the type of the job already in the system, denoted by $c_1$. If both jobs are of a different type, the new job will start being served immediately,  and hence the first term in the integral of $\overline D^{(1)}(0,b)=b$,  that is,  the first light-traffic derivative is equal to zero.  On the other hand, if both jobs have the same type, which happens with probability $\frac{1}{\binom{K}{d}}$, the job that is already in the system will have an impact  on the sojourn time of the tagged job. We note that the precise value of the impact will depend on $P$, which we quantify later on in the proof of  Lemma~\ref{lem:LT_2}. We thus have:
\begin{equation}
\label{eq:D1}
\bar{D}^{(1)}(0,b):= \frac{1}{\binom{K}{d}} \int_{-\infty}^{\infty}\left(\mathbb{E}\left( \bar{D}(0,b)\vert \tilde A(-\infty,\infty)=1,\tau=t, c_1=c\right)\right. \left.- \mathbb{E}\left( \bar{D}(0,b)\vert \tilde A(-\infty,\infty)=0\right)\right)\textrm dt,
\end{equation}
where we note that in the first term we conditioned on the job being of the same type as the tagged job. 

We note that if the scheduling policy does not depend on $d$, then $\mathbb{E}( \bar{D}(0,b)| \tilde A(-\infty,\infty)=1,\tau=t, c_1=c)$ will not either, hence the light-traffic approximation of order 1 is minimized when $d$ is set equal to  $d^*=\lfloor K/2 \rfloor$.

\subsubsection*{Proof of Lemma~\ref{lem:LT_2}:}	
In order to obtain an expression for $\bar{N}^{LT,P}(\lambda)$, we will calculate $\bar{D}^{LT,P}(\lambda,b)$, uncondition on $b$ and then apply Little's law. 
By \eqref{eq:LTs} and \eqref{eq:D1}, calculating $\bar{D}^{LT,P}(\lambda,b)$ reduces to calculating \\ $\mathbb{E}\left( \bar{D}(0,b)\vert \tilde A(-\infty,\infty)=1,\tau=t,  c_1=c\right)$. 
Below we do so for exponentially distributed service requirements and with $P$ equal to PS,  FCFS, or ROS. 

First consider FCFS.
If in addition to the tagged job, another job arrives in the interval $(-\infty,\infty)$, which is of the same type~$c_1=c$ and has service time $B_1=b_1$, then the sojourn time of the tagged job will be given by
\begin{eqnarray}
& \mathbb{E}\left( \bar{D}(0,b)\vert \tilde A(-\infty,\infty)=1,\tau=t, B_1=b_1, c_1=c\right) 
= \left\{\begin{array}{ll}
b & \textrm{ if } t\leq -b_1,\\

t + b_1 + b & \textrm{ if } -b_1\leq t\leq 0,\\

b & \textrm{ if } t\geq 0.\\
\end{array}\right.\nonumber
\end{eqnarray}
For example, the second equation is the case where the other job arrives before the tagged job, and has still $b_1+t$ remaining service left. Hence, the tagged job has to wait $b_1+t$, so that its sojourn time equals $b+b_1+t$. 
To calculate $\bar{D}^{(1)}(0,b)$, we subtract from the above  $\bar{D}^{(0)}(0,b) = b$ and we multiply with $\frac{1}{\binom{K}{d}}$, integrate over $t$, and uncondition on the service requirements $b_1$. Further unconditioning on  $b$ gives  $\frac{3\lambda}{2\mu^2} \frac{1}{\binom{K}{d}}$. On the other hand, unconditioning $\bar{D}^{(0)}(0,b)$ over $b$ readily yields $1/\mu$. Summing both terms, we get $\bar{D}^{LT,FCFS}(\lambda) = \frac{1}{\mu} + \frac{3\lambda}{2\mu^2} \frac{1}{\binom{K}{d}} $,  and multiplying by $\lambda$ (Little's law) we obtain the expression for $\bar N^{LT,FCFS}(\lambda)$.

The analysis of FCFS carries directly over to ROS,  since at most two jobs are considered to be present in the system, in which case jobs under ROS will be served in order of arrival. 

We now consider PS. We consider the case where in addition to the tagged job, another job arrives in the system in the interval $(-\infty,\infty)$. Let $b_1$ denote the service of this other job and $c_1$ its type.  We have 
\begin{eqnarray}
&\mathbb{E}\left(\bar{D}(0,b)\vert \tilde A(-\infty,\infty)=1,\tau=t, B_1=b_1 c_1=c\right)
= \left\{\begin{array}{ll}
b & \textrm{ if } t\leq -b_1,\\

b+b_1+t & \textrm{ if } -b_1\leq t\leq -b_1+b \textrm{ and } b\leq b_1,\\

2b & \textrm{ if } -b_1+b\leq t\leq 0 \textrm{ and } b\leq b_1,\\

b+b_1+t & \textrm{ if } -b_1\leq t\leq 0 \textrm{ and } b\geq b_1,\\

b+b_1 & \textrm{ if } 0\leq t\leq b-b_1 \textrm{ and } b\geq b_1,\\

2b-t & \textrm{ if } b-b_1\leq t\leq b \textrm{ and } b\geq b_1,\\

2b-t & \textrm{ if } 0\leq t\leq b \textrm{ and } b\leq b_1,\\

b & \textrm{ if } t\geq b.\\
\end{array}\right.\nonumber
\end{eqnarray}
The expression above takes into account all the possible events. For example, the first equation is the case when the other job arrives and leaves before the tagged job arrives. 
The second equation is the case where the other job arrives before the tagged job and leaves first. In that case, the sojourn time experienced by the tagged job is $b$ plus the capacity spend on serving the other job $b_1-(-t)$.  The third equation is the case where the other job arrives before the tagged job and leaves after the tagged job. In that case, the tagged job has shared during its whole stay the server, hence its sojourn time equals $2b$.

To calculate $\bar{D}^{(1)}(0,b)$, we combine all the cases, subtract $\bar{D}^{(0)}(0,b) = b$ and multiply by $\frac{1}{\binom{K}{d}}$,  integrate over $t$ and uncondition over $b_1$. Then, further unconditioning on $b$  we get the expression $\frac{\lambda}{\mu^2} \frac{1}{\binom{K}{d}}$. As in the case of FCFS, summing now with $1/\mu$, we get $\bar{D}^{LT,PS}(\lambda) =\frac{1}{\mu}+\frac{\lambda}{\mu^2} \frac{1}{\binom{K}{d}}$. Multiplying by $\lambda$ yields $\bar N^{LT,PS}(\lambda)$.

\end{appendix}

\end{document}